\documentclass[reqno]{amsart}

\usepackage[T1]{fontenc} 
\usepackage[utf8]{inputenc} 
\usepackage{amssymb}
\usepackage{graphicx} 
\graphicspath{{Figures/}} 
\usepackage{enumitem} 
	\setlist[enumerate]{label={\bf \arabic*.}, ref = {\bf\arabic*}, wide = 0pt} 
	\setlist[enumerate, 2]{label={\bf (\alph*)}, wide = 15pt, leftmargin = 15pt} 
\usepackage{subfig} 
\usepackage{varioref} 
\usepackage[authoryear,sort&compress]{natbib}
\usepackage{tikz}
\usepackage{tikz-cd}
\usetikzlibrary{arrows,arrows.meta}
\usetikzlibrary{calc}
\usepackage{esint} 
\usepackage{amsaddr} 
\theoremstyle{plain} 
\newtheorem{theorem}{Theorem}
\newtheorem{lemma}[theorem]{Lemma}
\newtheorem{prop}[theorem]{Proposition}

\newtheorem{oldtheorem}{Theorem}

\theoremstyle{definition} 
\newtheorem{defn}[theorem]{Definition}
\newtheorem{example}[theorem]{Example}

\theoremstyle{remark} 
\newtheorem{remark}[theorem]{Remark}

\numberwithin{equation}{section}
\numberwithin{theorem}{section}
\usepackage[bookmarks, colorlinks=true, pdfstartview=FitV, linkcolor=blue, citecolor=blue, urlcolor=blue]{hyperref}

\renewcommand{\d}{\mathrm d}

\newcommand{\bdy}{\partial}
\newcommand{\lap}{\Delta}
\newcommand{\Grad}{\nabla}
\newcommand{\weakconv}{\rightharpoonup}

\newcommand{\mc}{\mathcal}
\newcommand{\bb}{\mathbb}
\newcommand{\abs}[1]{\left|#1\right|}

\DeclareMathOperator{\card}{card}

\DeclareMathOperator{\supp}{supp}

\DeclareMathOperator{\dist}{dist}
\DeclareMathOperator{\loc}{loc}
\DeclareMathOperator{\diam}{diam}

\newif\ifdetails
\usepackage[pagewise]{lineno}

\title[Quantization for a nonlocal elliptic equation]{Quantization for sequences of blow-up solutions to an elliptic equation having nonlocal exponential nonlinearity}
\author[M. Gluck]{Mathew Gluck}
\subjclass[2020]{35J15, 35J61, 35J91}
\thanks{This material is based upon work supported by the National Science Foundation under Grant No. DMS-2418889.}
\address{Southern Illinois University \\
School of Mathematical and Statistical Sciences\\ 
Carbondale, IL, U.S.A.}
\email{mathew.gluck@siu.edu}
\keywords{Choquard, concentration compactness, quantization}
\date{\today}
\begin{document}
\begin{abstract}
This work provides a description of the asymptotic behavior of sequences of solutions to an elliptic equation with a nonlocal exponential nonlinearity of Choquard type. The equation under consideration is a nonlocal analog of the classical prescribed Gaussian curvature equation. A concentration-compactness alternative is established for sequences of solutions to the equation under consideration whenever suitable integrability assumptions on the solutions and the curvature functions are satisfied. Under further regularity assumptions on the curvature functions, and when blow-up occurs in the concentration-compactness alternative, an energy quantization result is established. 
\end{abstract}

\maketitle
\section{Introduction}
Elliptic equations with exponential nonlinearities arise in a variety of applications in both pure and applied mathematics. For example, the problem 
\begin{equation}
\label{eq:mean_field}
	\lap_g u + \rho\left(\frac{he^u}{\int_Mhe^u} - 1\right) = 0
	\qquad \text{ on }(M, g),
\end{equation}
where $(M, g)$ is a compact Riemannian surface without boundary, $h$ is a known positive function, and $\rho\in \bb R$ is a parameter arises in mathematical and physical contexts including the prescribed Gaussian curvature problem \cite{KazdanWarner1974} and Chern-Simon Higgs models \cite{Taubes1980, Taubes1980_1, HongKimPac1990, JackiwWeinberg1990, SpruckYang1995, CaffarelliYang1995, Tarantello1996, StruweTarantello1998, DingJostLiWang1997, DingJostLiWang1998}. For a bounded smooth domain $\Omega\subset \bb R^2$, the analogous problem
\begin{equation}
\label{eq:dirichlet_problem}
	\begin{cases}
	\lap u + \rho\frac{he^u}{\int_\Omega he^u} = 0 & \text{ in }\Omega\\
	u=0 & \text{ on }\bdy\Omega, 
	\end{cases}
\end{equation}
where $h$ is a sufficiently smooth positive function on $\Omega$ and $\rho\in \bb R$ is a parameter arises in connection to statistical mechanics of point vortices \cite{Caglioti1992, Caglioti1995, Kiessling1993}. Due in part to the concrete applications of problems \eqref{eq:mean_field} and \eqref{eq:dirichlet_problem}, the question of existence of solutions to these problems has been investigated by many authors using a variety of methods. For example, when $\rho< 8\pi$, existence of solutions to problem \eqref{eq:mean_field} can be routinely established using variational methods, see \cite{Moser1973, KazdanWarner1974, DingJostLiWang1997}. For $\rho\geq 8\pi$ the existence problem for \eqref{eq:mean_field} and \eqref{eq:dirichlet_problem} is much more delicate. To address this issue, a program for the computation of the Leray-Schauder degree for \eqref{eq:mean_field} and \eqref{eq:dirichlet_problem} was initiated in \cite{Li1999}. The program was completed through the combined works of \cite{BrezisMerle1991, LiShafrir1994, Li1999, ChenLin2002, ChenLin2003}. In particular, in \cite{ChenLin2003} it was shown that if $h$ is a suitably smooth positive function on $M$ and if $\rho\in (8\pi m, 8\pi(m + 1))$ for some positive integer $m$, then the Leray-Schauder degree $d_\rho$ for problem \eqref{eq:mean_field} is given by 
\begin{equation*}
	d_\rho= {{m -\chi(M)}\choose{m}},  
\end{equation*}
where ${m\choose k} = \frac{m!}{k!(m - k)!}$ is the binomial coefficient and $\chi(M)$ is the Euler characteristic of $M$. This equality implies that if $M$ is a compact Riemannian surface with genus $\gamma\geq 1$ then $d_\rho\neq 0$, and thus problem \eqref{eq:mean_field} has a solution whenever $\rho$ is not a positive integer multiple of $8\pi$. A similar degree-counting formula was given for the Dirichlet problem \eqref{eq:dirichlet_problem}, see \cite{ChenLin2003} for details.

From an analytical point of view, the primary obstacle in computing the Leray-Schauder degree for problem \eqref{eq:mean_field} (similarly for problem \eqref{eq:dirichlet_problem}) is to obtain sharp pointwise estimates for sequences of blow-up solutions $(V_k, u_k)_{k = 1}^\infty$ to the following local-coordinate model of the problem
\begin{equation}
\label{eq:local_problem}
\begin{cases}
	-\lap u = V e^u & \text{ in }\Omega\\
	\|e^u\|_{L^1(\Omega)} \leq c, 
\end{cases}
\end{equation}
where $\Omega\subset \bb R^2$ is a bounded domain. The purpose of this note is to describe the blow-up mechanism for a nonlocal analog of problem \eqref{eq:local_problem}. The main results of this work are analogous to the results of \cite{BrezisMerle1991, LiShafrir1994} for problem \eqref{eq:local_problem}. To introduce the problems to be considered, let $\mu\in (0, 2)$ and define the convolution operator $I_\mu$ by
\begin{equation*}
	I_\mu f(x) = \int_{\bb R^2}\frac{f(y)}{|x - y|^\mu}\; \d y. 
\end{equation*}
For a domain $\Omega\subset \bb R^2$ we consider problems of the form 
\begin{equation}
\label{eq:nonlocal_version}
\begin{cases}
	-\lap u = VI_\mu[e^{\lambda u}\chi_{\Omega}]e^{\lambda u} & \text{ in }\omega\\
	\|e^u\|_{L^1(\Omega)}\leq c_0,
\end{cases}
\end{equation}
where 
\begin{equation}
\label{eq:lambda}
	\lambda = \frac{4 - \mu}4\in \left(\frac 12, 1\right), 
\end{equation}
$\omega\subset \Omega$ is a subdomain, $V:\omega \to [0, \infty)$, and $\chi_A$ is the characteristic function of a measurable subset $A\subset \bb R^2$. 

The motivation for studying problem \eqref{eq:nonlocal_version} comes from the richness of problem \eqref{eq:local_problem} (and its global versions on manifolds) and the fact that problem \eqref{eq:nonlocal_version} enjoys the same ``limiting symmetries'' as problem \eqref{eq:local_problem}. More specifically, when blow-up occurs in problem \eqref{eq:local_problem}, after a suitable rescaling procedure, one obtains the globally-defined problem 
\begin{equation}
\label{eq:local_entire}
	\begin{cases}
	-\lap u = e^u & \text{ in }\bb R^2\\
	\int_{\bb R^2}e^u< \infty. 
	\end{cases}
\end{equation}
Similarly, as we will show, when blow-up occurs in problem \eqref{eq:nonlocal_version}, after a similar rescaling procedure, one obtains the globally-defined problem
\begin{equation}
\label{eq:entire_PDE_and_L1}
\begin{cases}
	-\lap u = I_\mu[e^{\lambda u}]e^{\lambda u} &\text{ in }\bb R^2\\
	\int_{\bb R^2}e^u <\infty.  
\end{cases}
\end{equation}
The statement that problems \eqref{eq:local_problem} and \eqref{eq:nonlocal_version} enjoy the same limiting symmetries is based on the following observation, which can be directly verified. For any $(x_0, \delta)\in \bb R^2\times (0, \infty)$, both of problems \eqref{eq:local_entire} and \eqref{eq:entire_PDE_and_L1} are invariant under the rescaling $u\mapsto u(\delta(\cdot - x_0)) + 2\log \delta$ and both of these problems are invariant under the logarithmic Kelvin transform $u\mapsto u_{x_0,  \sigma}$, where 
\begin{equation*}
	u_{x_0, \sigma}(x) 
	= u(x^{x_0, \sigma}) +4\log \frac{\sigma}{|x- x_0|}
\end{equation*}
and 
\begin{equation*}
	x^{x_0, \sigma}= x_0 + \frac{\sigma^2(x - x_0)}{|x - x_0|^2}
\end{equation*}
is the inversion of $x\in \bb R^2\setminus\{x_0\}$ through $\bdy B(x_0, \sigma)$. Moreover, the classification result of \cite{ChenLi1991} guarantees that if $u$ is any solution to \eqref{eq:local_entire} then $\|e^u\|_{L^1(\bb R^2)} = 8\pi$ while the classification result of \cite{Gluck2025classification} guarantees that if $u$ is any solution to \eqref{eq:entire_PDE_and_L1} then $\|I_\mu[e^{\lambda u}]e^{\lambda u}\|_{L^1(\bb R^2)} = 8\pi$. For convenience, these classification results are stated in Theorems \ref{theorem:CL_classification} and \ref{oldtheorem:classification} respectively in Appendix \ref{s:appendix}. These observations suggest that sequences of blow-up solutions to problem \eqref{eq:nonlocal_version} may share some similar properties as sequences of blow-up solutions to problem \eqref{eq:local_problem}. We will show that this is indeed the case.

To contextualize the main results of this work we first state the analogous results for problem \eqref{eq:local_problem}. The first such result is the concentration-compactness result in Theorem 3 of \cite{BrezisMerle1991}. In the statement of the theorem and throughout the manuscript $p'= p/(p -1)$ denotes the Lebesgue conjugate of $p$. 
\begin{oldtheorem}
\label{oldtheorem:BM_alternative}
Let $\Omega\subset \bb R^2$ be a bounded domain and let $p\in (1, \infty]$. If $(V_k)_{k = 1}^\infty$ and $(u_k)_{k = 1}^\infty$ be sequences of functions on $\Omega$ for which $V_k\geq 0$ for all $k$, for which there exists a constant $C_0>0$ such that 
\begin{equation*}
	\|V_k\|_{L^p(\Omega)} + \|e^{u_k}\|_{L^{p'}(\Omega)}\leq C_0
	\qquad \text{ for all }k, 
\end{equation*}
and for which 
\begin{equation}
\label{eq:local_PDE}
	-\lap u_k = V_ke^{u_k} \qquad \text{ in }\Omega
\end{equation}
is satisfied for all $k$, then there is a subsequence $(u_{k_\ell})_{\ell = 1}^\infty\subset (u_k)_{k = 1}^\infty$ for which one of the following holds: 
\begin{enumerate}[label = {\bf BM \arabic*.}, ref = {\bf BM \arabic*}, wide = 0pt]
	\item $(u_{k_\ell})_{\ell = 1}^\infty$ is bounded in $L^\infty_{\loc}(\Omega)$, 
	\item $u_{k_\ell}\to-\infty$ locally uniformly on $\Omega$, or
	\item \label{item:BM_finite_blowup} there is a finite nonempty set $S = \{a^1, \ldots, a^m\}\subset \Omega$ such that for each $i\in \{1, \ldots, m\}$ there is a sequence $(x_\ell^i)_{\ell = 1}^\infty$ such that $\lim_{\ell\to\infty}x_\ell^i= a^i$ and $u_{k_\ell}(x_\ell^i)\to\infty$. Moreover, $u_{k_\ell}\to-\infty$ locally uniformly on $\Omega \setminus S$ and 
	\begin{equation}
	\label{eq:BM_to_point_masses}
		V_{k_\ell}e^{u_{k_\ell}}\weakconv \sum_{i = 1}^m\alpha_i\delta_{a^i}
	\end{equation}
	weak-$*$ in the sense of measures with $\alpha_i\geq 4\pi$ for all $i$. 
\end{enumerate}
\end{oldtheorem}
The following quantization result due to \cite{LiShafrir1994} computes the values of $\alpha_1, \ldots, \alpha_m$ in alternative \ref{item:BM_finite_blowup} of Theorem \ref{oldtheorem:BM_alternative}. 
\begin{oldtheorem}
\label{oldtheorem:li_shafrir}
Let $\Omega\subset \bb R^2$ be a bounded domain and suppose $(V_k)_{k = 1}^\infty\subset C^0(\overline\Omega)$ is a sequence of nonnegative functions for which $V_k\to V$ for some $V\in C^0(\overline\Omega)$. If $(u_k)_{k = 1}^\infty$ is a sequence of solutions to \eqref{eq:local_PDE} for which there exists $C_0>0$ such that $\|e^{u_k}\|_{L^1(\Omega)}\leq C_0$ for all $k$ and if alternative \ref{item:BM_finite_blowup} holds in Theorem \ref{oldtheorem:BM_alternative} (with $p = +\infty$) then for every $i\in \{1, \ldots, m\}$ there is $N_i\in \bb N$ for which $\alpha_i = 8\pi N_i$. 
\end{oldtheorem}
The first result of this work is the following analog of Theorem \ref{oldtheorem:BM_alternative} for problem \eqref{eq:nonlocal_version}. Although the primary case of interest for Theorem \ref{theorem:BM_alternative} is the case $p = \infty$, we will state and prove this theorem for a broader range of $p$. 
\begin{theorem}
\label{theorem:BM_alternative}
Let $\omega\subset \bb R^2$ be a bounded domain, let $\mu\in (0, 2)$, let $p\in (\frac 2\mu, \infty]$, and let $q$ be given by 
\begin{equation}
\label{eq:integrability_power}
	\frac 1 q + \frac 1{2p}= \lambda, 
\end{equation}
where $\lambda$ is as in \eqref{eq:lambda}. Suppose $(\Omega_k)_{k = 1}^\infty$ is a sequence of bounded domains in $\bb R^2$ for which 
\begin{equation}
\label{eq:nested_domains}
	\omega\subset \Omega_1\subset \Omega_2\subset \ldots. 
\end{equation} 
If $u_k:\Omega_k\to \bb R$ and $V_k:\omega\to [0, \infty)$ are functions for which 
\begin{equation}
\label{eq:sequence_of_problems}
	-\lap u_k = V_kI_\mu[e^{\lambda u_k}\chi_{\Omega_k}]e^{\lambda u_k}
	\qquad \text{ in }\omega
\end{equation}
and for which there is $c_0>0$ such that
\begin{equation}
\label{eq:vkp_uklambdaq_energy}
	\|V_k\|_{L^p(\omega)} + \|e^{\lambda u_k}\|_{L^q(\Omega_k)}
	\leq c_0
	\qquad \text{ for all }k, 
\end{equation} 
then there is a subsequence $(u_{k_\ell})_{\ell = 1}^\infty\subset(u_k)_{k = 1}^\infty$ for which one of the following alternatives holds:
\begin{enumerate}[label = {\bf A\arabic*.}, ref = {\bf A\arabic*}, wide = 0pt]
	\item \label{item:uniformly_bounded} $(u_{k_\ell})_{\ell = 1}^\infty$ is bounded in $L^\infty_{\loc}(\omega)$, 
	\item  \label{item:uniformly_to_-infty} $u_{k_\ell}\to-\infty$ locally uniformly on $\omega$, or
	\item \label{item:finite_blow_up} there exists a finite nonempty set $S = \{a^1, \ldots, a^m\}\subset \omega$ such that, for all $i\in \{1, \ldots, m\}$ there is a sequence $(x_\ell^i)_{\ell = 1}^\infty\subset\omega$ with $x_\ell^i\to a^i$ and $u_{k_\ell}(x_\ell^i)\to \infty$. Moreover, $u_{k_\ell}\to-\infty$ locally uniformly on $\omega \setminus S$ and 
	\begin{equation}
	\label{eq:point_masses}
		V_{k_\ell}I_\mu[e^{\lambda u_{k_\ell}}\chi_{\Omega_{k_\ell}}]e^{\lambda u_{k_\ell}}\weakconv \sum_{i = 1}^m \alpha_i\delta_{a^i}
	\end{equation}
	weak-$*$ in the sense of measures on $\omega$ with $\alpha_i\geq  4\pi\left(1 - \frac{1}{2\lambda p}\right)$ for all $i$. 
\end{enumerate}
\end{theorem}
The second result of this work is an inequality of $\sup + \inf$ type in the spirit of \cite{Shafrir1992, BrezisLiShafrir1993}. See also \cite{EspositoLucia2021} for a $\sup +\inf$ inequality for the nonlinear $n$-Laplacian. To state the result, for $a\leq b$ and $S\subset \bb R^2$ we introduce the notation 
\begin{equation}
\label{eq:lambda_ab}
	\Lambda_{a,b}(S) = \{V\in C^0(S): a\leq V\leq b \text{ in }S\}. 
\end{equation}
\begin{theorem}
\label{theorem:sup_inf_inequality}
Let $\omega\subset \bb R^2$ be a bounded domain, let $\mu\in (0, 2)$ and let $\lambda$ be as in \eqref{eq:lambda}. Suppose $0< a\leq b< \infty$ and $\Lambda\subset\Lambda_{a,b}(\omega)$ is a subset that is equicontinuous at each point of $\omega$. For each compact subset $K\subset \omega$, each $c_0>0$, and each $C_1>1$ there is a constant $C_2= C_2(\Lambda, K, \omega, c_0, C_1)>0$ such that if $V\in \Lambda$, if $\Omega\subset \bb R^2$ is a bounded domain for which $\omega\subset\Omega$ and if $u$ is a distributional solution to \eqref{eq:nonlocal_version} then 
\begin{equation}
\label{eq:sup+inf_inequality}
	\max_Ku + C_1\inf_{\omega} u\leq C_2. 
\end{equation}
We emphasize that $C_2$ is independent of $\Omega$. 
\end{theorem}

Our final result is an analog of Theorem \ref{oldtheorem:li_shafrir}. It concerns the specification of the coefficients $\alpha_i$ in \eqref{eq:point_masses} when $p = \infty$ and alternative \ref{item:finite_blow_up} is assumed to hold. 
\begin{theorem}
\label{theorem:quantization}
Let $\Omega\subset \bb R^2$ be a bounded domain, let $\omega\subset \Omega$ be a subdomain, let $\mu\in (0, 2)$ and let $\lambda$ be as in \eqref{eq:lambda}. Suppose $V\in C^0(\overline\omega)$ and $(V_k)_{k = 1}^\infty$ is a sequence of nonnegative functions on $\overline \omega$ for which $V_k\to V$ in $C^0(\overline \omega)$. If $(u_k)_{k = 1}^\infty$ is a sequence of solutions to 
\begin{equation*}
	-\lap u_k = V_k I_\mu[e^{\lambda u_k}\chi_\Omega]e^{\lambda u_k} 
	\qquad \text{ in }\omega
\end{equation*}
for which there exists $c_0>0$ satisfying $\|e^{u_k}\|_{L^1(\Omega)}\leq c_0$ for all $k$ and if alternative \ref{item:finite_blow_up} holds in Theorem \ref{theorem:BM_alternative} with $p = +\infty$, then for each $i\in \{1, \ldots, m\}$ there is $N_i\in \bb N$ for which $\alpha_i = 8\pi N_i$. 
\end{theorem}
Our approach in proving Theorems \ref{theorem:BM_alternative} and \ref{theorem:quantization} mirrors the approaches of \cite{BrezisMerle1991} and \cite{LiShafrir1994} in proving Theorems \ref{oldtheorem:BM_alternative} and \ref{oldtheorem:li_shafrir} respectively. Techniques used in the proofs of Theorems \ref{oldtheorem:BM_alternative} and \ref{oldtheorem:li_shafrir} provide a general framework for the proofs of Theorems \ref{theorem:BM_alternative} and \ref{theorem:quantization}, but the presence of the nonlocal term in \eqref{eq:nonlocal_version} makes the application of these techniques more involved and new ideas must be introduced. For example, in the proof of Theorem \ref{theorem:quantization}, we employ a ``bubble selection process'' in the spirit of that employed by \cite{LiShafrir1994} in the proof of Theorem \ref{oldtheorem:li_shafrir}. However, compared to the proof of Theorem \ref{oldtheorem:li_shafrir}, we face the additional task of ruling out the possibility of nonlocal interactions (at the $L^1$ level) between distinct bubbles.

The paper is organized as follows. In Section \ref{s:preliminaries} we discuss some preliminary notions including the definition of a distributional solution and the basic regularity results for distributional solutions to the problem under consideration. The proof of Theorem \ref{theorem:BM_alternative} is provided in Section \ref{s:concentration_compactness}. The proof of Theorem \ref{theorem:sup_inf_inequality} is provided in Section \ref{s:sup_inf_inequality}. Section \ref{s:quantization} is devoted to the proof of Theorem \ref{theorem:quantization}. Finally, Section \ref{s:appendix} is an appendix where some computations and some statements of prior results are kept. 

The symbol $C$ will be used throughout the manuscript to denote various $k$-independent constants. We adopt the custom that the value of $C$ may change from line to line and even within the same line. 

\section{Preliminaries}
\label{s:preliminaries}
\begin{defn}
\label{defn:distributional_solution}
Let $\omega\subset\bb R^2$ be an open set (possibly unbounded) and let $f\in L^1(\omega)$. A \emph{distributional solution} to $-\lap u = f$ in $\omega$ is a function $u\in L^1_{\loc}(\omega)$ for which 
\begin{equation*}
	-\int_{\omega}u\lap \varphi = \int_{\omega}f\varphi
	\qquad\text{ for all }\varphi\in C_c^\infty(\omega). 
\end{equation*}
\end{defn}
We will make extensive use of the Hardy-Littlewood-Sobolev inequality which we recall here. For brevity, throughout the manuscript we refer to this inequality as the \emph{HLS inequality}. 
\begin{oldtheorem}
\label{theorem:HLS}
Let $n\geq 1$, let $\mu\in (0, n)$ and suppose $p, r\in (1,  \infty)$ satisfy $\frac 1 r = \frac 1 p - \frac{n- \mu}n$. There is an optimal constant $\mc H = \mc H(n, \mu, p)>0$ such that for all $f\in L^p(\bb R^n)$, 
\begin{equation*}
	\|I_\mu f\|_{L^r(\bb R^n)}
	\leq \mc H \|f\|_{L^p(\bb R^n)}. 
\end{equation*}	
\end{oldtheorem}
The following lemma justifies (in part) the integrability assumption on $e^{u_k}$ in Theorem \ref{theorem:BM_alternative}. We refer the reader to Lemma 2.2 of \cite{Gluck2025classification} for a proof. 
\begin{lemma}
\label{lemma:integrability_assumption}
Let $\Omega\subset \bb R^2$ be a (possibly unbounded) domain, let $\mu\in (0, 2)$ and let $\lambda$ be as in \eqref{eq:lambda}. If $p\in (\frac2\mu, \infty]$ and if $e^{\lambda u}\in L^q(\Omega)$, where $q$ is defined by \eqref{eq:integrability_power},
then $I_\mu[e^{\lambda u}\chi_\Omega]e^{\lambda u}\in L^{p'}(\Omega)$ and there is a constant $C = C(p, \mu)>0$ such that
\begin{equation*}
	\|I_\mu[e^{\lambda u}\chi_\Omega]e^{\lambda u}\|_{L^{p'}(\Omega)}
	\leq C\|e^{\lambda u}\|_{L^q(\Omega)}^2. 
\end{equation*}
In particular, under these hypotheses, if $\omega\subset \Omega$ is any subdomain, if $V\in L^p(\omega)$ and if $u\in L^1_{\loc}(\omega)$, then the problem 
\begin{equation}
\label{eq:subdomain_PDE}
	-\lap u  = VI_\mu[e^{\lambda u}\chi_\Omega]e^{\lambda u}
	\qquad \text{ in }\omega
\end{equation}
has a meaning in the sense of distributions. 
\end{lemma}
\subsection{Regularity for a single solution}
In this subsection we show that under suitable integrability assumptions on $V$ and $e^u$ the distributional solutions to $-\lap u = VI_\mu[e^{\lambda u}\chi_\Omega]e^{\lambda u}$ are in $C^{1, \alpha}$ for some $\alpha\in (0, 1)$. The following proposition is the main result of the subsection.
\begin{prop}
\label{prop:C1alpha_regularity}
Let $\Omega\subset \bb R^2$ be a bounded domain and let $\omega\subset \Omega$ be a subdomain. Let $\mu\in (0, 2)$, let $\lambda$ be as in \eqref{eq:lambda}, let $p\in \left(\max\{\frac 2\mu, \frac 2{2 - \mu}\}, \infty\right]$, define $q$ by \eqref{eq:integrability_power} and suppose $0\leq V\in L^p(\omega)$. If $u\in L^1_{\loc}(\omega)$ is a distributional solution to \eqref{eq:subdomain_PDE} for which $e^{\lambda u}\in L^q(\Omega)$ then for any $t$ that satisfies
\begin{equation}
\label{eq:admissible_large_t}
	\frac 1{2p} + \frac\mu4 < \frac 1t < \frac 12
\end{equation}
we have $u\in W^{2, t}_{\loc}(\omega)\subset C^{1, 1 - \frac 2t}(\omega)$. 
\end{prop}
The remainder of this subsection is devoted to the proof of Proposition \ref{prop:C1alpha_regularity}. Under the hypotheses of Lemma \ref{lemma:integrability_assumption}, the integrability of the function 
\begin{equation*}
	f = f(u) = VI_\mu[e^{\lambda u}\chi_\Omega]e^{\lambda u} 
\end{equation*}
that appears on the right-hand side of \eqref{eq:subdomain_PDE} is ensured. The proof of Proposition \ref{prop:C1alpha_regularity} relies on an improved integrability result for $f$, which we establish below in Lemma \ref{lemma:boost_integrability}. Before stating Lemma \ref{lemma:boost_integrability}, we first recall some basic properties of the convolution operator determined by the fundamental solution
\begin{equation}
\label{eq:fundamental_solution}
	\Gamma(y) = -\frac 1{2\pi}\log|y|
\end{equation} 
for $-\lap$ on $\bb R^2$. Since the first item in the following lemma is well-known and the remaining items can be established by making obvious modifications to the proof of Lemma 2.2 of \cite{Gluck2020classification}, we omit the proof. 
\begin{lemma}
\label{lemma:Gamma_mapping_properties}
Let $f\in L^1(\bb R^2)$ with $\supp f\subset B_r$ for some $r\geq 1$ and let $\Gamma$ be as in \eqref{eq:fundamental_solution}. The following hold: 
\begin{enumerate}
	\item $\Gamma*f\in L^1_{\loc}(\bb R^2)$ is a distributional solution to $-\lap u = f$ in $\bb R^2$ in the sense of Definition \ref{defn:distributional_solution}. 
	\item $\Gamma*f\in W^{1, 1}_{\loc}(\bb R^2)$ and for every $i\in \{1, 2\}$, 
	\begin{equation*}
		\partial_i(\Gamma*f)(x)
		= -\frac 1{2\pi}\int_{\supp f}\frac{x_i - y_i}{|x- y|^2}f(y)\; \d y, 
	\end{equation*}
	where equality holds in the sense of $L^1_{\loc}(\bb R^2)$. 
	\item If, in addition to the above hypotheses,  $f\in L^p(\bb R^2)$ for some $p> 2$ then $\Gamma*f\in W^{1, \infty}_{\loc}(\bb R^2)$. 
	\item If, in addition to the above hypotheses, $f\in L^\infty(\bb R^2)$ then $\Gamma*f\in C^1(\bb R^2)$. 
\end{enumerate}
\end{lemma}
The proof of the following basic inequality can be found in Theorem 1 of \cite{BrezisMerle1991}. 
\begin{lemma}
\label{lemma:BM_basic_inequality}
Let $\Omega\subset \bb R^2$ be a bounded domain, let $f\in L^1(\Omega)$ and suppose $u$ satisfies
\begin{equation*}
\begin{cases}
	-\lap u = f & \text{ in }\Omega\\
	u = 0 & \text{ on }\bdy \Omega. 
\end{cases}
\end{equation*}
For every $\delta\in (0, 4\pi)$ there holds
\begin{equation*}
	\int_\Omega \exp\left[\frac{(4\pi - \delta)|u(x)|}{\|f\|_{L^1(\Omega)}}\right]\; \d x
	\leq \frac{4\pi^2}\delta(\diam\Omega)^2. 
\end{equation*}
\end{lemma}
\begin{lemma}
\label{lemma:boost_integrability}
Let $\Omega\subset \bb R^2$ be a bounded domain and let $\omega\subset \Omega$ be a subdomain. Let $\mu\in (0, 2)$, let $p\in (\frac 2\mu, \infty]$, define $q$ by \eqref{eq:integrability_power} and suppose $0\leq V\in L^p(\omega)$. If $u\in L^1_{\loc}(\omega)$ is a distributional solution to \eqref{eq:subdomain_PDE} for which $e^{\lambda u}\in L^q(\Omega)$ then for every $t$ satisfying 
\begin{equation}
\label{eq:t_range}
	\frac 1{2p} + \frac\mu 4< \frac 1t< 1
\end{equation}
we have $VI_\mu[e^{\lambda u}\chi_\Omega]e^{\lambda u}\in L^t_{\loc}(\omega)$. 
\end{lemma}
\begin{proof}
Lemma \ref{lemma:integrability_assumption} guarantees that $I_\mu[e^{\lambda u}\chi_\Omega]e^{\lambda u}\in L^{p'}(\omega)$, so in view of the assumption $V\in L^p(\omega)$, H\"older's inequality guarantees that the function
\begin{equation}
\label{eq:the_nonlinearity}
	f:= VI_\mu[e^{\lambda u}\chi_\Omega]e^{\lambda u} 
\end{equation}
satisfies $f\in L^1(\omega)$. With $\Gamma$ as in \eqref{eq:fundamental_solution}, from Lemma \ref{lemma:Gamma_mapping_properties} we have $\Gamma*f\in W^{1, 1}_{\loc}(\bb R^2)\subset W^{1, 1}(\omega)$ and $-\lap(\Gamma*f) = f$ in the distributional sense on $\omega$. 
Since $\varphi:= u - \Gamma*f$ satisfies $-\lap\varphi = 0$ in the distributional sense on $\omega$, Weyl's Lemma guarantees that $\varphi\in C^\infty(\omega)\subset W^{1, 1}_{\loc}(\omega)$. We deduce that $u= \Gamma*f + \varphi\in W^{1, 1}_{\loc}(\omega)$. Next, let $t$ satisfy \eqref{eq:t_range}, let $\epsilon>0$ be sufficiently small so that 
\begin{equation*}
	\frac{1}{2p} + \frac \mu 4 + \lambda \epsilon < \frac 1t, 
\end{equation*}
and decompose $f$ as
\begin{equation*}
	f = f\chi_{\{f> M\}} + f\chi_{\{0\leq f\leq M\}}
	=: f_1 + f_2, 
\end{equation*}
where $M$ is chosen sufficiently large so that $\|f_1\|_{L^1(\omega)}\leq \epsilon$. By construction we have $\|f_2\|_{L^\infty(\omega)}\leq M$. 
%
In what follows we distinguish constants that depend on $M$ (and hence also on the distribution function of $u$) from those that do not by writing $C_M$ and $C$ respectively for any such constants. Let $u_1\in W_0^{1, 1}(\omega)$ be the weak solution to 
\begin{equation*}
\begin{cases}
	-\lap u_1 = f_1 & \text{ in }\omega\\
	u_1 = 0 & \text{ on }\bdy \omega
\end{cases}
\end{equation*}
and apply Lemma \ref{lemma:BM_basic_inequality} to $u_1$ with $\delta= 4\pi - 1$ to obtain 
\begin{equation}
\label{eq:fromBM_basic_inequality}
	\int_{\omega}\exp\left(\frac{|u_1|}\epsilon\right)
	\leq \int_{\omega}\exp\left(\frac{|u_1|}{\|f_1\|_{L^1(\omega)}}\right)
	\leq \frac{4\pi^2}{4\pi - 1}(\diam \omega)^2. 
\end{equation}
This estimate also implies 
\begin{equation*}
	\int_\omega|u_1|
	\leq \epsilon\int_{\omega}\exp\left(\frac{|u_1|}\epsilon\right)
	\leq C(\omega).
\end{equation*}
Let $u_2\in W_0^{1, 2}(\omega)$ be the weak solution to 
\begin{equation*}
\begin{cases}
	-\lap u_2 = f_2 & \text{ in }\omega\\
	u_2 = 0 & \text{ on }\bdy \omega. 
\end{cases}
\end{equation*}
Since $f_2\in L^\infty(\omega)$ standard elliptic estimates give
\begin{equation}
\label{eq:u2_uniform_bound}
	\|u_2\|_{L^\infty(\omega)}
	\leq C(\omega, \|f_2\|_{L^\infty(\bb R^2)})
	= C_M(\omega). 
\end{equation}
Defining $u_3 = u - u_1 - u_2$, we have $-\lap u_3 =0$ in the distributional sense on $\omega$. 
For any $x\in \omega$ and any $R>0$ for which $B_R(x)\subset\omega$ the Mean Value Theorem gives
\begin{equation*}
\begin{split}
	|B_R|u_3(x)
	& = \int_{B_R(x)}u_3(y)\; \d y\\
	& \leq\int_{\omega}u_3^+(y)\; \d y\\
	& \leq \int_{\omega}(u^+(y) + |u_1(y)| + |u_2(y)|)\; \d y\\
	& \leq \int_{\omega}e^{\lambda q u(y)}\; \d y + \|u_1\|_{L^1(\omega)} + |\omega|\|u_2\|_{L^\infty(\omega)}\\
	& \leq \|e^{\lambda u}\|_{L^q(\omega)}^q + C(\omega) + C_M(\omega).
\end{split}
\end{equation*}
This estimate implies that $u_3^+\in L^\infty_{\loc}(\omega)$ and for any compact subset $K\subset \omega$ there is a constant $C_M = C_M(\mu, p, K, \omega, \|e^{\lambda u}\|_{L^q(\omega)})>0$ for which 
\begin{equation}
\label{eq:u3_upper_bound}
	\|u_3^+\|_{L^\infty(K)}\leq C_M. 
\end{equation}
Fixing such a set $K$, for any $x\in K$, applying \eqref{eq:u2_uniform_bound} and \eqref{eq:u3_upper_bound} we have
\begin{equation}
\label{eq:elambdau_upper_bound}
	e^{u(x)}
	\leq e^{(|u_2(x)| + u_3^+(x))}e^{|u_1(x)|}
	\leq C_Me^{|u_1(x)|}
\end{equation}
for some constant $C_M = C_M(\mu, p, K, \omega, \|e^{\lambda u}\|_{L^q(\omega)})$. Now define $r$ by 
\begin{equation*}
	\frac 1r = \frac 1q - \frac{2 - \mu}2 
\end{equation*}
and note that from the definition of $q$ in \eqref{eq:integrability_power} there holds $\frac 1p + \frac 1r = \frac 1{2p} + \frac\mu 4$. Defining $\sigma = \frac 1t - (\frac 1p + \frac 1r + \lambda \epsilon)$ and using each of  \eqref{eq:elambdau_upper_bound}, H\"older's inequality, the HLS inequality, and \eqref{eq:fromBM_basic_inequality} we have 
\begin{equation*}
\begin{split}
	\|f\|_{L^t(K)}
	& \leq C_M\left[\int_KV^tI_\mu[e^{\lambda u}\chi_\Omega]^t e^{\lambda t|u_1|}\; \d x\right]^{1/t}\\
	&\leq C_M|K|^\sigma\|V\|_{L^p(\omega)}\|I_\mu[e^{\lambda u}\chi_\Omega]\|_{L^{r}(\Omega)}\left(\int_Ke^{\frac{|u_1|}\epsilon}\right)^{\lambda\epsilon}\\
	& \leq C_M
	\|V\|_{L^p(\omega)}\|e^{\lambda u}\|_{L^q(\Omega)}\\
	&< \infty. 
\end{split}
\end{equation*} 
\end{proof}
With Lemma \ref{lemma:boost_integrability} in hand we are ready to give the proof of Proposition \ref{prop:C1alpha_regularity}. 
\begin{proof}[Proof of Proposition \ref{prop:C1alpha_regularity}]
Fix any $t$ that satisfies \eqref{eq:admissible_large_t} and apply Lemma \ref{lemma:boost_integrability} to find that $f\in L^t_{\loc}(\omega)$, where $f$ is as in \eqref{eq:the_nonlinearity}. By standard estimates on the Newtonian potential we have $\Gamma*f\in W^{2, t}_{\loc}(\omega)$. The function $\varphi = u - \Gamma*f$ satisfies $-\lap \varphi = 0$ in the distributional sense on $\omega$ so Weyl's Lemma guarantees that $\varphi\in C^\infty(\omega)\subset W^{2, t}_{\loc}(\omega)$. Therefore $u= \Gamma * f + \varphi \in W^{2, t}_{\loc}(\omega)\subset C^{1, 1 - \frac 2t}(\omega)$. 
\end{proof}
\section{Proof of the Concentration-Compactness Alternative}
\label{s:concentration_compactness}
This section is devoted to the proof of Theorem \ref{theorem:BM_alternative}. Let $\mu$, $p$ and $q$ satisfy the hypotheses of Proposition \ref{prop:C1alpha_regularity}. We assume in this subsection that $\omega$ and $(\Omega_k)_{k = 1}^\infty$ are bounded open subsets of $\bb R^2$ for which \eqref{eq:nested_domains} holds. For each $k\in \bb N$ we consider functions $0\leq V_k\in L^p(\omega)$ and $u_k\in L^1_{\loc}(\omega)$ that satisfy \eqref{eq:sequence_of_problems} in the distributional sense.  We also assume the existence of $c_0>0$ for which \eqref{eq:vkp_uklambdaq_energy} holds. Under these assumptions Proposition \ref{prop:C1alpha_regularity} guarantees the existence of $\alpha\in (0, 1)$ for which the containment $u_k\in C^{1, \alpha}(\omega)$ holds for all $k$. 
\begin{defn}
\label{defn:blow_up_point}
Let $\omega\subset \bb R^2$ and let $(u_k)_{k = 1}^\infty$ be a sequence of real-valued functions on $\omega$. A \emph{blow-up point} for $(u_k)_{k = 1}^\infty$ is a point $x_0\in \overline\omega$ for which there is a sequence $(x_k)_{k = 1}^\infty\subset\omega$ such that both $x_k\to x_0$ and $(u_k(x_k))_{k = 1}^\infty$ is unbounded from above. The \emph{blow-up set} for $(u_k)_{k = 1}^\infty$ is the set of blow-up points. 
\end{defn}

The following lemma gives an imprecise $L^1$-energy threshold on $(-\lap u_k)_{k = 1}^\infty$ under which blowup cannot occur for solutions to \eqref{eq:sequence_of_problems}. It is an analog of Corollary 4 of \cite{BrezisMerle1991}. In the lemma and throughout this section we use the notation 
\begin{equation}
\label{eq:fk}
	f_k = V_kI_\mu[e^{\lambda u_k}\chi_{\Omega_k}]e^{\lambda u_k}. 
\end{equation}
\begin{lemma}
\label{lemma:uniform_upper_bound}
Let $\omega\subset \bb R^2$ be a bounded domain, let $\mu\in (0, 2)$, let $p\in (\frac 2\mu, \infty]$, and let $q$ be given by \eqref{eq:integrability_power}, where $\lambda$ is as in \eqref{eq:lambda}.  Suppose $(\Omega_k)_{k = 1}^\infty$ is a sequence of bounded domains in $\bb R^2$ for which \eqref{eq:nested_domains} holds. For each $k$, suppose $V_k:\omega\to [0, \infty)$ and $u_k:\Omega_k\to \bb R$ are functions for which $u_k\in L^1_{\loc}(\omega)$ and for which \eqref{eq:sequence_of_problems} is satisfied in the distributional sense. Suppose further that there exists a constant $c_0>0$ such that \eqref{eq:vkp_uklambdaq_energy} is satisfied. If there is $\beta>0$ for which 
\begin{equation}
\label{eq:uniform_energy_bound}
	\|f_k\|_{L^1(\omega)}
	\leq \beta
	< 4\pi\left(1 - \frac 1{2\lambda p}\right)
	\qquad \text{ for all }k, 
\end{equation}
where $f_k$ is as in \eqref{eq:fk}, then $(u_k^+)_{k = 1}^\infty$ is bounded in $L^\infty_{\loc}(\omega)$. 
\end{lemma}
%
\begin{proof}[Proof
	]
It suffices to assume $\omega = B_R$ and to prove that $(u_k^+)_{k = 1}^\infty$ is bounded in $L^\infty(B_{R/4})$. Assumption \eqref{eq:vkp_uklambdaq_energy} guarantees that $(\|u_k^+\|_{L^1(\Omega_k)})_{k = 1}^\infty$ is bounded in $\bb R$. 
Choose $\delta\in (0, 4\pi)$ for which 
\begin{equation}
\label{eq:choice_of_delta}
	\frac{\beta}{4\pi - \delta}
	< 1 - \frac 1{2\lambda p}
\end{equation}
and define
\begin{equation}
\label{eq:s_defn}
	s = \frac{4\pi - \delta}{\lambda \beta},  
\end{equation}
so that $1< q< s$. 
In what follows we use $C$ (or $C(a)$ for some parameters $a$) to denote various $k$-independent positive constants that may depend on $\mu$, $p$, $c_0$, $R$ and $\delta$ in addition to any explicitly indicated parameters. Decompose $u_k$ as $u_k = v_k + w_k$ where $v_k\in W_0^{1, 1}(B_R)$ is the weak solution to
\begin{equation}
\label{eq:vk_BVP}
\begin{cases}
	-\lap v_k = f_k & \text{ in }B_R\\
	v_k= 0 & \text{ on }\bdy B_R 
\end{cases}
\end{equation}
and $w_k$ is harmonic in $B_R$. Applying Lemma \ref{lemma:BM_basic_inequality} to $v_k$ and using assumption \eqref{eq:uniform_energy_bound} gives
\begin{equation*}
\begin{split}
	\frac{16\pi^2R^2}\delta
	& \geq \int_{B_R}\exp\left(\frac{(4\pi - \delta)|v_k|}{\|f_k\|_{L^1(B_R)}}\right)\; \d x\\
	& \geq \int_{B_R}\exp\left(\frac{(4\pi - \delta)|v_k|}{\beta}\right)\; \d x\\
	& = \int_{B_R}\exp\left(\lambda s|v_k|\right)\; \d x, 
\end{split}
\end{equation*}
and thus $(e^{\lambda|v_k|})_{k = 1}^\infty$ is bounded in $L^{s}(B_R)$. This bound 
implies that $(v_k)_{k = 1}^\infty$ is bounded in $L^s(B_R)$. 
 Since $w_k$ is harmonic in $B_R$, for each $x\in B_{R}$ and any $\rho\in (0, R - |x|)$, the Mean Value Theorem gives
\begin{equation*}
\label{eq:wk_pointwise_estimate}
\begin{split}
	|B_{\rho}|w_k(x)
	& = \int_{B_{\rho}(x)}w_k\\
	& \leq \int_{B_R}w_k^+\\
	& 
	\leq \int_{B_R}(u_k^+ + |v_k|)\\
	& \leq C.
\end{split}
\end{equation*}
Using this estimate, it is routine to show that $(w_k^+)_{k= 1}^\infty$ is bounded in $L^\infty_{\loc}(B_R)$ and, for any compact subset $K\subset B_R$ we have $\|w_k^+\|_{L^\infty(K)}\leq C(\dist(K, \bdy\Omega))$. 
Since, in addition, $u_k\leq |v_k| + w_k^+$ for all $k$, we find that $(e^{\lambda u_k})_{k = 1}^\infty$ is bounded in $L^{s}_{\loc}(B_R)$. In particular, 
\begin{equation}
\label{eq:elambdauk_Ls_bound}
	\|e^{\lambda u_k}\|_{L^s(B_{R/2})}
	\leq C. 
\end{equation}
Next, define $r\in (1, \infty)$ by $\frac 1r = \frac 1q - \frac{2 - \mu}2$ and define $t\in (1, \infty)$ by 
\begin{equation*}
	\frac 1t
	= \frac 1p + \frac 1r + \frac 1s
	= \frac 1{2p} + \frac\mu4+ \frac 1s. 
\end{equation*}
From H\"older's inequality, the HLS inequality, assumption \eqref{eq:vkp_uklambdaq_energy}, and estimate \eqref{eq:elambdauk_Ls_bound} we have
\begin{equation*}
\begin{split}
	\|f_k\|_{L^t(B_{R/2})}
	& \leq \|V_k\|_{L^p(B_R)}\|I_\mu[e^{\lambda u_k}\chi_{\Omega_k}]\|_{L^r(B_{R/2})}\|e^{\lambda u_k}\|_{L^s(B_{R/2})}\\
	& \leq C\|V_k\|_{L^p(B_R)}\|e^{\lambda u_k}\|_{L^q(\Omega_k)}\|e^{\lambda u_k}\|_{L^s(B_{R/2})}\\
	& \leq C. 
\end{split}
\end{equation*}
Since $(f_k)_{k = 1}^\infty$ is bounded in $L^t(B_{R/2})$ with $t > 1$ and since $v_k$ satisfies \eqref{eq:vk_BVP}, standard elliptic estimates imply that $(v_k)_{k = 1}^\infty$ is bounded in $L^\infty(B_{R/4})$. 
Finally, in view of the pointwise inequality
$u_k^+
\leq |v_k| + w_k^+$ and in view of the fact that $(w_k^+)_{k = 1}^\infty$ is bounded in $L^\infty_{\loc}(B_R)$ we obtain $\|u_k^+\|_{L^\infty(B_{R/4})}\leq C$. 
\end{proof}

\begin{example}
Let $\Omega = B_2\subset \bb R^2$ and let $\omega = B(e_1, \frac 34)$, where $e_1 = (1, 0)\in \bb R^2$. Proposition \ref{prop:C1alpha_regularity} applied with $p = \infty$ and $q = 1/\lambda$ guarantees that if $0\leq V\in L^\infty(B(e_1, \frac 34))$ and if $u$ is a distributional solution to 
\begin{equation*}
	-\lap u = VI_\mu[e^{\lambda u}\chi_{B_2}]e^{\lambda u}
	\qquad \text{ in }B(e_1, \frac 34)
\end{equation*}
for which $e^u\in L^1(B_2)$, then $u\in C^{1, \alpha}(\overline B(e_1, \frac 12))$. This example shows that we cannot estimate either of $u$ or the function $\varphi = u - \Gamma* f$ from below in terms of $\|V\|_{L^\infty}$ and $\|e^{\lambda u}\|_{L^q(B_2)} = \|e^u\|_{L^1(B_2)}^\lambda$, even if an arbitrarily small positive bound on $\|f\|_{L^1(\omega)}$ is imposed. Here $f = VI_\mu[e^{\lambda u}\chi_{B_2}]e^{\lambda u}$ and $\Gamma$ is as in \eqref{eq:fundamental_solution}. For $k\in \bb N$ consider 
\begin{equation*}
	u_k(x) = 2\log \frac{Ak}{1 + k^2|x|^2}, 
\end{equation*}
where $A = A(\mu) = \left(\frac{4(2 - \mu)}{\pi}\right)^{1/(4 - \mu)}$ and define 
\begin{equation*}
	V_k = \frac{I_\mu[e^{\lambda u_k}]}{I_\mu[e^{\lambda u_k}\chi_{B_2}]}. 
\end{equation*}	
These functions satisfy 
\begin{equation}
\label{eq:rigged_solutions}
	-\lap u_k(x) 
	= \frac{8k^2}{(1 + k^2|x|^2)^2} 
	= V_kI_\mu[e^{\lambda u_k}\chi_{B_2}](x)e^{\lambda u_k(x)}
	\qquad \text{ for }x\in\bb R^2.   
\end{equation}
The second of these equalities can be verified using the fact that $e^{\lambda u_k}$ is an extremal function for the sharp HLS inequality and thus
\begin{equation}
\label{eq:standard_Imu}
	I_\mu[e^{\lambda u_k}] = ce^{\mu u_k/4}
\end{equation}
for some constant $c>0$. See the proof of Theorem 1.1 of \cite{Gluck2025classification} for details. Theorem \ref{oldtheorem:classification} in Appendix \ref{s:appendix} gives the equality $\|e^{u_k}\|_{L^1(\bb R^2)}= \left(4(2 -\mu)\right)^{\frac 2{4 - \mu}}\pi^{\frac{2 - \mu}{4 - \mu}}$. Setting $F_k = V_kI_\mu[e^{\lambda u_k}\chi_{B_2}]e^{\lambda u_k}$, from the explicit expression of $F_k$ in \eqref{eq:rigged_solutions} we see that 
\begin{equation}
\label{eq:rigged_energy_smallness}
	\|F_k\|_{L^1(B(e_1, \frac 34))}\leq Ck^{-2}.
\end{equation}
Next we verify a uniform $L^\infty$ bound for the coefficient functions $V_k$. Evidently $V_k\geq 1$ in $\bb R^2$ for all $k$. Moreover, there is a $k$-independent constant $C>0$ such that
\begin{equation}
\label{eq:rigged_Vk_bound}
	V_k(x)\leq C
	\qquad \text{ for all }x\in B_{7/4}. 
\end{equation}
To verify \eqref{eq:rigged_Vk_bound}, note that for $x\in B_{7/4}$ we have
\begin{equation*}
\begin{split}
	I_\mu[e^{\lambda u_k}\chi_{\bb R^2\setminus B_2}](x)
	& 
	\leq 8^\mu\left(\frac Ak\right)^{2\lambda}\int_{\bb R^2\setminus B_2}|y|^{-(\mu + 4\lambda)}\; \d y\\
	& 
	\leq Ck^{-2\lambda}. 
\end{split}
\end{equation*}
Therefore, using \eqref{eq:standard_Imu}
when $k$ is sufficiently large for every $x\in B_{7/4}$ we have
\begin{equation*}
\begin{split}
	I_\mu[e^{\lambda u_k}\chi_{B_2}](x)
	& = I_\mu[e^{\lambda u_k}](x) - I_\mu[e^{\lambda u_k}\chi_{\bb R^2\setminus B_2}](x)\\
	& \geq C\left[\left(\frac{Ak}{1 + k^2|x|^2}\right)^{\mu/2} - k^{-2\lambda}\right]\\
	& 
	\geq C\left(\frac{Ak}{1 + k^2|x|^2}\right)^{\mu/2}\\
	& = CI_\mu[e^{\lambda u_k}](x)
\end{split}
\end{equation*}
from which \eqref{eq:rigged_Vk_bound} follows. Finally we show that in spite of the uniform bounds $1\leq V_k\leq C$ in $B_{7/4}$ and $\|e^{u_k}\|_{L^1(B_2)}
	\leq \|e^{u_k}\|_{L^1(\bb R^2)}\leq C$, and in spite of the energy smallness \eqref{eq:rigged_energy_smallness}, we still have $u_k\to-\infty$ uniformly on $\overline B(e_1, \frac 12)$ and $\varphi_k := u_k - \Gamma*F_k\to-\infty$ uniformly on $\overline B(e_1, \frac 12)$. The first of these limits follows immediately from the explicit expression of $u_k$. To verify the second of these limits, we note that Proposition 2.14 of \cite{Gluck2025classification} guarantees that for each $k$, there is a constant $C_k$ for which $\varphi_k \equiv C_k$ (i.e., $\varphi_k$ is independent of $x$). Thanks to Lemma \ref{lemma:bubble_potential} of Appendix \ref{s:appendix} we can compute the values of these constants. Indeed, using said lemma and the change of variable $y\mapsto ky$ we have 
\begin{equation*}
\begin{split}
	u_k(x)&  - 2\log\left(\frac{Ak}{2}\right)\\
	& = 2\log\left(\frac 2{1 + k^2|x|^2}\right)\\
	& = \frac 4\pi\int_{\bb R^2}\log\left(\frac{\sqrt 2}{|kx - y|}\right)\frac{1}{(1 + |y|^2)^2}\; \d y\\
	& = \frac 4\pi\int_{\bb R^2}\log\left(\frac{\sqrt 2}{k|x - y|}\right)\frac{k^2}{(1 + k^2|y|^2)^2}\; \d y\\
	& 
	= 4\log\left(\frac{\sqrt 2}k\right) + \Gamma* F_k(x)
\end{split}
\end{equation*}
which, upon rearranging, yields
\begin{equation*}
	u_k(x) - \Gamma*F_k(x)
	= 2\log A - 2\log k
	\to-\infty. 
\end{equation*}
\end{example}
The following proof is based on the proof of Theorem 3 of \cite{BrezisMerle1991} but includes necessary adjustments to handle the nonlocality. 
\begin{proof}[Proof of Theorem \ref{theorem:BM_alternative}]
Combining assumption \eqref{eq:vkp_uklambdaq_energy} with Lemma \ref{lemma:integrability_assumption} shows that the sequence of functions $f_k$ in \eqref{eq:fk} is bounded in $L^1(\omega)$. Therefore, there is a non-negative bounded measure $\eta$ on $\omega$ and a subsequence of $(f_k)_{k = 1}^\infty$ (whose members are still denoted $f_k$) along which $f_k\weakconv \eta$ weak-$*$ in the sense of measures. That is, 
\begin{equation*}
	\int_\omega f_k\psi \; \d x
	\to \int_\omega \psi\; \d\eta
\end{equation*}
for all $\psi\in C_c(\omega)$. 
To continue, we define a \emph{regular point} for $\eta$ as any point $x_0\in \omega$ for which there exists $\psi\in C_c(\omega)$ such that $0\leq \psi\leq 1$ and $\psi\equiv 1$ in a neighborhood of $x_0$ and for which 
\begin{equation*}
	\int_\omega \psi\; \d \eta < 4\pi\left(1 - \frac 1{2\lambda p}\right). 
\end{equation*}
Let $\Sigma\subset \omega$ denote the collection of non-regular points for $\eta$. Evidently $x_0\in \Sigma$ if and only if $\eta(\{x_0\})\geq 4\pi\left(1 - \frac 1{2\lambda p}\right)$. 
Since $\eta$ is a bounded measure with $\int_\omega \; \d \eta\leq \sup_k\|f_k\|_{L^1(\omega)}$, we deduce that $\Sigma$ is a finite set with 
\begin{equation}
\label{eq:finitely_many_nonregular}
	\card(\Sigma)\leq \left(4\pi(1 - \frac{1}{2\lambda p})\right)^{-1}\sup_k\|f_k\|_{L^1(\omega)},
\end{equation} 
where $\card(\Sigma)$ denotes the cardinality of $\Sigma$. The remainder of the proof is divided into three steps. \\
\begin{enumerate}[label = {\bf Step \arabic*.}, ref = {\bf Step \arabic*}, wide = 0pt]
	\item \label{item:blowup_iff_nonregular} We show that the blow-up set $S$ for $(u_k)_{k = 1}^\infty$ relative to $\omega$ as defined in Definition \ref{defn:blow_up_point} coincides with the set of non-regular points for $\eta$ (i.e., $S = \Sigma$). To verify the containment $S\subset \Sigma$, suppose $x_0\in \omega\setminus \Sigma$ (i.e., $x_0$ is regular for $\eta$), and choose $r>0$ small and $\psi\in C_c(\omega)$ for which both $0\leq \psi\leq 1$ and $\psi|_{B_{2r}(x_0)}\equiv 1$ and for which
	\begin{equation*}
		4\pi\left(1 - \frac 1{2\lambda p}\right)
		> \int_\omega \psi\; \d \eta
		= \lim_k\int_\omega f_k\psi \; \d x. 
	\end{equation*}
	Passing to a further subsequence, and from the properties of $\psi$ we obtain 
	\begin{equation*}
	\begin{split}
		4\pi\left(1 - \frac 1{2\lambda p}\right)
		& > \sup_k\int_\omega f_k \psi \; \d x\\
		& \geq \sup_k\int_{B_{2r}(x_0)} V_kI_\mu[e^{\lambda u_k}\chi_{\Omega_k}]e^{\lambda u_k}\; \d x. 
	\end{split}
	\end{equation*}
	Applying Lemma \ref{lemma:uniform_upper_bound} 
	implies $(u_k^+)_{k = 1}^\infty$ is bounded in $L^\infty(B_r(x_0))$ and hence $x_0\not\in S$. Next we show that $\Sigma\subset S$. Let $x_0\in \Sigma$. First we show by way of contradiction that 
	\begin{equation}
	\label{eq:concentration_point}
		\|u_k^+\|_{L^\infty(B_R(x_0))}\to\infty
		\qquad \text{ for all } R\in (0, \dist(x_0, \bdy\omega)).  
	\end{equation}
	Accordingly, suppose there is $R_0\in (0, \dist(x_0, \bdy\omega))$ and a subsequence of $(u_k)_{k = 1}^\infty$ along which $\|u_k^+\|_{L^\infty(B_{R_0}(x_0))}\leq C$ for all $k$ and pass to such a subsequence. In particular, along this subsequence we have $\sup_k\|e^{\lambda u_k}\|_{L^\infty(B_{R_0}(x_0))}\leq C$. Fix any $R\in (0, R_0)$, let $\psi\in C_c(\omega)$ satisfy $0\leq \psi\leq 1$, $\psi\equiv 1$ on $B_{R/2}(x_0)$, and $\supp \psi\subset B_R(x_0)$ and let $r$ be given by $\frac 1r = \frac 1q - \frac{2 - \mu}2$. 
	For each $k$, using H\"older's inequality 
	and the HLS inequality we have
	\begin{equation*}
	\begin{split}
		\int_\omega f_k\psi\; \d x
		& \leq \int_{B_R(x_0)}f_k\; \d x\\
		& 
		\leq |B_R|^{1 - \frac 1p - \frac 1 r}\|e^{\lambda u_k}\|_{L^\infty(B_R(x_0))}\|V_k\|_{L^p(\omega)}\|I_\mu[e^{\lambda u_k}]\|_{L^r(\Omega_k)}\\
		& \leq C|B_R|^{\frac 1q }\|e^{\lambda u_k}\|_{L^\infty(B_R(x_0))}\|V_k\|_{L^p(\omega)}\|e^{\lambda u_k}\|_{L^q(\Omega_k)}\\
		& \leq CR^{2/q}.
	\end{split}
	\end{equation*}
	Letting $k\to\infty$ yields 
	\begin{equation*}
		\int_\omega \psi \; \d \eta \leq CR^{2/q}, 
	\end{equation*}
	so by choosing $R$ sufficiently small we deduce that $x_0$ is a regular point for $\eta$. This contradicts the containment $x_0\in\Sigma$ and thereby establishes \eqref{eq:concentration_point}. To complete the proof that $x_0\in S$, choose $R>0$ small so that $\Sigma \cap B_{2R}(x_0) = \{x_0\}$. Choose $(x_k)_{k = 1}^\infty\subset B_R(x_0)$ for which $u_k^+(x_k)\to \infty$. For such $(x_k)_{k = 1}^\infty$ we must have $x_k\to x_0$. Indeed, if $x_k\not\to x_0$, then one may select $\tilde x\in \overline B_R(x_0)\setminus \{x_0\}$ and a subsequence of $(x_k)_{k = 1}^\infty$ along which $x_k\to \tilde x$. By the smallness of $R$ we have $\tilde x\not\in \Sigma$. On the other hand, by construction of $\tilde x$ we have $\tilde x\in S\subset \Sigma$, where the containment $S\subset \Sigma$ was established at the beginning of Step 1. This is a contradiction, so we deduce that $x_k\to x_0$. The containment  $x_0\in S$ is established. 
	\item \label{item:blowup_set_empty} We show that if $S=\emptyset$ then one of alternatives \ref{item:uniformly_bounded} or \ref{item:uniformly_to_-infty} holds. Accordingly suppose $S= \emptyset$ so that by Step 1, $\Sigma = \emptyset$ and $(u_k^+)_{k= 1}^\infty$ is bounded in $L^\infty_{\loc}(\omega)$. As a consequence of this bound we find that $(f_k)_{k = 1}^\infty$ is bounded in $L^p_{\loc}(\omega)$. To verify this, fix a compact set $K\subset\omega$ and choose $\delta>0$ such that $K_\delta \subset \omega$, where 
	\begin{equation*}
		K_\delta 
		= \{x\in \bb R^2: \dist(x, K)\leq \delta\}. 
	\end{equation*}
	For any $x\in K$ we have both
	\begin{equation*}
		I_\mu[e^{\lambda u_k}\chi_{K_\delta}](x)
		\leq C\|e^{u_k}\|_{L^\infty(K_\delta)}^\lambda (\diam \omega)^{2 - \mu}
	\end{equation*}
	and
	\begin{equation*}
	\begin{split}
		I_\mu[e^{\lambda u_k}\chi_{\Omega\setminus K_\delta}](x)
		& \leq \|e^{\lambda u_k}\|_{L^q(\Omega_k)}\left(\int_{\bb R^2\setminus B_\delta(x)}|x - y|^{-\mu q'}\; \d y\right)^{1/q'}\\
		& \leq C(p,\mu)c_0\delta^{-(\mu - \frac 2{q'})}, 
	\end{split}
	\end{equation*}
	where the assumption $p>2/\mu$ was used to ensure that $\mu q'>2$. Combining these two estimates shows that $(I_\mu[e^{\lambda u_k}\chi_{\Omega_k}])_{k = 1}^\infty$ is bounded in $L^\infty(K)$, so assumption \eqref{eq:vkp_uklambdaq_energy} guarantees that $(f_k)_{k = 1}^\infty$ is bounded in $L^p(K)$. Since $K$ is an arbitrary compact subset of $\omega$, we deduce that $(f_k)_{k = 1}^\infty$ is bounded in $L^1(\omega)\cap L^p_{\loc}(\omega)$. In particular, $\eta\in L^1(\omega)\cap L^p_{\loc}(\omega)$.
	Let $v_k$ be the solution to the problem 
	\begin{equation}
	\label{eq:vk_BVP_2}
	\begin{cases}
		-\lap v_k = f_k & \text{ in }\omega\\
		v_k = 0 &\text{ on }\bdy \omega. 
	\end{cases}
	\end{equation}
	The maximum principle ensures that $v_k\geq 0$ in $\omega$. 
	Moreover, $v_k\to v$ uniformly on compact subsets of $\omega$, where $v$ is the solution to 
	\begin{equation}
	\label{eq:v_BVP}
	\begin{cases}
		-\lap v = \eta & \text{ in }\omega\\
		v = 0 & \text{ on }\bdy\omega. 
	\end{cases}
	\end{equation}
	Define $w_k = u_k -v_k$ so that $w_k$ satisfies 
	\begin{equation}
	\label{eq:wk_BVP}
	\begin{cases}
		-\lap w_k = 0 & \text{ in }\omega\\
		w_k = u_k & \text{ on }\bdy \omega, 
	\end{cases}
	\end{equation}
	and by the non-negativity of $v_k$, we have $w_k\leq u_k$. In particular $w_k^+\leq u_k^+$ so $(w_k^+)_{k = 1}^\infty$ is bounded in $L^\infty_{\loc}(\omega)$. A routine argument involving Harnack's inequality for positive harmonic functions shows that either $w_k\to-\infty$ uniformly on compact subsets of $\omega$ or there is a subsequence of $(w_k)_{k= 1}^\infty$ that is bounded in $L^\infty_{\loc}(\omega)$. 
	These two alternatives for $(w_k)_{k = 1}^\infty$ correspond to alternatives \ref{item:uniformly_to_-infty} or \ref{item:uniformly_bounded} respectively for $(u_k)_{k = 1}^\infty$. 
	\item \label{item:blowup_set_nonempty} We show that if $S\neq \emptyset$ then alternative \ref{item:finite_blow_up} holds. Accordingly, suppose $S\neq \emptyset$ and note that by \ref{item:blowup_iff_nonregular} and inequality \eqref{eq:finitely_many_nonregular}, there are finitely many points $a^1, \ldots, a^m\in \Omega$ for which $S = \{a^1, \ldots, a^m\}$. In particular $\eta(\{a^i\})\geq 4\pi\left(1 - \frac 1{2\lambda p}\right)$ for all $i\in \{1, \ldots, m\}$. Next, an argument similar to the one carried out at the beginning of \ref{item:blowup_set_empty} shows that $(f_k)_{k = 1}^\infty$ is bounded in $L^p_{\loc}(\omega\setminus S)$ and thus $\eta\in L^p_{\loc}(\Omega\setminus S)$. 
	As in \ref{item:blowup_set_empty}, we decompose $u_k = v_k + w_k$, where $v_k$ and $w_k$ satisfy \eqref{eq:vk_BVP_2} and \eqref{eq:wk_BVP} respectively. Similarly to \ref{item:blowup_set_empty} we have $v_k\geq 0$ and $v_k\to v$ uniformly on compact subsets of $\omega\setminus S$, where $v$ is the solution to \eqref{eq:v_BVP}. 
	Moreover, as in \ref{item:blowup_set_empty} a routine argument involving Harnack's inequality implies that either $w_k\to-\infty$ uniformly on compact subsets of $\omega \setminus S$ or there is a subsequence of $(w_k)_{k = 1}^\infty$ that is bounded in $L^\infty_{\loc}(\omega\setminus S)$. To complete the proof of \ref{item:blowup_set_nonempty} it suffices to show that there is no subsequence of $(w_k)_{k = 1}^\infty$ that is bounded in $L^\infty_{\loc}(\omega\setminus S)$. Proceeding by way of contradiction, suppose $(w_k)_{k = 1}^\infty$ is such a subsequence. Fix $a\in S$ and $R>0$ sufficiently small so that $B_{2R}(a)\cap S = \{a\}$. Evidently both of $(w_k)_{k = 1}^\infty$ and $(v_k)_{k = 1}^\infty$ are bounded in $L^\infty(\bdy B_R(a))$, so there is a constant $C_1>0$ such that $\inf_{\bdy B_R(a)}u_k\geq -C_1$ for all $k$. Letting $z_k$ be the solution to 
	\begin{equation*}
	\begin{cases}
		-\lap z_k = f_k & \text{ in }B_R(a)\\
		z_k = -C_1 & \text{ on }\bdy B_R(a), 
	\end{cases}
	\end{equation*}
	the maximum principle guarantees that $u_k\geq z_k$ in $B_R(a)$ and thus
	\begin{equation}
	\label{eq:zk_energy_bound}
	\begin{split}
		\|I_\mu[e^{\lambda z_k}\chi_{B_R(a)}]\|_{L^{p'}(B_R(a))}
		\leq \|I_\mu[e^{\lambda u_k}\chi_{\Omega_k}\|_{L^{p'}(\Omega_k)}
		\leq C. 
	\end{split}
	\end{equation}
	Since $(f_k)_{k = 1}^\infty$ is bounded in $L^p_{\loc}(B_R(a)\setminus\{a\})$, by a standard argument, after passing to a suitable subsequence we have $z_k\to z$ uniformly on compact subsets of $B_R(a)\setminus\{a\}$ (and hence a.e. on $B_R(a)$), where $z$ is the solution to 
	\begin{equation*}
	\begin{cases}
		-\lap z = \eta & \text{ in }B_R(a)\\
		z = -C_1 &\text{ on }\bdy B_R(a). 
	\end{cases}
	\end{equation*}
	We proceed to show that this convergence, combined with the assumption that $a$ is not a regular point for $\eta$ is incompatible with estimate \eqref{eq:zk_energy_bound}. Since $a$ is not a regular point for $\eta$ we have $\eta(\{a\})\geq 4\pi\left(1 - \frac{1}{2\lambda p}\right)$ and thus $\eta \geq 4\pi\left(1 - \frac{1}{2\lambda p}\right) \delta_{a}$. Setting $\Gamma_{a}(x) = -\frac 1{2\pi}\log|x - a|$ we have 
	\begin{equation*}
	\begin{cases}
		-\lap(z - 4\pi\left(1 - \frac{1}{2\lambda p}\right)\Gamma_{a}) 
		\geq 0 & \text{ in }B_R(a)\\
		z - 4\pi\left(1 - \frac{1}{2\lambda p}\right) \Gamma_{a} = -C_1 + 2\left(1 - \frac 1{2\lambda p}\right)\log R & \text{ on }\bdy B_R(a). 
	\end{cases}
	\end{equation*}
	The maximum principle yields
	\begin{equation*}
	\begin{split}
		z(y)
		& 
		\geq -2\left(1 - \frac 1{2\lambda p}\right)\log|y - a| + O(1)
	\end{split}
	\end{equation*}
	as $y\to a$ and thus
	\begin{equation}
	\label{eq:elambdaz_lower}
		e^{\lambda z(y)}
		\geq C|y - a|^{\frac 1 p - 2\lambda}
		\qquad \text{ for }y\in B_R(a)\setminus\{a\}. 
	\end{equation}
\end{enumerate}
For $x\in B_R(a)\setminus \{a\}$ and $y\in B_{|x - a|/2}(a)$ we have $|x- y| \leq 
2|x - a|$ and thus, 
\begin{equation*}
\begin{split}
	I_\mu[e^{\lambda z}\chi_{B_R(a)}](x)
	& \geq C\int_{B_{|x - a|/2}(a)}\frac{1}{|y - a|^{2\lambda - \frac 1 p}}\cdot \frac{1}{|x- y|^\mu}\; \d y\\
	& \geq \frac{C}{|x - a|^\mu}\int_{B_{|x - a|/2}(a)}\frac{1}{|y - a|^{2\lambda - \frac 1 p}}\; \d y\\
	& \geq C|x - a|^{\frac 1p - \frac \mu 2}.
\end{split}
\end{equation*}
Combining this estimate with estimate \eqref{eq:elambdaz_lower} we find that 
\begin{equation*}
\begin{split}
	I_\mu[e^{\lambda z}\chi_{B_R(a)}](x)e^{\lambda z(x)}
	& 
	\geq C|x - a|^{-\frac2{p'}} 
\end{split}
\end{equation*}
whenever $x\in B_R(a)\setminus\{a\}$ and thus $\|I_\mu[e^{\lambda z}\chi_{B_R(a)}]e^{\lambda z}\|_{L^{p'}(B_R(a))} = +\infty$. Finally, two applications of Fatou's Lemma give
\begin{equation*}
\begin{split}
	\liminf_k & \int_{B_R(a)}\left[I_\mu[e^{\lambda z_k}\chi_{B_R(a)}]e^{\lambda z_k}\right]^{p'}\; \d x\\
	& 
	\geq \int_{B_R(a)}\left[I_\mu[e^{\lambda z}\chi_{B_R(a)}]e^{\lambda z}\right]^{p'}\; \d x\\
	& = +\infty, 
\end{split}
\end{equation*}
which contradicts estimate \eqref{eq:zk_energy_bound}. 
\end{proof}
\section{Proof of the $\sup +\inf$ inequality}
\label{s:sup_inf_inequality}
This section is devoted to the proof of Theorem \ref{theorem:sup_inf_inequality}. Our proof will follow the general strategy outlined in Remark 2.3 of \cite{EspositoLucia2021}. Theorem \ref{theorem:sup_inf_inequality} is a consequence of the following proposition. 
\begin{prop}
\label{prop:reduced_sup+inf_inequality}
Let $\omega\subset \bb R^2$ be a bounded domain, let $\mu\in (0, 2)$, and let $\lambda$ be as in \eqref{eq:lambda}. Suppose $0< a\leq b< \infty$ and $\Lambda\subset\Lambda_{a,b}(\omega)$ is a subset that is equicontinuous at each point of $\omega$. For each compact subset $K\subset\omega$, each $c_0> 0$, and each $C_1>1$ there is a constant $C = C(\Lambda, K, \omega, c_0, C_1)>0$ such that for all $V\in \Lambda$, if $\Omega\subset \bb R^2$ is a bounded domain for which $\omega\subset \Omega$ and if $u$ is a distributional solution to \eqref{eq:nonlocal_version} for which
\begin{equation}
\label{eq:sup+inf_lower_bound_assumption}
	\max_K u + C_1\inf_\omega u\geq 0, 
\end{equation}
then $\max_Ku \leq C$. We emphasize that $C$ is independent of $\Omega$.  
\end{prop}
Before providing the proof of Proposition \ref{prop:reduced_sup+inf_inequality} we show that it implies Theorem \ref{theorem:sup_inf_inequality}. 
\begin{proof}[Proof of Theorem \ref{theorem:sup_inf_inequality}]
Fix a compact set $K\subset \omega$, $c_0>0$, $C_1>1$, and $V\in \Lambda$. Let $\Omega \subset \bb R^2$ be a bounded domain for which $\omega\subset \Omega$ and let $u$ be a solution to \eqref{eq:nonlocal_version}. If $\max_Ku + C_1\inf_\omega u< 0$ then there is nothing to prove. Otherwise, 
Proposition \ref{prop:reduced_sup+inf_inequality} guarantees the existence of a positive constant constant $C = C(\Lambda, K, \omega, c_0, C_1)$ for which $\max_K u\leq C$. For any such $C$ we have
\begin{equation*}
\begin{split}
	\max_K u + C_1\inf_\omega u
	& 
	\leq C + C_1\max_Ku\\
	& \leq C(1 + C_1). 
\end{split}
\end{equation*}

\end{proof}
The remainder of this section is devoted to the proof of Proposition \ref{prop:reduced_sup+inf_inequality}. We will need the following lemma which yields a uniform bound on the norms $\|-\lap u_k\|_{L^1(B_k)}$ over suitably chosen shrinking balls $B_k$ for any sequence $u_k$ that violates the assertion of Proposition \ref{prop:reduced_sup+inf_inequality}, see \eqref{eq:small_ball_limiting_energy} below. For convenience we provide a proof of the lemma in Appendix \ref{s:appendix}. 
\begin{lemma}
\label{lemma:sup+inf_driving_estimate}
Let $\omega\subset \bb R^2$, let $0\leq f\in L^1(\omega)$ in $\omega$ and suppose $u$ satisfies $-\lap u = f$ in $\omega$. For any $B_\rho(x_0)\subset \omega$ there holds
\begin{equation*}
	u(x_0) - \inf_\omega u
	\geq \frac 1{2\pi}\int_{B_r(x_0)}f(x)\; \d x \; \log\frac\rho r
	\qquad \text{ for all }r\in (0, \rho). 
\end{equation*}
\end{lemma}
\begin{proof}[Proof of Proposition \ref{prop:reduced_sup+inf_inequality}]
Proceeding by way of contradiction, suppose the proposition is false and choose a compact subset $K\subset\omega$, $c_0> 0$, $C_1>1$, and sequences $(V_k)_{k = 1}^\infty\subset \Lambda$, $(\Omega_k)_{k = 1}^\infty$ with $\omega\subset \Omega_k$ and $(u_k)_{k = 1}^\infty$ of distributional solutions to 
\begin{equation}
\label{eq:contradictive_sequence_PDE_and_energy}
\begin{cases}
	-\lap u_k = V_kI_\mu[e^{\lambda u_k}\chi_{\Omega_k}]e^{\lambda u_k} & \text{ in }\omega\\
	\|e^{u_k}\|_{L^1(\Omega_k)}\leq c_0&\text{ for all }k\\ 
\end{cases}
\end{equation}
for which 
\begin{equation}
\label{eq:max_inf_bounded_below}
	\max_Ku_k + C_1\inf_{\omega} u_k\geq 0
\end{equation}
but 
\begin{equation}
\label{eq:K_maximizer}
	u_k(\tilde x_k)\geq k, 
\end{equation}
where $\tilde x_k\in K$ satisfies $u_k(\tilde x_k) = \max_Ku_k$. Fix $\rho>0$ for which $K_{3\rho}\subset \omega$, where for $r>0$ we define
\begin{equation*}
	K_r:= \{x\in \omega: \dist(x, K)\leq r\}. 
\end{equation*}
For each $k\in \bb N$, applying Lemma \ref{lemma:selection_process} of Appendix \ref{s:appendix} with $a = 2$ to the function $\varphi = \varphi_k= e^{u_k}$ over $B_\rho(\tilde x_k)$ yields $x_k\in B_\rho(\tilde x_k)$ for which both 
\begin{equation}
\label{eq:first_selection_property}
	u_k(x_k) + 2\log 2\geq \max\{u_k(x): x\in \overline B(x_k, r_k)\}
\end{equation}
and
\begin{equation}
\label{eq:second_selection_property}
	u_k(\tilde x_k) 
	\leq u_k(x_k) + 2\log\frac{2r_k}\rho, 
\end{equation}
where $r_k = \frac 12(\rho - |x_k - \tilde x_k|)< \frac \rho 2$. In particular, property \eqref{eq:second_selection_property} guarantees that $u_k(x_k)\geq u_k(\tilde x_k)\geq k$, so defining $\delta_k = \exp(-u_k(x_k)/2)$ we have $\delta_k\to 0$. Moreover, properties \eqref{eq:K_maximizer} and \eqref{eq:second_selection_property} guarantee that $\frac{r_k}{\delta_k}\to\infty$, so for any $R\gg 1$ we have $B_{R\delta_k}(x_k)\subset B_{r_k}(x_k)\subset K_{2\rho}$ whenever $k$ is sufficiently large. After passing to a suitable subsequence of $(u_k)_{k = 1}^\infty$ (and continuing to denote the members of such a subsequence by $u_k$), applying Lemma \ref{lemma:sup+inf_driving_estimate} to $B_\rho(x_k)$ with $r = R\delta_k< \rho$ we obtain 
\begin{equation*}
	u_k(x_k) - \inf_\omega u_k
	\geq \frac 1{2\pi}\int_{B(x_k, R\delta_k)}V_kI_\mu[e^{\lambda u_k}\chi_{\Omega_k}]e^{\lambda u_k}\; \d x\; \log \frac{\rho}{R\delta_k}. 
\end{equation*}
Upon rearranging this estimate and in view of assumption \eqref{eq:max_inf_bounded_below} we have 
\begin{equation*}
\begin{split}
	\int_{B(x_k, R\delta_k)}V_kI_\mu[e^{\lambda u_k}\chi_{\Omega_k}]e^{\lambda u_k}\; \d x
	& \leq \frac{4\pi u_k(x_k)}{2\log\frac\rho R + u_k(x_k)}\left(1 - \frac{\inf_\omega u_k}{u_k(x_k)}\right)\\
	& 
	\leq \frac{4\pi u_k(x_k)}{2\log\frac\rho R + u_k(x_k)}\left(1 + \frac 1{C_1}\right). 
\end{split}
\end{equation*}
Combining this estimate with the inequality $u_k(x_k)\geq k$ guarantees that for any $R\gg 1$ there holds
\begin{equation}
\label{eq:small_ball_limiting_energy}
	\limsup_{k\to\infty}\int_{B(x_k, R\delta_k)}V_kI_\mu[e^{\lambda u_k}\chi_{\Omega_k}]e^{\lambda u_k}\; \d x
	\leq 4\pi\left(1 + \frac 1{C_1}\right). 
\end{equation}
For each $k$ define 
\begin{equation*}
	v_k(y)
	= u_k(x_k + \delta_k y) + 2\log \delta_k
	\qquad \text{ for }y\in \Sigma_k:= \frac{\Omega_k - x_k}{\delta_k}. 
\end{equation*}
Evidently $v_k$ satisfies
\begin{equation}
\label{eq:rescaled_PDE_and_others}
\begin{cases}
	-\lap v_k = H_k I_\mu[e^{\lambda v_k}\chi_{\Sigma_k}]e^{\lambda v_k} &\text{ in }\omega_k\\
	v_k(0) = 0\\
	v_k(y)\leq 2\log 2 & \text{ in }\overline B_{r_k\delta_k^{-1}} \subset \omega_k\\
	\int_{\Sigma_k}e^{v_k}\leq c_0, 
\end{cases}
\end{equation}
where $H_k(y) = V_k(x_k + \delta_ky)$, $\omega_k = (\omega - x_k)/\delta_k$, and the uniform upper bound on $v_k$ in $\overline B_{r_k\delta_k^{-1}}$ follows from \eqref{eq:first_selection_property}. The compactness of $K_\rho$ 
guarantees the existence of $x_*\in K_\rho$ and a subsequence of $(x_k)_{k = 1}^\infty$ along which $x_k\to x_*$. Passing to a further subsequence we assume in addition that $V_k(x_*)\to \tau$ for some $\tau\in [a,b]$. Moreover, the equiconuity of $\Lambda$ at $x_*$ guarantees that $H_k\to \tau$ locally uniformly on $\bb R^2$. Using this local uniform convergence we find that for any $R\gg 1$, along a suitable subsequence there holds
\begin{equation*}
\begin{split}
	\int_{B(x_k, R\delta_k)}V_kI_\mu[e^{\lambda u_k}\chi_{\Omega_k}]e^{\lambda u_k}
	& = \int_{B_R}H_k I_\mu[e^{\lambda v_k}\chi_{\Sigma_k}]e^{\lambda v_k}\\
	& = \tau\int_{B_R}I_\mu[e^{\lambda v_k}\chi_{\Sigma_k}]e^{\lambda v_k} + \circ(1), 
\end{split}
\end{equation*}
where $\circ(1)\to 0$ as $k\to\infty$. 
Combining this estimate with \eqref{eq:small_ball_limiting_energy} gives 
\begin{equation}
\label{eq:rescaled_small_ball_limiting_energy}
	\limsup_{k\to\infty}\int_{B_R}I_\mu[e^{\lambda v_k}\chi_{\Sigma_k}]e^{\lambda v_k}
	\leq \frac{4\pi}{\tau}\left(1 + \frac 1{C_1}\right). 
\end{equation}
Independently, for any $R\gg 1$, applying Theorem \ref{theorem:BM_alternative} with $p = +\infty$ and with $\Omega = B_R\subset\Sigma_k$ yields a subsequence of $(v_k)_{k = 1}^\infty$ that is bounded in $L^\infty_{\loc}(B_R)$. 
Routine computations now show that $(H_kI_\mu[e^{\lambda v_k}\chi_{\Sigma_k}]e^{\lambda v_k})_{k = 1}^\infty$ is bounded in $L^\infty(B_{R/2})$. Indeed, for any $y\in B_{R/2}$ we have both 
\begin{equation*}
	H_k(y)I_\mu[e^{\lambda v_k}\chi_{\Sigma_k}](y)e^{\lambda v_k(y)}
	\leq be^{2\lambda \log 2}I_\mu[e^{\lambda v_k}\chi_{\Sigma_k}](y)
\end{equation*}
and 
\begin{equation*}
\begin{split}
	I_\mu[e^{\lambda v_k}\chi_{\Sigma_k}](y)
	& = \int_{B_R}\frac{e^{\lambda v_k(z)}}{|y - z|^\mu}\; \d z + \int_{\Sigma_k\setminus B_R}\frac{e^{\lambda v_k(z)}}{|y - z|^\mu}\; \d z\\
	& \leq e^{2\lambda\log 2}\int_{B_{2R}(y)}|y - z|^{-\mu}\; \d z \\
	& \quad + \|e^{v_k}\|_{L^1(\Sigma_k)}^\lambda \left(\int_{\bb R^2\setminus B_{R/2}(y)}|y - z|^{-4}\; \d z\right)^{1- \lambda} 
	\\
	& 
	\leq C(\mu, c_0, R). 
\end{split}
\end{equation*}
Standard elliptic estimates guarantee the existence of $\alpha\in(0, 1)$ for which $(v_k)_{k = 1}^\infty$ is bounded in $C^{1, \alpha}(B_{R/4})$. 
The Arzel\`a-Ascoli Theorem implies the existence of $v^{(R)}\in C^{1, \alpha}(B_{R/4})$ and a subsequence of $(v_k)_{k = 1}^\infty$ along which $v_k\to v^{(R)}$ in $C^{1, \alpha}(B_{R/4})$. Upon repeating this argument with $R$ replaced by a sequence $(R_\ell)_{\ell = 1}^\infty$ for which $R_\ell\to\infty$ and employing a standard diagonal subsequence argument, one finds that there is $v\in C^{1,  \alpha}(\bb R^2)$ satisfying
\begin{equation}
\label{eq:sup+inf_lemma_limiting_problem}
\begin{cases}
	-\lap v = \tau I_\mu[e^{\lambda v}]e^{\lambda v} & \text{ in }\bb R^2\\
	v(0) = 0\\
	v(y)\leq 2\log 2 & \text{ for all }y\in \bb R^2\\
	\int_{\bb R^2}e^v\leq c_0
\end{cases}
\end{equation}
and a subsequence of $(v_k)_{k = 1}^\infty$ along which $v_k\to v$ in $C^{1, \alpha}_{\loc}(\bb R^2)$. We remark that due to the nonlocality in the nonlinearity, the verification that $v$ satisfies the first equality in \eqref{eq:sup+inf_lemma_limiting_problem} is more involved than in the local setting (e.g., for problems of the form $-\lap u = e^u$). For the reader's convenience, we provide the details of the verification immediately after the conclusion of the present proof. Applying Theorem \ref{oldtheorem:classification} of Appendix \ref{s:appendix} to the function $v + \frac{1}{2\lambda}\log \tau$ guarantees the existence of $(y_0, \delta)\in \bb R^2\times (0, \infty)$ for which 
\begin{equation*}
	v(y) = 2\log\frac{\delta}{1 + \delta^2|y - y_0|^2} + \frac{2}{4 - \mu}\log\frac{4(2 - \mu)}{\pi \tau}. 
\end{equation*}
In particular, in view of the second equality in \eqref{eq:asserted_energies} we have
\begin{equation}
\label{eq:global_limiting_energy}
	\int_{\bb R^2}I_\mu[e^{\lambda v}]e^{\lambda v} = \frac{8\pi}{\tau}. 
\end{equation}
In the remainder of the proof we show that equation \eqref{eq:global_limiting_energy} is incompatible with the $C^{1, \alpha}_{\loc}(\bb R^2)$ convergence $v_k\to v$ and inequality \eqref{eq:rescaled_small_ball_limiting_energy}. Let $\epsilon>0$ and choose $R = R(\epsilon)\gg 1$ for which both 
\begin{equation*}
	\int_{B_R}I_\mu[e^{\lambda v}]e^{\lambda v} > \frac{8\pi(1 - \epsilon)}{\tau}
\end{equation*}
and
\begin{equation*}
	\int_{\bb R^2}I_\mu[e^{\lambda v}\chi_{\bb R^2\setminus B_R}]e^{\lambda v}
	< \frac{8\pi\epsilon}{\tau}.  
\end{equation*}
For any $k$ large enough to satisfy $B_R\subset\Sigma_k$ we have
\begin{equation}
\label{eq:limiting_function_large_ball_energy}
\begin{split}
	\frac{8\pi(1 - \epsilon)}{\tau}
	& < \int_{B_R}I_\mu[e^{\lambda v}]e^{\lambda v}\\
	& \leq \int_{B_R}I_\mu[e^{\lambda v}\chi_{B_R}]e^{\lambda v} + \int_{\bb R^2}I_\mu[e^{\lambda v}\chi_{\bb R^2\setminus B_R}]e^{\lambda v}\\
	& \leq \int_{B_R}I_\mu[e^{\lambda v_k}\chi_{\Sigma_k}]e^{\lambda v_k}+ E_k + \frac{8\pi\epsilon}{\tau}, 
\end{split}
\end{equation}
where 
\begin{equation*}
	E_k
	= \int_{B_R}|I_\mu[(e^{\lambda v}-e^{\lambda v_k})\chi_{B_R}]|e^{\lambda v}
	+ \int_{B_R}I_\mu[e^{\lambda v_k}\chi_{B_R}]|e^{\lambda v} - e^{\lambda v_k}|. 
\end{equation*}
From H\"older's inequality and the HLS inequality we have 
\begin{equation*}
\begin{split}
	E_k
	\leq & \; \|e^v\|_{L^1(\bb R^2)}^\lambda \|I_\mu[(e^{\lambda v} - e^{\lambda v_k})\chi_{B_R}\|_{L^{4/\mu}(B_R)}\\
	&  + \|I_\mu[e^{\lambda v_k}\chi_{B_R}]\|_{L^{4/\mu}(B_R)}\|e^{\lambda v} - e^{\lambda v_k}\|_{L^{1/\lambda}(B_R)}\\
	\leq & \; C\left(\|e^v\|_{L^1(\bb R^2)}^\lambda + \|e^{v_k}\|_{L^1(\Sigma_k)}^\lambda\right)\|e^{\lambda v} - e^{\lambda v_k}\|_{L^{1/\lambda}(B_R)}, 
\end{split}
\end{equation*}
so 
since $v_k\to v$ in $C^{1, \alpha}(B_R)$ we have $E_k\to 0$ as $k\to\infty$. Now coming back to \eqref{eq:limiting_function_large_ball_energy} and in view of \eqref{eq:rescaled_small_ball_limiting_energy}, for any $\epsilon>0$ we have
\begin{equation*}
	\frac{8\pi(1 - 2\epsilon)}{\tau}
	\leq \limsup_{k\to\infty}\int_{B_R}I_\mu[e^{\lambda v_k}\chi_{\Sigma_k}]e^{\lambda v_k}
	\leq \frac{4\pi}{\tau}\left(1 + \frac 1{C_1}\right).   
\end{equation*}
In view of the assumption $C_1>1$, we may choose $\epsilon$ for which $4\epsilon< 1 - \frac 1{C_1}$ and thereby obtain a contradiction. 
\end{proof}
\begin{proof}[Verification of the first item in \eqref{eq:sup+inf_lemma_limiting_problem}]
\label{proof:verification}
In the following we verify that the $C^{1, \alpha}_{\loc}(\bb R^2)$-limit $v$ of the sequence $(v_k)_{k = 1}^\infty$ constructed immediately prior to equations \eqref{eq:sup+inf_lemma_limiting_problem} satisfies the PDE in \eqref{eq:sup+inf_lemma_limiting_problem} in the distributional sense. Fix $\varphi\in C_c^\infty(\bb R^2)$ and $\epsilon>0$. Choose $R = R(\epsilon, \varphi)\gg1$ so that both $\supp \varphi\subset B_R$ and $c_0^{2\lambda}\|\varphi\|_{L^{4/\mu}(\bb R^2)}R^{-\frac\mu 2}< \epsilon$. Since $\|v_k - v\|_{C^0(B_R)} + \|H_k - \tau\|_{C^0(B_R)} = \circ(1)$ we have 
\begin{equation}
\label{eq:initial_limiting_solution_test}
	\bigg|\int_{\bb R^2} v\lap \varphi + \tau\int_{\bb R^2}I_\mu[e^{\lambda v}]e^{\lambda v}\varphi\bigg|
	\leq  \circ(1) + \tau A_k,  
\end{equation}
where 
\begin{equation*}
	A_k
	= \abs{\int_{\bb R^2}\left(I_\mu[e^{\lambda v_k}\chi_{\Sigma_k}]e^{\lambda v_k} - I_\mu[e^{\lambda v}]e^{\lambda v}\right)\varphi}.
\end{equation*}
To estimate $A_k$, we first observe that $A_k\leq \sum_{j= 1}^4A_k^j$, where 
\begin{equation*}
\begin{split}
	A_k^1 & = \int_{\bb R^2}I_\mu[e^{\lambda v}]|e^{\lambda v} - e^{\lambda v_k}||\varphi|\\
	A_k^2 & = \int_{\bb R^2}I_\mu[|e^{\lambda v_k} - e^{\lambda v}|\chi_{B_{2R}}]e^{\lambda v_k}|\varphi|\\
	A_k^3 & = \int_{\bb R^2}I_\mu[|e^{\lambda v_k} - e^{\lambda v}|\chi_{\Sigma_k\setminus B_{2R}}]e^{\lambda v_k}|\varphi|\\
	A_k^4& = \int_{\bb R^2}I_\mu[e^{\lambda v}\chi_{\bb R^2\setminus \Sigma_k}]e^{\lambda v_k}|\varphi|,
\end{split}	
\end{equation*}
and then we separately estimate $A_k^j$ for $j = 1, \ldots, 4$. To estimate $A_k^1$, we use the Mean-Value Theorem, the fact that each of $v_k$ and $v$ is uniformly bounded above (independently of $k$) on $\supp \varphi$, H\"older's inequality and the HLS inequality to obtain 
\begin{equation*}
	A_k^1
	\leq C R^{2\lambda}\|\varphi\|_{C^0}\|e^v\|_{L^1(\bb R^2)}^\lambda \|v - v_k\|_{C^0(B_R)}
	= \circ(1). 
\end{equation*}
To estimate $A_k^2$, observe that for $y\in B_R$ there holds
\begin{equation*}
\begin{split}
	I_\mu[\abs{e^{\lambda v_k}- e^{\lambda v}}\chi_{B_{2R}}](y)
	& \leq C\|v_k- v\|_{C^0(B_{2R})}\int_{B_{4R}(y)}|y - z|^{-\mu}\; \d z\\
	& \leq CR^{2- \mu}\|v_k - v\|_{C^0(B_{2R})}
\end{split}
\end{equation*}
and therefore an application of H\"older's inequality gives
\begin{equation*}
	A_k^2
	\leq CR^{2 - \frac\mu 2}c_0^\lambda\|v_k - v\|_{C^0(B_{2R})}
	= \circ(1). 
\end{equation*}
To estimate $A_k^3$ observe that for $y\in \supp \varphi\subset B_R$ we have 
\begin{equation*}
\begin{split}
	I_\mu[\abs{e^{\lambda v_k} - e^{\lambda v}}\chi_{\Sigma_k \setminus B_{2R}}](y)
	& = \int_{\Sigma_k\setminus B_{2R}}\frac{|e^{\lambda v_k(z)} - e^{\lambda v(z)}|}{|y - z|^\mu}\; \d z\\
	& \leq \|e^{\lambda v_k} - e^{\lambda v}\|_{L^{1/\lambda}(\Sigma_k)}\left(\int_{\bb R^2\setminus B_R(y)}|y - z|^{-4}\; \d z\right)^{1 - \lambda}\\
	& \leq C\left(\|e^{\lambda v_k}\|_{L^{1/\lambda}(\Sigma_k)} +  \|e^{\lambda v}\|_{L^{1/\lambda}(\bb R^2)}\right)R^{-2(1 - \lambda)}\\
	& \leq Cc_0^\lambda R^{-\mu/2}. 
\end{split}
\end{equation*}
Therefore, an application of H\"older's inequality gives
\begin{equation*}
\begin{split}
	A_k^3
	& \leq Cc_0^\lambda R^{-\mu/2}\int_{B_R}e^{\lambda v_k}|\varphi|\\
	& \leq Cc_0^{2\lambda}\|\varphi\|_{L^{4/\mu}(\bb R^2)}R^{-\mu/2}\\
	& \leq C\epsilon, 
\end{split}
\end{equation*}
where the final inequality holds by the largeness assumption on $R$. 
Finally, the estimate of $A_k^4$ is similar to that of $A_k^3$. For every $y\in B_R$ we have the pointwise estimate
\begin{equation*}
\begin{split}
	I_\mu[e^{\lambda v}\chi_{\bb R^2\setminus \Sigma_k}](y)
	\leq I_\mu[e^{\lambda v}\chi_{\bb R^2\setminus B_{2R}}](y)
	\leq Cc_0^\lambda R^{-\frac\mu2}, 
\end{split}
\end{equation*}
so H\"older's inequality and the largeness assumption on $R$ give $A_k^4\leq C\epsilon$. 
Combining the estimates of $A_k^1,\ldots, A_k^4$ we find that $A_k\leq C\epsilon + \circ(1)$.  Since $\epsilon>0$ is arbitrary bringing this estimate back to \eqref{eq:initial_limiting_solution_test} shows that $v$ is indeed a distributional solution to the PDE in \eqref{eq:sup+inf_lemma_limiting_problem}. 
\end{proof} 
\section{Quantization}
\label{s:quantization}
In this section we provide a proof of Theorem \ref{theorem:quantization}. The primary task in the proof is to establish the following proposition. 
\begin{prop}
\label{prop:energy_quantization}
Let $\Omega\subset \bb R^2$ be a bounded domain for which $\overline B_R\subset\Omega$, let $\mu\in (0, 2)$ and let $\lambda$ be as in \eqref{eq:lambda}. Let $V\in C^0(\overline B_R)$ be a nonnegative function and let $(V_k)_{k = 1}^\infty$ be a sequence of nonnegative functions on $\overline B_R$ for which $V_k\to V$ in $C^0(\overline B_R)$. If $(u_k)_{k = 1}^\infty$ is a sequence of distributional solutions to 
\begin{equation}
\label{eq:BR_PDE_and_energy}
\begin{cases}
	-\lap u_k = V_kI_\mu[e^{\lambda u_k}\chi_\Omega]e^{\lambda u_k}
	& \text{ in }B_R\\
	\|e^{u_k}\|_{L^1(\Omega)}\; \d x \leq c_0
\end{cases}
\end{equation}
for which 
\begin{equation}
\label{eq:max_to_infty}
	\max_{\overline B_R}u_k\to \infty, 
\end{equation}
for which
\begin{equation}
\label{eq:uniformly_down_away_from_origin}
	\max_{\overline B_R\setminus B_r}u_k\to -\infty
	\qquad \text{ for all }r\in (0, R), 
\end{equation}
and for which
\begin{equation}
\label{eq:BR_energy_limit}
	\lim_{k\to\infty}\int_{B_R}V_kI_\mu[e^{\lambda u_k}\chi_\Omega]e^{\lambda u_k}\; \d x 
	= \alpha,
\end{equation}
then there is a positive integer $N$ for which $\alpha = 8\pi N$. 
\end{prop}
Before proving Proposition \ref{prop:energy_quantization}, let us show that it implies Theorem \ref{theorem:quantization}. 
\begin{proof}[Proof of Theorem \ref{theorem:quantization}]
Let $S = \{a^1, \ldots, a^m\}\subset\omega$ be as in alternative \ref{item:finite_blow_up} in Theorem \ref{theorem:BM_alternative}. For any $\varphi\in C_c(\omega)$, and any $R\in (0, 1)$ for which both $\dist(a^i, \bdy\omega)> 2R$ and $B_{2R}(a^i)\cap B_{2R}(a^j)= \emptyset$ whenever $i\neq j$ we have 
\begin{equation*}
\begin{split}
	\int_{\omega\setminus \bigcup_{i = 1}^mB_{R}(a^i)} \varphi f_k\; \d x
	& 
	\leq \|\varphi\|_{L^\infty(\omega)}\int_{\supp\varphi\setminus \bigcup_{i = 1}^mB_{R}(a^i)} f_k\; \d x, 
\end{split}
\end{equation*}
where $f_k = V_kI_\mu[e^{\lambda u_k}\chi_\Omega]e^{\lambda u_k}$. Moreover, since $u_k\to -\infty$ uniformly on $\supp\varphi\setminus \bigcup_{i = 1}^mB_{R}(a^i)$ we have 
\begin{equation*}
\begin{split}
	\int_{\supp\varphi\setminus \bigcup_{i = 1}^mB_{R}(a^i)} f_k\; \d x
	& 
	\leq C\|V_k\|_{L^\infty(\omega)}\|e^{u_k}\|_{L^1(\Omega)}^\lambda\|e^{u_k}\|_{L^1(\supp\varphi\setminus \bigcup_{i = 1}^mB_{R}(a^i))}^\lambda\\
	& = \circ(1). 
\end{split}
\end{equation*}
Therefore, for any such $\varphi$ and any such $R$, 
\begin{equation*}
	\int_\omega \varphi f_k
	= \int_{\bigcup_{i = 1}^m B_R(a_i)}\varphi f_k + \circ(1)
	= \sum_{i =1 }^m\int_{B_R(a^i)}\varphi f_k + \circ(1). 
\end{equation*}
In view of the uniform continuity of $\varphi$ on $\omega$, for any $\epsilon>0$ we may choose $R \in (0,1)$ 
sufficiently small such that $|\varphi(a^i) - \varphi(x)|< \epsilon$ whenever $i\in \{1, \ldots, m\}$ and $x\in B_R(a^i)$. For any such $R$ and for any $i\in\{1, \ldots, m\}$ we have
\begin{equation*}
	\abs{\int_{B_R(a^i)}\varphi f_k - \varphi(a^i)\int_{B_R(a^i)}f_k}
	\leq \epsilon \|f_k\|_{L^1(\omega)}
	\leq C\epsilon. 
\end{equation*}
For each $i\in \{1, \ldots, m\}$, Proposition \ref{prop:energy_quantization} guarantees the existence of a positive integer $N_i$ such that
\begin{equation*}
	\int_{B_R(a^i)}f_k
	= 8\pi N_i + \circ(1), 
\end{equation*}
so we conclude that 
\begin{equation*}
	\lim_k\int_\omega\varphi f_k
	= 8\pi \sum_{i = 1}^m N_i\varphi(a^i). 
\end{equation*}
\end{proof}
The remainder of this section is devoted to the proof of Proposition \ref{prop:energy_quantization}. The strategy is to employ a ``bubble selection process'' near the origin whereby a maximal number of functions approximately of the form \eqref{eq:rescaled_bubble} whose centers of symmetry are converging to the origin are selected (each such function is referred to as a ``bubble''). As suggested by the second equality in \eqref{eq:asserted_energies}, each bubble contributes $8\pi$ to the limit in \eqref{eq:BR_energy_limit}. Moreover, for each bubble, the contribution of $8\pi$ comes entirely from the spatial region in the immediate vicinity of the center of symmetry of the bubble. The spatial regions corresponding to distinct bubbles are disjoint and there is no nonlocal interaction among distinct bubbles that contributes to the limit in \eqref{eq:BR_energy_limit}. The selection process, together with the verification that each bubble contributes $8\pi$ to the limit in \eqref{eq:BR_energy_limit} and the fact that there are not nonlocal interactions at the $L^1$ level between distinct bubbles in carried out in Lemma \ref{lemma:N_bubble_selection}. In Lemma \ref{lemma:neck_energy_vanish} it is shown that there is no contribution to the limit in \eqref{eq:BR_energy_limit} coming from regions outside the regions in the immediate vicinities of the centers of symmetry of the bubbles. \\

The following lemma is the analog of Lemma 1 in \cite{LiShafrir1994}. It guarantees that the limiting coefficient function $V$ in Theorem \ref{theorem:quantization} cannot vanish at a blow-up point and that each blow-up point carries at least $8\pi$ in energy. 
\begin{lemma}
\label{lemma:blow_up_minimal_energy}
Under the hypotheses of Proposition \ref{prop:energy_quantization}, the inequalities $V(0)>0$ and $\alpha \geq 8\pi$ both hold. 
\end{lemma}
\begin{proof}
Let $(x_k)_{k = 1}^\infty\subset B_R$ with $u_k(x_k)= \max_{\overline B_R}u_k$. Assumptions \eqref{eq:max_to_infty} and \eqref{eq:uniformly_down_away_from_origin} guarantee that $u_k(x_k)\to\infty$ and $x_k\to 0$ respectively. Setting $\delta_k = e^{-u_k(x_k)/2}$ and defining 
\begin{equation*}
	v_k(y)
	= u_k(x_k + \delta_k y) + 2\log \delta_k
	\qquad \text{ for }y\in \Omega_k:= \frac{\Omega - x_k}{\delta_k} 
\end{equation*}	
we find that
\begin{equation*}
\begin{cases}
	-\lap v_k = H_kI_\mu[e^{\lambda v_k}\chi_{\Omega_k}]e^{\lambda v_k} & \text{ in } \Sigma_k\\
	v_k\leq v_k(0)= 0  & \text{ in }\Sigma_k\\
	\int_{\Omega_k}e^{v_k}\leq c_0, 
\end{cases}
\end{equation*}
where $H_k(y) = V_k(x_k + \delta_k y)$ and $\Sigma_k = \frac{B_R- x_k}{\delta_k} 
$. In view of the containment $B_{R/(2\delta_k)}\subset \Sigma_k$ (which holds whenever $k$ is sufficiently large), for any $\rho\gg1$, the sequence $(v_k)_{k = 1}^\infty$ is well-defined in $B_\rho$ whenever $k$ is sufficiently large. 
Fixing $\rho\gg 1$ we have $H_k\to V(0)$ uniformly on $B_\rho$ and Theorem \ref{theorem:BM_alternative} 
guarantees that (along a subsequence) $(v_k)_{k = 1}^\infty$ is bounded in $L^\infty_{\loc}(B_\rho)$. By standard elliptic estimates there is $\alpha\in (0, 1)$ for which $(v_k)_{k = 1}^\infty$ is bounded in $C^{1, \alpha}(B_\rho)$. For any sequence $(\rho_\ell)_{\ell = 1}^\infty\subset(0, \infty)$ for which $\rho_\ell\to\infty$ we repeat this argument with $\rho$ replaced by $\rho_\ell$ and then employ a standard diagonal subsequence argument to find that there is $v\in C^{1, \alpha}(\bb R^2)$ that satisfies 
\begin{equation}
\label{eq:limit_function_properties}
\begin{cases}
	-\lap v = V(0)I_\mu[e^{\lambda v}]e^{\lambda v}& \text{ in }\bb R^2\\
	v\leq v(0) = 0 & \text{ in }\bb R^2\\
	\|e^v\|_{L^1(\bb R^2)}\leq c_0 
\end{cases}
\end{equation} 
and a subsequence of $(v_k)_{k = 1}^\infty$ along which $v_k\to v$ in $C^{1, \alpha}_{\loc}(\bb R^2)$. The verification that $v$ satisfies the PDE in \eqref{eq:limit_function_properties} is similar to the verification that the PDE in \eqref{eq:sup+inf_lemma_limiting_problem} is satisfied, as carried out on page \pageref{proof:verification}. From \eqref{eq:limit_function_properties} and the assumption that $V_k\geq 0$ for all $k$, we see that $V(0)>0$. 
To see that $\alpha \geq 8\pi$, apply Theorem \ref{oldtheorem:classification} of Appendix \ref{s:appendix} to the function $v + \frac{1}{2\lambda}\log V(0)$ to obtain
\begin{equation}
\label{eq:centered_v_gamma}
	v(y) = - 2\log(1 + \gamma^2|y|^2) 
	\quad \text{ with } \quad 
	\gamma = \left(\frac{\pi V(0)}{4(2 - \mu)}\right)^{\frac 1{4 - \mu}}. 
\end{equation}
In particular, the second item of \eqref{eq:asserted_energies} guarantees that $V(0)\int_{\bb R^2}I_\mu[e^{\lambda v}]e^{\lambda v} = 8\pi$. Let $\epsilon\in (0, \frac{8\pi}{V(0)})$ and choose $r\gg 1$ such that 
\begin{equation*}
	V(0)\int_{B_r}I_\mu[e^{\lambda v}\chi_{B_{2r}}]e^{\lambda v}
	> 8\pi - \epsilon. 
\end{equation*}
For any such $r$ and for $k$ sufficiently large so that $\Sigma_k\supset B_{2r}$ we have
\begin{equation*}
\begin{split}
	\int_{B_R}V_kI_\mu[e^{\lambda u_k}\chi_\Omega]e^{\lambda u_k}
	& = \int_{\Sigma_k}H_kI_\mu[e^{\lambda v_k}\chi_{\Omega_k}]e^{\lambda v_k}\\
	& \geq \int_{B_r}H_kI_\mu[e^{\lambda v_k}\chi_{B_{2r}}]e^{\lambda v_k}\\
	& = V(0)\int_{B_r}I_\mu[e^{\lambda v}\chi_{B_{2r}}]e^{\lambda v} + \circ(1)\\
	& \geq 8\pi - \epsilon +\circ(1). 
\end{split}
\end{equation*}
Since $\epsilon>0$ is arbitrary, letting $k\to\infty$ gives $\alpha\geq 8\pi$. 
\end{proof}
The following lemma specifies the bubble selection process. 
\begin{lemma}
\label{lemma:N_bubble_selection}
Let $\Omega\subset \bb R^2$ be a bounded domain for which $\overline B_R\subset\Omega$, let $\mu\in (0, 2)$ and let $\lambda$ be as in \eqref{eq:lambda}. Let $(V_k)_{k = 1}^\infty\subset C^0(\overline B_R)$ be a sequence of nonnegative functions satisfying $V_k\to V$ in $C^0(\overline B_R)$ for some nonnegative $V\in C^0(\overline B_R)$, and let $(\rho_k)_{k = 1}^\infty\subset (0, \infty)$ be any sequence for which $\rho_k\to\infty$. If $(u_k)_{k = 1}^\infty$ is a sequence of distributional solutions to \eqref{eq:BR_PDE_and_energy} for which both \eqref{eq:max_to_infty} and \eqref{eq:uniformly_down_away_from_origin} hold, then $V(0)>0$ and there exists $N\in\bb N$, there exists a collection of sequences $\{(x_k^{(j)})_{k = 1}^\infty:j = 0, \ldots, N -1\}$ in $B_R$, and there exists a subsequence of $(u_k)_{k = 1}^\infty$ along which all of the following hold with $\delta_k^{(j)} = e^{-u_k(x_k^{(j)})/2}$ and $r_k^{(j)} = \rho_k \delta_k^{(j)}$: 
\begin{enumerate}[label = {\bf \arabic*.}, ref = {\bf \arabic*}, wide = 0pt]
	\item \label{item:xkj_local_maximizer} For every $j\in \{0, \ldots, N - 1\}$, we have $r_k^{(j)}= \circ(1)$, $x_k^{(j)}\to 0$ and 
	\begin{equation}
	\label{eq:xkj_blowup}
		u_k(x_k^{(j)}) = \max_{\overline B(x_k^{(j)}, r_k^{(j)})}u_k
	\to \infty. 
	\end{equation} 
	\item \label{item:discreasing_about_center} For every $j\in \{0, \ldots, N - 1\}$ and every $x\in B(0, 16r_k^{(j)})\setminus B(0, \delta_k^{(j)})$ there holds
	\begin{equation}
	\label{eq:decreasing_about_xkj}
		\left.\frac{\d}{\d t}\right|_{t = 1}u_k(x_k^{(j)} + tx) < 0.
	\end{equation}
	\item \label{item:disjoint_energy_balls} For every pair of distinct indices $i, j\in \{0, \ldots, N - 1\}$, there holds
	\begin{equation}
	\label{eq:disjoint_energy_balls}
		B(x_k^{(j)}, 16r_k^{(j)})\cap B(x_k^{(i)}, 16r_k^{(i)}) = \emptyset
	\end{equation}
	and
	\begin{equation}
	\label{eq:j_relative_radius_vanishing}
		r_k^{(j)} = \circ(1)|x_k^{(i)} - x_k^{(j)}|
		\qquad \text{ whenever }j> i. 
	\end{equation}
	\item \label{item:no_energy_ring} For every $j\in \{0, \ldots, N - 1\}$ we have both
	\begin{equation}
	\label{eq:full_energy_region}
	\begin{split}
		\lim_{k\to\infty}& \int_{B(x_k^{(j)}, 16r_k^{(j)})}V_kI_\mu[e^{\lambda u_k}\chi_{B(x_k^{(j)}, 16r_k^{(j)})}]e^{\lambda u_k}\\
		= \; & \lim_{k\to\infty}\int_{B(x_k^{(j)}, r_k^{(j)})}V_kI_\mu[e^{\lambda u_k}\chi_{B(x_k^{(j)}, r_k^{(j)})}]e^{\lambda u_k}\\
		= \; & 8\pi
	\end{split}
	\end{equation}
	and 
	\begin{equation}
	\label{eq:no_energy_region}
		\lim_{k\to\infty} \int_{B(x_k^{(j)}, 16r_k^{(j)})}V_kI_\mu[e^{\lambda u_k}\chi_{\Omega\setminus B_{R/4}}]e^{\lambda u_k}\\
		= 0. 
	\end{equation}
	If in addition $i\in \{0, \ldots, N - 1\}\setminus\{j\}$ then 
	\begin{equation}
	\label{eq:no_close_range_interaction}
		\lim_{k\to\infty} \int_{B(x_k^{(j)}, 15r_k^{(j)})}V_kI_\mu[e^{\lambda u_k}\chi_{B(x_k^{(i)}, 15r_k^{(i)})}]e^{\lambda u_k}
		= 0.
	\end{equation}
	\item \label{item:termination_condition} There is a constant $C>0$ such that for every $k$ there holds
	\begin{equation}
	\label{eq:termination_condition}
		\max_{x\in \overline B_R}\{u_k(x) + 2\log\min_{0\leq j\leq N - 1}|x - x_k^{(j)}|\}\leq C. 
	\end{equation}
\end{enumerate}
\end{lemma}
\begin{remark}
Conditions \eqref{eq:no_energy_region} and \eqref{eq:no_close_range_interaction} are nonlocal conditions that have no analogs in the local setting (e.g., for problems of the form \eqref{eq:local_problem}). Condition \eqref{eq:no_energy_region} guarantees, among other things, that bubbles near $0\in B_R$ do not have sufficient ``long-range'' nonlocal interaction with bubbles near any other blow-up points $a^i\in \Omega\setminus B_R$ to contribute to the limit in \eqref{eq:BR_energy_limit}. Condition \eqref{eq:no_close_range_interaction} guarantees that the ``close-range'' nonlocal bubble interactions among distinct bubbles near $0\in B_R$ do not contribute to the limit in \eqref{eq:BR_energy_limit}. 
\end{remark}
\begin{proof}[Proof of Lemma \ref{lemma:N_bubble_selection}]
The inequality $V(0)>0$ follows as in Lemma \ref{lemma:blow_up_minimal_energy}. The idea of the rest of the proof is as follows. First, we select a sequence $(x_k^{(0)})_{k = 1}^\infty\subset B_R$ and a subsequence of $(u_k)_{k = 1}^\infty$  for which $x_k^{(0)}\to 0$ and for which all of items \ref{item:xkj_local_maximizer}, \ref{item:discreasing_about_center}, \ref{item:disjoint_energy_balls}, and \ref{item:no_energy_ring} hold for $N = 1$. 
Item \ref{item:termination_condition} is a termination condition. If it holds with $N = 1$ then we terminate the bubble selection process and set $N= 1$. Otherwise, the failure of item \ref{item:termination_condition} allows us to select a sequence $(x_k^{(1)})_{k = 1}^\infty\subset B_R$ and a further subsequence of $(u_k)_{k = 1}^\infty$ for which $x_k^{(1)}\to 0$ and for which all of items  \ref{item:xkj_local_maximizer}, \ref{item:discreasing_about_center}, \ref{item:disjoint_energy_balls}, and \ref{item:no_energy_ring} hold with $N = 2$. If after selecting such sequences, item \ref{item:termination_condition} holds with $N= 2$, we declare $N = 2$ and terminate the bubble selection process. Otherwise, we continue the selection process until we have found a positive integer $n$ for which all of  \ref{item:xkj_local_maximizer}, \ref{item:discreasing_about_center}, \ref{item:disjoint_energy_balls}, \ref{item:no_energy_ring}, and \ref{item:termination_condition} hold with $N = n$. To see that the selection process must terminate after finitely many steps, observe that on one hand \eqref{eq:full_energy_region} guarantees that each bubble contributes at least $8\pi$ to the limit in \eqref{eq:BR_energy_limit}, while on the other hand, the uniform integrability assumption in \eqref{eq:BR_PDE_and_energy} guarantees that $(\|f_k\|_{L^1(B_R)})_{k = 1}^\infty$ is bounded in $\bb R$, where $f_k = V_kI_\mu[e^{\lambda u_k}\chi_\Omega]e^{\lambda u_k}$.

In what follows we indicate the selection of the $0^{\text{th}}$ bubble, and for $n \geq 1$ we indicate the selection of the $n^{\text{th}}$ bubble in the event that bubbles $0, \ldots, n- 1$ have already been selected and the termination condition \eqref{eq:termination_condition} fails for $N = n$. 
\begin{enumerate}[label = {\bf Step \arabic*.}, ref = {Step \arabic*}, wide = 0pt]
	\item \label{item:selection_step1} In this step we select a sequence $(x_k^{(0)})_{k = 1}^\infty\subset B_R$ for which $x_k^{(0)}\to 0$ and we select a subsequence of $(u_k)_{k = 1}^\infty$ along which all of items \ref{item:xkj_local_maximizer},  \ref{item:discreasing_about_center}, \ref{item:disjoint_energy_balls}, and \ref{item:no_energy_ring} hold for $N = 1$. Let $x_k^{(0)}\in\overline B_R$ satisfy $u_k(x_k^{(0)}) = \max_{\overline B_R}u_k$. Assumptions \eqref{eq:max_to_infty} and \eqref{eq:uniformly_down_away_from_origin} guarantee that both $u_k(x_k^{(0)})\to\infty$ and $x_k^{(0)}\to 0$. Set $\delta_k^{(0)} = e^{-u_k(x_k^{(0)})/2}$ and define
\begin{equation*}
	v_k^{(0)}(y)
	= u_k(x_k^{(0)} + \delta_k^{(0)} y) + 2\log \delta_k^{(0)}
	\qquad \text{ for }y\in \Omega_k^{(0)}:= \frac{\Omega - x_k^{(0)}}{\delta_k^{(0)}}. 
\end{equation*}
Since $u_k$ satisfies \eqref{eq:BR_PDE_and_energy} and by the choice of $x_k^{(0)}$, we find that $v_k^{(0)}$ satisfies
\begin{equation*}
\begin{cases}
	-\lap v_k^{(0)} = H_k^{(0)} I_\mu[e^{\lambda v_k^{(0)}}\chi_{\Omega_k^{(0)}}]e^{\lambda v_k^{(0)}} & \text{ in } \Sigma_k^{(0)}\\
	v_k^{(0)}\leq 0 = v_k(0) & \text{ in }\Sigma_k^{(0)}\\
	\|e^{v_k^{(0)}}\|_{L^1(\Omega_k^{(0)})}\leq c_0, 
\end{cases}
\end{equation*}
where $H_k^{(0)}(y) = V_k(x_k^{(0)} + \delta_k^{(0)} y)$ and $\Sigma_k^{(0)} = (B_R - x_k^{(0)})/\delta_k^{(0)}$. A routine argument, similar to the one carried out in the proof of Lemma \ref{lemma:blow_up_minimal_energy} shows that for every $\rho\gg 1$ and for $k$ sufficiently large so that $B_{4\rho}\subset \Sigma_k^{(0)}$, (a subsequence of) the sequence $(-\lap v_k^{(0)})_{k = 1}^\infty$ is bounded in $L^\infty(B_\rho)$. By standard elliptic estimates and a diagonal subsequence argument we find that there is $\alpha\in (0, 1)$, there is $v\in C^{1, \alpha}(\bb R^2)$ satisfying
\begin{equation*}
\begin{cases}
	-\lap v = V(0) I_\mu[e^{\lambda v}]e^{\lambda v} & \text{ in }\bb R^2\\
	v\leq 0 = v(0) & \text{ in }\bb R^2\\
	\|e^v\|_{L^1(\bb R^2)}\leq c_0, 
\end{cases}
\end{equation*}
and there is a subsequence of $(u_k)_{k = 1}^\infty$ (whose members we continue to denote by $u_k$) along which both $r_k^{(0)}:= \rho_k\delta_k^{(0)} \to 0$ and 
\begin{equation}
\label{eq:j=0_C1alpha_convergence}
	\|v_k^{(0)}- v\|_{C^{1, \alpha}(B(0, 17\rho_k))}
	= \circ(1).
\end{equation} 
Applying Theorem \ref{oldtheorem:classification} of Appendix \ref{s:appendix} to $v + \frac 1{2\lambda}\log V(0)$, shows that $v$ is as in \eqref{eq:centered_v_gamma}. In particular, the second equality in \eqref{eq:asserted_energies} guarantees that $V(0)\|I_\mu[e^{\lambda v}]e^{\lambda v}\|_{L^1(\bb R^2)} = 8\pi$. Item \ref{item:xkj_local_maximizer} is clearly satisfied with $N = 1$. In the remainder of \ref{item:selection_step1} we verify that each of items \ref{item:discreasing_about_center}, \ref{item:disjoint_energy_balls}, and \ref{item:no_energy_ring} holds with $N = 1$. To see that item \ref{item:discreasing_about_center} holds for $N = 1$, observe that \eqref{eq:j=0_C1alpha_convergence} and the explicit form of $v$ in \eqref{eq:centered_v_gamma} give 
\begin{equation*}
	\left.\frac{\d}{\d t}\right|_{t = 1}u_k(x_k^{(0)} + tx)
	= \left.\frac{\d}{\d t}\right|_{t = 1}v_k^{(0)}\left(\frac{tx}{\delta_k^{(0)}}\right)
	< 0, 
\end{equation*}
whenever $x\in B(0, 16r_k^{(0)})\setminus B(0, \delta_k^{(0)})$. Item \ref{item:disjoint_energy_balls} holds trivially when $N = 1$. It remains to show that item \ref{item:no_energy_ring} holds with $N = 1$. To verify that \eqref{eq:full_energy_region} holds for $j = 0$, first observe that since $V_k\to V$ uniformly on $\overline B_R$ and by using the change of variable $x = x_k^{(0)} + \delta_k^{(0)} y$ we have
\begin{equation}
\label{eq:grab_8pi_decompose}
\begin{split}
	\int_{B(x_k^{(0)},16r_k^{(0)})}& V_kI_\mu[e^{\lambda u_k}\chi_{B(x_k^{(0)},16r_k^{(0)})}]e^{\lambda u_k} + \circ(1)\\
	& = V(0)\int_{B(x_k^{(0)}, 16r_k^{(0)})}I_\mu[e^{\lambda u_k}\chi_{B(x_k^{(0)},16r_k^{(0)})}]e^{\lambda u_k}\\
	& = V(0)\int_{B(0, 16\rho_k)}I_\mu[e^{\lambda v_k^{(0)}}\chi_{B(0, 16\rho_k)}]e^{\lambda v_k^{(0)}}. 
\end{split}
\end{equation}
Moreover, since $\rho_k\to\infty$ and in view of \eqref{eq:j=0_C1alpha_convergence} we have
\begin{equation}
\label{eq:grab_8pi_major}
\begin{split}
	\int_{B(0, 16\rho_k)}& I_\mu[e^{\lambda v_k^{(0)}}\chi_{B(0, 16\rho_k)}]e^{\lambda v_k^{(0)}}\\
	= & \; \int_{\bb R^2}I_\mu[e^{\lambda v}]e^{\lambda v} - \int_{\bb R^2\setminus B(0, 16\rho_k)}I_\mu[e^{\lambda v}\chi_{\bb R^2\setminus B(0, 16\rho_k)}]e^{\lambda v} \\
	& - 2\int_{\bb R^2\setminus B(0, 16\rho_k)}I_\mu[e^{\lambda v}\chi_{B(0, 16\rho_k)}]e^{\lambda v}\\ 
	& + \int_{B(0, 16\rho_k)}I_\mu[(e^{\lambda v_k^{(0)}} - e^{\lambda v})\chi_{B(0, 16\rho_k)}]e^{\lambda v_k^{(0)}}\\
	&  + \int_{B(0, 16\rho_k)}I_\mu[e^{\lambda v}\chi_{B(0, 16\rho_k)}]\left(e^{\lambda v_k^{(0)}} - e^{\lambda v}\right)\\
	= & \;  \frac{8\pi}{V(0)} + \circ(1). 
\end{split}
\end{equation}
Bringing \eqref{eq:grab_8pi_major} back to \eqref{eq:grab_8pi_decompose} verifies the limit in \eqref{eq:full_energy_region} corresponding to the larger balls $B(x_k^{(0)}, 16r_k^{(0)})$. The limit in \eqref{eq:full_energy_region} corresponding to the smaller balls $B(x_k^{(0)},r_k^{(0)})$ can be verified in a similar manner. Next we verify that \eqref{eq:no_energy_region} holds with $j = 0$. When $k$ is sufficiently large, for every $y\in B(0, 16\rho_k)$ have $\Omega_k^{(0)}\setminus\Sigma_k^{(0)}\subset \bb R^2\setminus B(y, R/(2\delta_k^{(0)}))$.  
In view of this containment, for any such $y$ an application of H\"older's inequality and the upper bound on $\|e^{v_k^{(0)}}\|_{L^1(\Omega_k^{(0)})}$ yields the pointwise estimate
\begin{equation*}
\begin{split}
	I_\mu[e^{\lambda v_k^{(0)}}\chi_{\Omega_k^{(0)}\setminus\Sigma_k^{(0)}}](y)
	& 
	\leq c_0^\lambda \left(\int_{\bb R^2\setminus B(y, R/(2\delta_k^{(0)}))}|y - \zeta|^{-4}\; \d \zeta\right)^{1- \lambda}\\
	& \leq C\left(\delta_k^{(0)}\right)^{\mu/2},  
\end{split}
\end{equation*}
and consequently
\begin{equation}
\label{eq:L4mu_far_away_estimate}
\begin{split}
	\|I_\mu[e^{\lambda v_k^{(0)}}\chi_{\Omega_k^{(0)}\setminus\Sigma_k^{(0)}}]\|_{L^{4/\mu}(B(0, 16\rho_k))}
	& 
	\leq C(r_k^{(0)})^{\mu/2} = \circ(1). 
\end{split}
\end{equation}
Now using \eqref{eq:uniformly_down_away_from_origin}, the change of variable $x= x_k^{(0)}+ \delta_k^{(0)} y$, H\"older's inequality, and estimate \eqref{eq:L4mu_far_away_estimate} we obtain
\begin{equation*}
\begin{split}
	\int_{B(x_k^{(0)},16r_k^{(0)})}& V_kI_\mu[e^{\lambda u_k}\chi_{\Omega\setminus B_{R/4}}]e^{\lambda u_k} + \circ(1)\\
	& \leq \|V_k\|_{L^\infty(\Omega)} \int_{B(x_k^{(0)},16r_k^{(0)})}I_\mu[e^{\lambda u_k}\chi_{\Omega\setminus B_R}]e^{\lambda u_k} \\
	& = \|V_k\|_{L^\infty(\Omega)} \int_{B(0, 16\rho_k)}I_\mu[e^{\lambda v_k^{(0)}}\chi_{\Omega_k^{(0)}\setminus\Sigma_k^{(0)}}]e^{\lambda v_k^{(0)}}\\
	& \leq \|V_k\|_{L^\infty(\Omega)}\|I_\mu[e^{\lambda v_k^{(0)}}\chi_{\Omega_k^{(0)}\setminus\Sigma_k^{(0)}}]\|_{L^{4/\mu}(B(0, 16\rho_k))}\|e^{v_k^{(0)}}\|_{L^1(\Omega_k^{(0)})}^\lambda\\
	& \leq C(r_k^{(0)})^{\mu/2} = \circ(1),  
\end{split}
\end{equation*}
which completes the verification of \eqref{eq:no_energy_region} for $j = 0$. 
	\item \label{item:selection_step2} Suppose that for some $n \geq 1$ and for each $j \in\{ 0, 1, \ldots, n - 1\}$ a sequence $(x_k^{(j)})_{k = 1}^\infty\subset B_R$ has been selected so that $x_k^{(j)}\to 0$, and a subsequence of $(u_k)_{k = 1}^\infty$ has been selected so that all of items \ref{item:xkj_local_maximizer}, \ref{item:discreasing_about_center}, \ref{item:disjoint_energy_balls}, and \ref{item:no_energy_ring} are satisfied with $N = n$, where $\delta_k^{(j)} = e^{-u_k(x_k^{(j)})/2}$ and $r_k^{(j)} = \rho_k\delta_k^{(j)}$. Suppose further that the termination condition in item \eqref{eq:termination_condition} fails with $N = n$. 
	We will construct a sequence $(x_k^{(n)})_{k = 1}^\infty\subset B_R$ for which $x_k^{(n)}\to 0$ and we will construct a further subsequence of $(u_k)_{k = 1}^\infty$ along which all of items \ref{item:xkj_local_maximizer}, \ref{item:discreasing_about_center}, \ref{item:disjoint_energy_balls},  and \ref{item:no_energy_ring} are satisfied with $N = n + 1$, where $\delta_k^{(n)} = e^{-u_k(x_k^{(n)})/2}$ and $r_k^{(n)} = \rho_k\delta_k^{(n)}$. Let $z_k^{(n)}\in \overline B_R \setminus\{x_k^{(0)}, \ldots, x_k^{(n - 1)}\}$ satisfy 
\begin{equation}
\label{eq:zkn}
	u_k(z_k^{(n)}) + 2\log\min_{0\leq j\leq n - 1}|z_k^{(n)} - x_k^{(j)}|
	= M_k^{(n)}, 
\end{equation}
where 
\begin{equation}
\label{eq:Mkl}
	M_k^{(n)}= \max_{x\in \overline B_R}\{u_k(x) + 2\log\min_{0\leq j\leq n - 1}|x - x_k^{(j)}|\} \to\infty. 
\end{equation}
Note that the sequence $(z_k^{(n)})_{k = 1}^\infty$ may not satisfy \eqref{eq:xkj_blowup} so in what follows we carefully choose $x_k^{(n)}$ near $z_k^{(n)}$ so that \eqref{eq:xkj_blowup} is satisfied. Setting $\sigma_k^{(n)} = e^{-u_k(z_k^{(n)})/2}$, the condition $M_k^{(n)}\to \infty$ is equivalent to the condition
\begin{equation*}
	\min_{0\leq j\leq n - 1}\frac{|z_k^{(n)} - x_k^{(j)}|}{\sigma_k^{(n)}} 
	\to \infty.
\end{equation*}
In particular, $u_k(z_k^{(n)})\to \infty$, so assumption \eqref{eq:uniformly_down_away_from_origin} forces $z_k^{(n)}\to 0$. Define 
\begin{equation*}
	\tilde v_k^{(n)}(y)
	= u_k\left(z_k^{(n)} + \sigma_k^{(n)} y\right) + 2\log\sigma_k^{(n)}
	\qquad \text{ for }y\in \tilde \Omega_k^{(n)} := \frac{\Omega - z_k^{(n)}}{\sigma_k^{(n)}}. 
\end{equation*}	
Observe that if $|y|< \min_{0\leq j\leq n - 1}\frac{|z_k^{(n)} - x_k^{(j)}|}{2\sigma_k^{(n)}}$ then for every $j\in\{0, \ldots, n - 1\}$ we have $|z_k^{(n)} + \sigma_k^{(n)} y - x_k^{(j)}|
	> |z_k^{(n)} - x_k^{(j)}|/2. $
For any such $y$, choosing $j\in \{0, \ldots, n - 1\}$ for which the minimum in \eqref{eq:Mkl} is attained we have 
\begin{equation*}
\begin{split}
	\tilde v_k^{(n)}(y)
	& 
	\leq M_k^{(n)} + 2\log\frac{2\sigma_k^{(n)} }{|z_k^{(n)} - x_k^{(j)}|}\\
	& = 2\log 2. 
\end{split}
\end{equation*}
Thus, $\tilde v_k^{(n)}$ satisfies
\begin{equation*}
\begin{cases}
	-\lap \tilde v_k^{(n)} = \tilde H_k^{(n)}I_\mu[e^{\lambda \tilde v_k^{(n)}}\chi_{\tilde \Omega_k^{(n)}}]e^{\lambda \tilde v_k^{(n)}} & \text{ in }\frac{B_R - z_k^{(n)}}{\sigma_k^{(n)}}\\
	\tilde v_k^{(n)}(0) = 0\\
	\tilde v_k^{(n)}(y) \leq 2\log 2
	& \text{ for }|y|< \min_{0\leq j\leq n - 1}\frac{|z_k^{(n)} - x_k^{(j)}|}{2\sigma_k^{(n)}}
	\\
	\|e^{\tilde v_k^{(n)}}\|_{L^1(\tilde \Omega_k^{(n)})}\leq c_0,  
\end{cases}
\end{equation*}
where $\tilde H_k^{(n)}(y) = V_k(z_k^{(n)} + \sigma_k^{(n)}y)$. By standard elliptic estimates and a diagonal subsequence argument, there is $\tilde v\in C^{1, \alpha}(\bb R^2)$ satisifying
\begin{equation*}
\begin{cases}
	-\lap \tilde v = V(0)I_\mu[e^{\lambda \tilde v}]e^{\lambda \tilde v}& \text{ in }\bb R^2\\
	\tilde v(0) = 0\\
	\tilde v\leq 2\log 2 & \text{ in }\bb R^2\\
	\|e^{\tilde v}\|_{L^1(\bb R^2)}\leq c_0, 
\end{cases}
\end{equation*}
and there is a further subsequence of $(u_k)_{k = 1}^\infty$ 
for which both 
\begin{equation}
\label{eq:rhokm_grow_slowly}
	\frac{\rho_k\sigma_k^{(n)}}{\min_{0\leq j\leq n- 1}|z_k^{(n)} - x_k^{(j)}|}
	\to 0
\end{equation}
and
\begin{equation}
\label{eq:vkl_to_tildev}
	\|\tilde v_k^{(n)} - \tilde v\|_{C^{1, \alpha}(B(0, 17\rho_k))}\to 0.
\end{equation}	
Applying Theorem \ref{oldtheorem:classification} to $\tilde v + \frac{1}{2\lambda}\log V(0)$ we find that there is $(y_0, \delta)\in \bb R^2\times (0, \infty)$ for which 
\begin{equation*}
	\tilde v(y) = 2\log\frac{\delta}{\gamma(1 + \delta^2|y - y_0|^2)}, 
\end{equation*}
where $\gamma$ is as in \eqref{eq:centered_v_gamma}. In particular, the second equality in line \eqref{eq:asserted_energies} of the appendix implies that
\begin{equation}
\label{eq:v_tilde_energy}
	V(0)\int_{\bb R^2}I_\mu[e^{\lambda \tilde v}]e^{\lambda \tilde v}
	= 8\pi. 
\end{equation}
Since $\tilde v(0) = 0$ and $\tilde v\leq 2\log 2$ we find that
\begin{equation}
\label{eq:estimate_maximizer_location}
	\delta\leq 2\gamma
	\qquad \text{ and }\qquad
	|y_0|\leq (\delta\gamma)^{-1/2}.
\end{equation}
Let $y_k^{(n)}\in B(0, 16\rho_k)$ satisfy 
\begin{equation}
\label{eq:ykm_maximizer}
	\tilde v_k^{(n)}(y_0 + y_k^{(n)}) 
	= \max_{y\in \overline B(0, 17\rho_k)}\tilde v_k^{(n)}(y_0 + y)
\end{equation}
and define $x_k^{(n)} = z_k^{(n)} + \sigma_k^{(n)}(y_0 + y_k^{(n)})$. In view of \eqref{eq:estimate_maximizer_location} and \eqref{eq:vkl_to_tildev} we find that $y_k^{(n)}\to 0$. Setting $A = 1 + \max_{\overline B} \tilde v- \min_{\overline B}\tilde v$, where $B = B(0, 1 + (\delta\gamma)^{-1/2})$, in view of  \eqref{eq:estimate_maximizer_location} and \eqref{eq:vkl_to_tildev} we have
\begin{equation}
\label{eq:last_max_is_comparable}
	u_k(z_k^{(n)})\leq u_k(x_k^{(n)})\leq u_k(z_k^{(n)}) + A.  
\end{equation}
In particular, $u_k(x_k^{(n)})\to\infty$ so assumption \eqref{eq:uniformly_down_away_from_origin} implies that $x_k^{(n)}\to 0$. Setting $\delta_k^{(n)} = e^{-u_k(x_k^{(n)})/2}$, \eqref{eq:last_max_is_comparable} gives
\begin{equation}
\label{eq:deltakm_mukm_comparable}
	e^{-A/2}\sigma_k^{(n)}\leq \delta_k^{(n)}\leq \sigma_k^{(n)}. 
\end{equation} 
Setting $r_k^{(n)} = \rho_k\delta_k^{(n)}$, line \eqref{eq:rhokm_grow_slowly} guarantees that $r_k^{(n)} = \circ(1)$. Moreover, since $y_k^{(n)}$ satisfies \eqref{eq:ykm_maximizer}, $x_k^{(n)}$ satisfies 
\begin{equation*}
	u_k(x_k^{(n)}) = \max_{\overline B(x_k^{(n)}, r_k^{(n)})} u_k
	\to\infty.  
\end{equation*}
Thus, 
item \ref{item:xkj_local_maximizer} holds with $N = n + 1$. 
Defining
\begin{equation*}
	v_k^{(n)}(y)
	= u_k(x_k^{(n)} + \delta_k^{(n)}y) + 2\log\delta_k^{(n)}
	\qquad \text{ for } y\in \Omega_k^{(n)} := \frac{\Omega - x_k^{(n)}}{\delta_k^{(n)}},  
\end{equation*}	
we conclude that along a subsequence
\begin{equation}
\label{eq:vkm_convergence}
	\|v_k^{(n)} - \left(\tilde v(y_0 + b\;\cdot )+ 2\log b\right)\|_{C^{1, \alpha}(B(0, \frac{33\rho_k}2))} = \circ(1), 
\end{equation}
where $b = \lim_k\frac{\delta_k^{(n)}}{\sigma_k^{(n)}}\in [e^{-A/2}, 1]$. In the remainder of \ref{item:selection_step2} we verify that items \ref{item:discreasing_about_center}, \ref{item:disjoint_energy_balls}, and \ref{item:no_energy_ring} hold with $N = n + 1$. To verify that item \ref{item:discreasing_about_center} holds with $N = n + 1$, we only need to verify that \eqref{eq:decreasing_about_xkj} holds with $j =n$. From \eqref{eq:vkm_convergence} and from the explicit form of $\tilde v$, for any $x\in B(0, 16r_k^{(n)})\setminus B(0, \delta_k^{(n)})$ we have
\begin{equation*}
	\left.\frac{\d}{\d t}\right|_{t = 1}u_k(x_k^{(n)} + tx)
	= \left.\frac{\d}{\d t}\right|_{t = 1}v_k^{(n)}\left(\frac{tx}{\delta_k^{(n)}}\right)
	< 0. 
\end{equation*}
To verify that item \ref{item:disjoint_energy_balls} holds for $N = n+ 1$ we only need to verify that \eqref{eq:disjoint_energy_balls} and \eqref{eq:j_relative_radius_vanishing} hold with $j = n$ and $i\in \{0, \ldots, n -1\}$. With $j$ and $i$ as such, \eqref{eq:disjoint_energy_balls} follows immediately from the fact that item \ref{item:discreasing_about_center} holds with $N = n+ 1$. The relation 
\begin{equation}
\label{eq:nth_relative_radius_vanishes}
	r_k^{(n)} = \circ(1)|x_k^{(n)} - x_k^{(i)}|
	\qquad \text{ whenever }i\in \{0, \ldots, n - 1\}
\end{equation}
is a consequence of \eqref{eq:rhokm_grow_slowly}, \eqref{eq:deltakm_mukm_comparable}, and the estimate $2|x_k^{(n)} - x_k^{(i)}|\geq |z_k^{(n)} - x_k^{(i)}|$. 
Finally, we verify that item \ref{item:no_energy_ring} holds with $N = n+ 1$. Using \eqref{eq:vkm_convergence}, the verification that \eqref{eq:full_energy_region} is satisfied with $j = n$ is similar to the verification that \eqref{eq:full_energy_region} is satisfied with $j = 0$ as carried out in \eqref{eq:grab_8pi_decompose}, \eqref{eq:grab_8pi_major} so we omit the details. 
Similarly, the verification of \eqref{eq:no_energy_region} for $j = n$ is the same as the verification of \eqref{eq:no_energy_region} for $j = 0$ as carried out above so we omit the details.
Next we verify that for every $i\in\{0, \ldots, n - 1\}$,  \eqref{eq:no_close_range_interaction} is satisfied with $j = n$.  Fix any such $i$ and set
\begin{equation*}
	\Gamma_k^{(i)} = \frac{B(x_k^{(i)}, 15r_k^{(i)}) - x_k^{(n)}}{\delta_k^{(n)}}. 
\end{equation*}
Since $B(x_k^{(i)}, 16r_k^{(i)})\cap B(x_k^{(n)}, 16r_k^{(n)}) = \emptyset$ 
for any $x\in B(x_k^{(n)}, 15r_k^{(n)})$ and any $z\in B(x_k^{(i)}, 15r_k^{(i)})$, there holds $|x- z|\geq |x_k^{(n)} - x_k^{(j)}|/16$. Consequently, for any $y\in B(0, 15\rho_k)$, the containment $\Gamma_k^{(i)}\subset \bb R^2\setminus B(y, |x_k^{(n)} - x_k^{(i)}|/(16\delta_k^{(n)}))$ holds and thus we obtain the pointwise estimate
\begin{equation}
\label{eq:Imu_close_range_estimate}
\begin{split}
	I_\mu[e^{\lambda v_k^{(n)}}\chi_{\Gamma_k^{(i)}}](y)
	& 
	\leq \|e^{v_k^{(n)}}\|_{L^1(\Omega_k^{(n)})}^\lambda\left(\int_{\bb R^2\setminus B(y, \frac{|x_k^{(n)} - x_k^{(i)}|}{16\delta_k^{(n)}})}|y - \zeta|^{-4}\; \d \zeta\right)^{1 - \lambda}\\
	& \leq C\left(\frac{\delta_k^{(n)}}{|x_k^{(n)} - x_k^{(j)}|}\right)^{\mu/2}. 
\end{split}
\end{equation}
Using the change of variable $x = x_k^{(n)} + \delta_k^{(n)}y$ together with H\"older's inequality and estimate \eqref{eq:Imu_close_range_estimate} we obtain 
\begin{equation}
\label{eq:no_m_i_close_range_interaction}
\begin{split}
	\int_{B(x_k^{(n)}, 15r_k^{(n)})}& I_\mu[e^{\lambda u_k}\chi_{B(x_k^{(i)}, 15r_k^{(i)})}]e^{\lambda u_k}\\
	& = \int_{B(0, 15\rho_k)}I_\mu[e^{\lambda v_k^{(n)}}\chi_{\Gamma_k^{(i)}}]e^{\lambda v_k^{(n)}}\\
	& \leq c_0^\lambda \|I_\mu[e^{\lambda v_k^{(n)}}\chi_{\Gamma_k^{(i)}}]\|_{L^{4/\mu}(B(0, 15\rho_k))}\\
	& 
	\leq C\left(\frac{r_k^{(n)}}{|x_k^{(n)} - x_k^{(i)}|}\right)^{\frac\mu 2}\\
	& =\circ(1),  
\end{split}
\end{equation}
where the convergence to zero in the final estimate follows from \eqref{eq:nth_relative_radius_vanishes}. We also note that by the symmetry of $I_\mu$ estimate \eqref{eq:no_m_i_close_range_interaction} implies 
\begin{equation*}
\begin{split}
	\int_{B(x_k^{(i)}, 15r_k^{(i)})}& I_\mu[e^{\lambda u_k}\chi_{B(x_k^{(n)}, 15r_k^{(n)})}]e^{\lambda u_k}\\
	& = \int_{B(x_k^{(n)}, 15r_k^{(n)})} I_\mu[e^{\lambda u_k}\chi_{B(x_k^{(i)}, 15r_k^{(i)})}]e^{\lambda u_k}\\
	& 
	= \circ(1).   
\end{split}
\end{equation*}
\end{enumerate}
\end{proof}
\begin{lemma}[Spherical Harnack Inequality]
\label{lemma:spherical_harnack}
Let $\mu\in (0, 2)$, let $\lambda$ be as in \eqref{eq:lambda}, let $R>0$ and let $R_0\in (0, \frac R4]$. There is a universal constant $\beta\in (0, 1)$ and for every $c_0, c_1, c_2>0$ there is a constant $\bar C = \bar C(c_0, c_1, c_2)>0$ such that for all $V\in L^\infty(B_R\setminus B_{R_0})$ satisfying $\|V\|_{L^\infty(B_R\setminus B_{R_0})}\leq c_1$, for all bounded domains $\Omega\subset \bb R^2$ for which $B_R\subset \Omega$, and for all solutions $u$ to 
\begin{equation*}
	-\lap u = VI_\mu[e^{\lambda u}\chi_\Omega]e^{\lambda u}
	\qquad \text{ in }B_R\setminus \overline B_{R_0}
\end{equation*}
for which both $\|e^u\|_{L^1(\Omega)}\leq c_0$ and 
\begin{equation}
\label{eq:implies_spherical_harnack}
	u(x) + 2\log|x| \leq c_2
	\qquad \text{ in }B_R\setminus B_{R_0}, 
\end{equation}
there holds
\begin{equation}
\label{eq:spherical_harnack}
	\sup_{\bdy B_r}u
	\leq \bar C + \beta\inf_{\bdy B_r}u - 2(1 - \beta)\log r
	\qquad \text{ for all }r\in [2R_0, \frac R2]. 
\end{equation}
We emphasize that $\beta$ and $\bar C$ are independent of $R$, $R_0$ and $\Omega$. 
\end{lemma}
\begin{proof}
For $r\in[2R_0, \frac R2]$ the function $\tilde u(y) = u(ry) + 2\log r$ is well-defined for $y\in r^{-1}\Omega$ and satisfies 
\begin{equation}
\label{eq:spherical_harnack_rescaled_PDE}
\begin{cases}
	-\lap \tilde u = \tilde V I_\mu[e^{\lambda\tilde u}\chi_{r^{-1}\Omega}]e^{\lambda \tilde u}
	& \text{ in }B_2\setminus\overline B_{1/2}\\
	\|e^{\tilde u}\|_{L^1(r^{-1}\Omega)}\leq c_0,  
\end{cases} 
\end{equation}
where $\tilde V(y) = V(ry)$. 
From \eqref{eq:implies_spherical_harnack}, for every $y\in B_2\setminus \overline B_{1/2}$ we have $\tilde u(y) \leq c_2 + 2\log 2$ and therefore 
$e^{\lambda\tilde u}\leq 4e^{c_2}$. Combining this estimate with the $L^\infty$ assumption on $V$ and the HLS inequality we find that
\begin{equation}
\label{eq:spherical_harnack_source_integrability}
\begin{split}
	\|\tilde VI_\mu[e^{\lambda\tilde u}\chi_{r^{-1}\Omega}]e^{\lambda \tilde u}\|_{L^{\frac 4\mu}(B_2\setminus B_{1/2})}
	& 
	\leq 4\mc Hc_1e^{c_2}\|e^{\lambda\tilde u}\|_{L^{1/\lambda}(r^{-1}\Omega)}\\
	& 
	\leq 4\mc Hc_1e^{c_2} c_0^\lambda. 
\end{split}
\end{equation}
Let $w$ be the solution of 
\begin{equation*}
\begin{cases}
	-\lap w = \tilde V I_\mu[e^{\lambda \tilde u}\chi_{r^{-1}\Omega}]e^{\lambda \tilde u} & \text{ in }B_2\setminus \overline B_{1/2}\\
	w = 0 & \text{ on } \bdy(B_2\setminus B_{1/2}). 
\end{cases}
\end{equation*}
Estimate \eqref{eq:spherical_harnack_source_integrability} gives $\|\lap w\|_{L^{4/\mu}(B_2\setminus B_{1/2})}\leq C(\mu, c_0, c_1, c_2)$, so since $\frac 4\mu> 2$ standard elliptic estimates guarantee that $w\in L^\infty(B_2\setminus B_{1/2})$ with 
\begin{equation}
\label{eq:spherical_harnack_w_uniform_bound}
	\|w\|_{L^\infty(B_2\setminus B_{1/2})}
	\leq C(\mu, c_0, c_1, c_2). 
\end{equation}
Now, the function $h = w - \tilde u$ is harmonic in $B_2\setminus \overline B_{1/2}$ and satisfies $h+ C>0$ in $B_2\setminus \overline B_{1/2}$ for some $C = C(\mu, c_0, c_1, c_2)>0$, so Harnack's inequality ensures the existence of a universal constant $\beta\in (0, 1)$ for which 
\begin{equation}
\label{eq:applied_harnack}
	\sup_{\bdy B_1}(h + C)
	\leq \frac 1\beta\inf_{\bdy B_1}(h + C). 
\end{equation}
Unraveling the definitions of $h$ and $\tilde u$ in estimate \eqref{eq:applied_harnack} and using estimate \eqref{eq:spherical_harnack_w_uniform_bound} yields estimate \eqref{eq:spherical_harnack}. 
\end{proof}
The following lemma is a consequence of the $\sup + \inf$ inequality of Theorem \ref{theorem:sup_inf_inequality}. We remind the reader of the notation $\Lambda_{a,b}(\overline B_R)$ defined in \eqref{eq:lambda_ab}. 
\begin{lemma}
\label{lemma:rescaled_sup+inf_inequality}
Let $\Omega\subset \bb R^2$ be a bounded domain for which  $\overline B_R\subset\Omega$, let $\mu\in (0, 2)$, and let $\lambda$ be as in \eqref{eq:lambda}. Suppose $0< a\leq b< \infty$ and let $\Lambda \subset \Lambda_{a,b}(\overline B_R)$ be equicontinuous at each point of $\overline B_R$. For each $c_0>0$ and each $C_1>1$ there is a constant $C_2= C_2(\Lambda, R, c_0, C_1)> 0$ such that for all $V\in \Lambda$, all solutions $u$ to 
\begin{equation*}
\begin{cases}
	-\lap u = VI_\mu[e^{\lambda u}\chi_\Omega]e^{\lambda u} &  \text{ in }B_R\\
	\|e^u\|_{L^1(\Omega)}\leq c_0, 
\end{cases}
\end{equation*}
and any $r\in (0, R)$, there holds
\begin{equation*}
	u(0) + C_1\inf_{B_r}u + 2(1 + C_1)\log r \leq C_2. 
\end{equation*}
We emphasize that $C_2$ is independent of $\Omega$. 
\end{lemma}
\begin{proof}
For $r\in (0, R)$, the function $\tilde u(y) = u(ry) + 2\log r $ defined for $y\in r^{-1}\Omega$ satisfies
\begin{equation*}
\begin{cases}
	-\lap \tilde u = \tilde VI_\mu[e^{\lambda \tilde u}\chi_{r^{-1}\Omega}]e^{\lambda \tilde u}
	& \text{ in } B_1\subset B_{R/r}\\
	\|e^{\tilde u}\|_{L^2(r^{-1}\Omega)} \leq c_0, 
\end{cases}
\end{equation*}
where $\tilde V(y) = V(ry)$. Applying 
Theorem \ref{theorem:sup_inf_inequality} to $\tilde u$ with $K = \{0\}$ and $\omega = B_1$ yields $\tilde u(0) + C_1\inf_{B_1}\tilde u\leq C_2$. The asserted inequality follows. 
\end{proof}
The following lemma implies that with $N$, $\{(x_k^{(j)}, r_k^{(j)}):j\in \{0, \ldots, N - 1\}\}$ as in Lemma \ref{lemma:N_bubble_selection}, the energy of $V_kI_\mu[e^{\lambda u_k}\chi_\Omega]e^{\lambda u_k}$ in $B_R\setminus\bigcup_{j = 0}^{N - 1}B(x_k^{(j)}, r_k^{(j)})$ is negligible in the sense that it does not contribute to the limit in \eqref{eq:BR_energy_limit}. 
\begin{lemma}
\label{lemma:neck_energy_vanish}
Let $\Omega\subset \bb R^2$ be a bounded domain for which $\overline B_R\subset \Omega$, let $\mu\in (0, 2)$, and let $\lambda$ be as in \eqref{eq:lambda}. Let $(V_k)_{k = 1}^\infty\subset C^0(\overline B_R)$ be a sequence of positive functions for which $V_k\to V$ in $C^0(\overline B_R)$ for some positive function $V\in C^0(\overline B_R)$ and suppose $(u_k)_{k = 1}^\infty$ is a sequence of solutions to \eqref{eq:BR_PDE_and_energy} for which \eqref{eq:uniformly_down_away_from_origin} holds. Suppose $N\in \bb N$ and $\{(x_k^{(j)}, r_k^{(j)})_{k = 1}^\infty: j\in \{0, \ldots, N - 1\}\}$ is a collection of sequences in $B_R\times (0, \infty)$ for which $u_k(x_k^{(j)})\to\infty$, and for which the following four items hold: 
\begin{enumerate}
	\item \label{item:rk_rate} For every $j\in \{0, \ldots, N - 1\}$ 
	\begin{equation}
	\label{eq:rk_rate}
		\frac{r_k^{(j)}}{\delta_k^{(j)}}\to \infty, 
	\end{equation}
	where $\delta_k^{(j)} = e^{-u_k(x_k^{(j)})/2}$.
	\item \label{item:the_balls_disjoint} If $i, j\in \{0,  \ldots, N - 1\}$ are distinct then  
	\begin{equation}
	\label{eq:the_balls_disjoint}
		B(x_k^{(j)}, r_k^{(j)})\cap B(x_k^{(i)}, r_k^{(i)})=\emptyset, 
	\end{equation}
	and there is $\ell\in \{i, j\}$ for which 
	\begin{equation}
	\label{eq:exists_vanishing_relative_radius}
		r_k^{(\ell)} = \circ(1)|x_k^{(i)} - x_k^{(j)}|. 
	\end{equation}
	\item \label{item:termination_condition_satisfied} There exists a constant $C$ 	for which
	\begin{equation}
	\label{eq:termination_condition_satisfied}
		\max_{x\in \overline B_R\setminus \bigcup_{j = 0}^{N - 1}B(x_k^{(j)}, r_k^{(j)})}\{u_k(x) + 2\log\min_{0\leq i\leq N- 1}|x - x_k^{(i)}|\}\leq C. 
	\end{equation}
	\item \label{item:beta_masses} There exist positive constants $\beta_0, \ldots, \beta_{N -1}$ such that both
	\begin{enumerate}
		\item \label{subitem:self_energy_limits} for any $j\in \{0, \ldots, N -1\}$ we have both
\begin{equation}
\label{eq:self_energy_limits}
\begin{split}
	\lim_{k\to\infty}& \int_{B(x_k^{(j)}, 16r_k^{(j)})}V_kI_\mu[e^{\lambda u_k}\chi_{B(x_k^{(j)}, 16r_k^{(j)})}]e^{\lambda u_k}\\
	= & \; \lim_{k\to\infty}\int_{B(x_k^{(j)}, r_k^{(j)})}V_kI_\mu[e^{\lambda u_k}\chi_{B(x_k^{(j)}, r_k^{(j)})}]e^{\lambda u_k}\\
	= & \; \beta_j,
\end{split}
\end{equation}
	and 
	\begin{equation}
	\label{eq:no_long_range}
	\lim_{k\to\infty} \int_{B(x_k^{(j)}, 16r_k^{(j)})}V_kI_\mu[e^{\lambda u_k}\chi_{\Omega\setminus B_{R/4}}]e^{\lambda u_k}
	= 0.
\end{equation}
	\item \label{subitem:no_close_interactions} for every pair of distinct indices $i,j\in\{0, \ldots, N - 1\}$,  
	\begin{equation}
	\label{eq:assumption_no_close_range}
	\lim_{k\to\infty} \int_{B(x_k^{(j)}, 15r_k^{(j)})}V_kI_\mu[e^{\lambda u_k}\chi_{B(x_k^{(i)}, 15r_k^{(i)})}]e^{\lambda u_k}
	= 0,
\end{equation}
	\end{enumerate} 
\end{enumerate} 
Then 
\begin{equation}
\label{eq:full_ BR_limit_neck_argument}
	\int_{B_R}V_kI_\mu[e^{\lambda u_k}\chi_\Omega]e^{\lambda u_k}\to \sum_{j = 0}^{N - 1} \beta_j. 
\end{equation}
\end{lemma}
\begin{proof}
To simplify the notation we set $f_k = V_kI_\mu[e^{\lambda u_k}\chi_\Omega]e^{\lambda u_k}$ and to clarify the presentation we refer to the points $x_k^{(i)}$ as \emph{centers}, the numbers $\delta_k^{(i)}$ as \emph{scales}, the numbers $r_k^{(i)}$ as \emph{radii} and the numbers $\beta^{(j)}$ as \emph{masses}. We proceed by induction on $N$. 
\begin{enumerate}[label = {\bf Step \arabic*.}, ref = {Step \arabic*}, wide = 0pt]
	\item \label{item:assume_N1}
	Assume $N = 1$.
	In the confines of \ref{item:assume_N1}, we use the simplified notation $x_k = x_k^{(0)}$, $\delta_k = \delta_k^{(0)}$, etc. We assume without losing generality that $x_k= 0$ for all $k$. Thus assumption \eqref{eq:termination_condition_satisfied} becomes 
	\begin{equation}
	\label{eq:step1_upper_bound}
		\max_{x\in \overline B_R\setminus B_{r_k}}\{u_k(x) + 2\log|x|\}\leq C. 
	\end{equation}
	Writing $\chi_\Omega = \chi_{B_{2r_k}} + \chi_{B_{R/4}\setminus B_{2r_k}} + \chi_{\Omega\setminus B_{R/4}}$, using both of \eqref{eq:self_energy_limits} and \eqref{eq:no_long_range}, and using the symmetry of $I_\mu$ we have
	\begin{equation*}
	\begin{split}
		\beta_0 + \circ(1)
		& \leq \int_{B_{2r_k}}f_k\\
		& = \beta_0 + \int_{B_{2r_k}}V_kI_\mu[e^{\lambda u_k}\chi_{B_{R/4}\setminus B_{2r_k}}]e^{\lambda u_k} + \circ(1)\\
		& \leq \beta_0 + \|V_k\|_{L^\infty(B_{2r_k})}\int_{B_R\setminus B_{2r_k}}I_\mu[e^{\lambda u_k}\chi_\Omega]e^{\lambda u_k} + \circ(1).
	\end{split}
	\end{equation*}
	Combining this estimate with the estimate
	\begin{equation*}
		\beta_0 +\circ(1)
		\leq \int_{B_R}f_k
		\leq \int_{B_{2r_k}}f_k + \|V_k\|_{L^\infty(B_R)}\int_{B_R\setminus B_{2r_k}}I_\mu[e^{\lambda u_k}\chi_\Omega]e^{\lambda u_k},
	\end{equation*}
	we find that to complete the proof of the lemma in the case $N = 1$, it suffices to show that 
	\begin{equation}
	\label{eq:m1_goal}
		\lim_{k\to\infty}\int_{B_R\setminus B_{2r_k}}I_\mu[e^{\lambda u_k}\chi_\Omega]e^{\lambda u_k}= 0. 
	\end{equation}
	Assumption \eqref{eq:uniformly_down_away_from_origin} guarantees that \eqref{eq:m1_goal} holds whenever (along a subsequence) $r_k\to r> 0$. Indeed, in this case for $k$ large
	\begin{equation*}
	\begin{split}
		\int_{B_R\setminus B_{2r_k}}f_k
		& \leq \int_{B_R\setminus B_r}V_kI_\mu[e^{\lambda u_k}\chi_\Omega]e^{\lambda u_k}\\
		& \leq \|V_k\|_{L^\infty(B_R)}\|I_\mu[e^{\lambda u_k}\chi_\Omega]\|_{L^{4/\mu}(B_R)}\|e^{u_k}\|_{L^1(B_R\setminus B_r)}^\lambda\\
		& \leq \mc Hc_0^\lambda\|V_k\|_{L^\infty(B_R)} \|e^{u_k}\|^\lambda_{C^0(B_R\setminus B_r)}\\
		& = \circ(1). 
	\end{split}
	\end{equation*}
	Assume for the remainder of \ref{item:assume_N1} that $r_k\to 0$. Lemma \ref{lemma:spherical_harnack} guarantees the existence of a universal constant $\beta\in (0, 1)$ and a positive constant $\bar C$ 
	such that
	\begin{equation*}
		\sup_{\bdy B_r}u_k
		\leq \bar C + \beta\inf_{\bdy B_r} u_k -2(1 - \beta)\log r
		\qquad \text{ for all }r\in [2r_k, \frac R2]. 
	\end{equation*}
	Moreover, for each $C_1>1$ (e.g., $C_1 = 2$), Lemma \ref{lemma:rescaled_sup+inf_inequality} guarantees the existence of a positive constant $C_2$ 
	such that  that
	\begin{equation*}
		u_k(0) + C_1\inf_{B_r}u_k + 2(1 + C_1)\log r\leq C_2
	\end{equation*}
	for all $r\in (0, R)$. Fixing $\beta$, $\bar C$, $C_1$ and $C_2$ as such, since $u_k$ is superharmonic 
	combining the previous two estimates gives
	\begin{equation}
	\label{eq:combine_inequalities}
		u_k(x)
		\leq C - \frac{\beta}{C_1}u_k(0) - 2\left(1+ \frac\beta{C_1}\right)\log|x| \qquad \text{ for all } x\in \overline B_{R/2}\setminus B_{2r_k},  
	\end{equation}
	where $C = C(\bar C, \beta, C_1)$. 
	Using this estimate and in view of \eqref{eq:rk_rate} we have
	\begin{equation*}
	\begin{split}
		\int_{B_{R/2}\setminus B_{2r_k}}e^{u_k}
		& \leq C\delta_k^{\frac{2\beta}{C_1}}\int_{\bb R^2\setminus B_{2r_k}}|x|^{-2(1 + \frac \beta{C_1})}\; \d x\\
		& \leq C\left(\frac{\delta_k}{r_k}\right)^{\frac{2\beta}{C_1}}\\
		& = \circ(1).  
	\end{split}
	\end{equation*}
	Therefore, since the integrability assumption in \eqref{eq:BR_PDE_and_energy} and the HLS inequality guarantee that $(I_\mu[e^{\lambda u_k}\chi_\Omega])_{k = 1}^\infty$ is bounded in $L^{4/\mu}(B_R)$, H\"older's inequality gives
	\begin{equation}
	\label{eq:base_case_harnack_region}
	\begin{split}
		\int_{B_{R/2} \setminus B_{2r_k}}& I_\mu[e^{\lambda u_k}\chi_\Omega]e^{\lambda u_k}\\
		& \leq \|I_\mu[e^{\lambda u_k}\chi_\Omega]\|_{L^{4/\mu}(B_R)}\left(\int_{B_{R/2}\setminus B_{2r_k}}e^{u_k}\right)^\lambda\\
		& 
		= \circ(1). 
	\end{split}
	\end{equation}
	Independently, assumption \eqref{eq:uniformly_down_away_from_origin} and the fact that $(I_\mu[e^{\lambda u_k}\chi_\Omega])_{k = 1}^\infty$ is bounded in $L^{4/\mu}(B_R)$ gives
	\begin{equation*}
		\int_{B_R\setminus B_{R/4}}I_\mu[e^{\lambda u_k}\chi_\Omega]e^{\lambda u_k}
		\leq \|I_\mu[e^{\lambda u_k}\chi_\Omega]\|_{L^{4/\mu}(B_R)}\|e^{u_k}\|_{L^1(B_R\setminus B_{R/4})}^\lambda
		= \circ(1)  
	\end{equation*}
	which, when combined with \eqref{eq:base_case_harnack_region} establishes \eqref{eq:m1_goal}. 
	\item \label{item:inductive_step} Suppose $N\geq 2$ and the lemma holds for $1, \ldots, N -1$. We proceed to show that it holds for $N$.  By relabeling the indices and passing to a suitable subsequence we may assume that 
	\begin{equation*}
		d_k := \min\{|x_k^{(j)} - x_k^{(i)}|: i, j\in \{0, \ldots, N - 1\}, \text{ and }i\neq j\}
	\end{equation*}
	satisfies $d_k = |x_k^{(0)} - x_k^{(1)}|$. As in \ref{item:assume_N1} we continue to assume that $x_k^{(0)} = 0$. We separately consider the case where the distances between $x_k^{(i)}$ and $x_k^{(j)}$ are all comparable and the case where these distances are not comparable. 
	\begin{enumerate}[label = {\bf Case \arabic*.}, ref = {Case \arabic*}]
		\item \label{case:A_exists} Assume there is $M\geq 1$ for which 
		\begin{equation}
		\label{eq:comparable_centers}
			d_k \leq |x_k^{(i)} - x_k^{(j)}|\leq Md_k 
		\end{equation}
		for all distinct indices $i, j\in \{0, \ldots, N - 1\}$ and all $k$. We will apply the lemma in the case $N = 1$ 
		to $u_k$ on $B_R$ with centers $x_k'^{(0)}= x_k^{(0)}= 0$, with scales $\delta_k'^{(0)} = \delta_k^{(0)}$, with radii $r_k'^{(0)} = 2Md_k$, and with mass $\beta_0' = \sum_{j = 0}^{N- 1}\beta_j$ to conclude that \eqref{eq:full_ BR_limit_neck_argument} holds. The remainder of the proof in \ref{case:A_exists} is devoted to verifying that the hypotheses of the lemma hold in this setting. To see that \eqref{eq:rk_rate} holds (with $r_k'^{(0)}$ and $\delta_k'^{(0)}$ in place of $r_k^{(j)}$ and $\delta_k^{(j)}$ respectively), note that assumption \eqref{eq:the_balls_disjoint} guarantees that $r_k^{(0)}\leq d_k$ so the limit $\frac{r_k'^{(0)}}{\delta_k'^{(0)}}\to\infty$ follows from \eqref{eq:rk_rate} with $j = 0$. 
		Item \ref{item:the_balls_disjoint} needs no verification when $N = 1$. To show that item \ref{item:termination_condition_satisfied} holds with $N =1$ and with $(x_k'^{(0)}, r_k'^{(0)})$ in place of $(x_k^{(0)}, r_k^{(0)})$, observe that if $x\in \overline B_R\setminus B_{2Md_k}$ and $j\in \{0, \ldots, N - 1\}$ then we have $|x - x_k^{(j)}|
		\geq |x|/2$ and therefore
		\begin{equation*}
			u_k(x) + 2\log|x|
			\leq u_k(x) + 2\log\min_{0\leq j\leq N - 1}|x - x_k^{(j)}| + 2\log 2
			\leq C,  
		\end{equation*}
		where the uniform upper bound follows from the containment $\bigcup_{j = 0}^{N - 1}B(x_k^{(j)}, r_k^{(j)})\subset B_{2Md_k}$ and assumption \eqref{eq:termination_condition_satisfied}.	
		To apply the Lemma with $N = 1$ it remains to verify that subitem \ref{item:beta_masses}\ref{subitem:self_energy_limits} holds with $(x_k'^{(0)}, r_k'^{(0)})$ in place of $(x_k^{(0)}, r_k^{(0)})$. 
		Equivalently, we must verify both
		\begin{equation}
		\label{eq:also_verify_for_N1_lemma}
			\int_{B_{32Md_k}}V_kI_\mu[e^{\lambda u_k}\chi_{\Omega\setminus B_{R/4}}]e^{\lambda u_k} = \circ(1) 
		\end{equation}
		and
		\begin{equation}
		\label{eq:verify_for_N1_lemma}
		\begin{split}
			\lim_k\int_{B_{32Md_k}} & V_kI_\mu[e^{\lambda u_k}\chi_{B_{32Md_k}}]e^{\lambda u_k}\\
			= & \; \lim_k\int_{B_{2Md_k}}V_kI_\mu[e^{\lambda u_k}\chi_{B_{2Md_k}}]e^{\lambda u_k}\\
			= &\; \sum_{j = 0}^{N - 1}\beta_j.
		\end{split}
		\end{equation} 
		To verify \eqref{eq:also_verify_for_N1_lemma}, first note that when $k$ is sufficiently large, for any $x\in B_{32Md_k}$ and any $z\in \Omega\setminus B_{R/4}$ we have $|x - z|
		\geq \frac R8$ and thus for any such $x$ an application of H\"older's inequality yields the following pointwise estimate:
		\begin{equation}
		\label{eq:far_away_pointwise}
		\begin{split}
			I_\mu[e^{\lambda u_k}\chi_{\Omega\setminus B_{R/4}}](x)
			& 
			\leq \|e^{u_k}\|_{L^1(\Omega)}^\lambda\left(\int_{\bb R^2\setminus B(x, R/8)}|x - z|^{-4}\; \d z\right)^{1 - \lambda}\\
			& 
			\leq C(c_0, \mu, R). 
		\end{split}
		\end{equation}
		Using this estimate together with H\"older's inequality gives
		\begin{equation}
		\label{eq:far_away_integral}
		\begin{split}
			\int_{B_{32Md_k}}& V_kI_\mu[e^{\lambda u_k}\chi_{\Omega\setminus B_{R/4}}]e^{\lambda u_k}\\
			& \leq \|V_k\|_{L^\infty(B_R)}\|e^{u_k}\|_{L^1(\Omega)}^\lambda \|I_\mu[e^{\lambda u_k}\chi_{\Omega\setminus B_{R/4}}]\|_{L^{4/\mu}(B_{32Md_k})}\\
			& 
			\leq C d_k^{\mu/2} = \circ(1), 
		\end{split}
		\end{equation}
		which is \eqref{eq:also_verify_for_N1_lemma}. To verify \eqref{eq:verify_for_N1_lemma} we first note that assumptions \eqref{eq:the_balls_disjoint}, \eqref{eq:exists_vanishing_relative_radius} and \eqref{eq:self_energy_limits} guarantee that for every pair of distinct indices $i, j\in \{0, \ldots,N -1\}$ we have
		\begin{equation}
		\label{eq:relative_radius_limsup}
			\limsup_k\frac{r_k^{(j)}}{|x_k^{(i)} - x_k^{(j)}|}\leq \frac 1{16}. 
		\end{equation}
		Indeed, \eqref{eq:comparable_centers} guarantees that for any such $i$ and $j$ we have either
		\begin{equation}
		\label{eq:one_relative_radius_vanishes}
			\frac{r_k^{(j)}}{|x_k^{(i)} - x_k^{(j)}|} = \circ(1)
			\quad \text{ or }\quad
			\frac{r_k^{(i)}}{|x_k^{(i)} - x_k^{(j)}|} = \circ(1). 
		\end{equation} 
		If the first equality in \eqref{eq:one_relative_radius_vanishes} holds then \eqref{eq:relative_radius_limsup} holds. If the second equality in \eqref{eq:one_relative_radius_vanishes} holds and if there is a subsequence of $k$ along which $\lim_k\frac{r_k^{(j)}}{|x_k^{(i)} - x_k^{(j)}|}> \frac 1{16}$ then we have $B(x_k^{(i)}, r_k^{(i)})\subset B(x_k^{(j)}, 16 r_k^{(j)})\setminus B(x_k^{(j)}, r_k^{(j)})$ whenever $k$ is sufficiently large. 
		Using this containment and in view of \eqref{eq:self_energy_limits} we have
		\begin{equation*}
		\begin{split}
			\beta_i + \circ(1)
			= & \; \int_{B(x_k^{(i)}, r_k^{(i)})}V_kI_\mu[e^{\lambda u_k}\chi_{B(x_k^{(i)}, r_k^{(i)})}]e^{\lambda u_k}\\
			\leq & \; \int_{B(x_k^{(j)}, 16 r_k^{(j)})\setminus B(x_k^{(j)}, r_k^{(j)})}V_kI_\mu[e^{\lambda u_k}\chi_{B(x_k^{(j)}, 16 r_k^{(j)})\setminus B(x_k^{(j)}, r_k^{(j)})}]e^{\lambda u_k}\\
			\leq &\; \int_{B(x_k^{(j)}, 16 r_k^{(j)})}V_kI_\mu[e^{\lambda u_k}\chi_{B(x_k^{(j)}, 16 r_k^{(j)})}]e^{\lambda u_k}\\
			& - \int_{B(x_k^{(j)}, r_k^{(j)})}V_kI_\mu[e^{\lambda u_k}\chi_{B(x_k^{(j)}, r_k^{(j)})}]e^{\lambda u_k}\\
			= & \; \circ(1),
		\end{split}
		\end{equation*}
		which contradicts the positivity of $\beta_i$. Estimate \eqref{eq:relative_radius_limsup} is established. Define 
		\begin{equation*}
			U_k = \bigcup_{j = 0}^{N -1} B(x_k^{(j)}, 2r_k^{(j)})
			\quad \text{ and }\quad
			E_k = B_{32Md_k}\setminus U_k. 
		\end{equation*}	
		Combining assumptions \eqref{eq:self_energy_limits} and \eqref{eq:assumption_no_close_range} shows that 
		\begin{equation*}
			\int_{U_k}V_kI_\mu[e^{\lambda u_k}\chi_{U_k}]e^{\lambda u_k}
			= \sum_{j = 0}^{N - 1}\beta_j + \circ(1). 
		\end{equation*}
		Combining this equality with the estimate
		\begin{equation*}
		\begin{split}
			0
			& \leq \int_{B_{32Md_k}}V_kI_\mu[e^{\lambda u_k}\chi_{B_{32Md_k}}]e^{\lambda u_k} - \int_{U_k}V_kI_\mu[e^{\lambda u_k}\chi_{U_k}]e^{\lambda u_k}\\
			& \leq C\left(\int_{U_k}I_\mu[e^{\lambda u_k}\chi_{E_k}]e^{\lambda u_k} + \int_{E_k}I_\mu[e^{\lambda u_k}\chi_{E_k}]e^{\lambda u_k}\right)\\
			& \leq C\int_{E_k}I_\mu[e^{\lambda u_k}\chi_\Omega]e^{\lambda u_k}, 
		\end{split}
		\end{equation*}
		we find that to establish the equality
		\begin{equation}
		\label{eq:bigg_balls_limit}
			\lim_k\int_{B_{32Md_k}}V_kI_\mu[e^{\lambda u_k}\chi_{B_{32Md_k}}]e^{\lambda u_k}
			= \sum_{j = 0}^{N -1}\beta_j, 
		\end{equation}
		it suffices to show that
		\begin{equation}
		\label{eq:light_purple}
			\int_{E_k}I_\mu[e^{\lambda u_k}\chi_\Omega]e^{\lambda u_k} = \circ(1). 
		\end{equation}
		We will do so via a rescaling argument. Define
		\begin{equation*}
		\tilde u_k(y) = u_k(d_ky) + 2\log d_k
		\qquad \text{ for }y\in d_k^{-1}\Omega=:\tilde \Omega_k
	\end{equation*}
	and for each $j\in \{0, \ldots, N- 1\}$ set $\tilde x_k^{(j)} = d_k^{-1}x_k^{(j)}$, $\tilde r_k = d_k^{-1}r_k^{(j)}$ and 
	\begin{equation*}
		\tilde \delta_k^{(j)} 
		= e^{-\tilde u_k(\tilde x_k^{(j)})/2}
		= \frac{\delta_k^{(j)}}{d_k}. 
	\end{equation*}
	Combining the inequality $r_k^{(j)}\leq Md_k$ 
	with assumption \eqref{eq:rk_rate} gives $\tilde u_k(\tilde x_k^{(j)})\to\infty$. 
	Morevoer, since $u_k$ satisfies \eqref{eq:BR_PDE_and_energy}, $\tilde u_k$ satisfies
	\begin{equation*}
	\begin{cases}
		-\lap \tilde u_k = \tilde V_kI_\mu[e^{\lambda \tilde u_k}\chi_{\tilde \Omega_k}]e^{\lambda \tilde u_k} & \text{ in } B_{d_k^{-1}R}\\
		\|e^{\tilde u_k}\|_{L^1(\tilde \Omega_k)}\leq c_0,
	\end{cases}
	\end{equation*}
	where $\tilde V_k(y)= V_k(d_k^{-1}y)$. Assumption \eqref{eq:comparable_centers} together with the equality $\tilde x_k^{(0)} = 0$ guarantees that $(\tilde x_k^{(j)})_{k = 1}^\infty\subset \overline B_{4M}$ for all $j\in \{0, \ldots, N - 1\}$, so after passing to a subsequence we assume the existence of $\{\tilde x^{(j)}\}\subset \overline B_{4M}$ satisfying both $\tilde x^{(0)} = 0$ and $|\tilde x^{(i)} - \tilde x^{(j)}|\geq 1$ whenever $i\neq j$ and such that $\tilde x_k^{(j)}\to \tilde x^{(j)}$ for all $j\in \{0, \ldots, N - 1\}$. An application of Theorem \ref{theorem:BM_alternative} guarantees that $\tilde u_k\to-\infty$ locally uniformly on $\bb R^2\setminus \bigcup_{j = 0}^{N - 1}\{\tilde x^{(j)}\}$. For each $j\in\{0, \ldots, N - 1\}$, let $m_j\in \{0, \ldots, N - 1\}\setminus \{j\}$ be any index for which 
	\begin{equation*}
		|\tilde x^{(j)} - \tilde x^{(m_j)}|
		= \min\{|\tilde x^{(j)} - \tilde x^{(i)}| :i \in \{0, \ldots, N - 1\}\setminus \{j\}\}. 
	\end{equation*}
	Define the compact set
	\begin{equation*}
		K = \overline B_{4M}\setminus \bigcup_{j = 0}^{N - 1}B\left(\tilde x^{(j)}, \frac{|\tilde x^{(j)} - \tilde x^{(m_j)}|}{8}\right)
	\end{equation*}
	and for each $j\in \{0, \ldots, N - 1\}$ define the annulus
	\begin{equation*}
		\tilde A_k^{(j)}
		= B\left(\tilde x_k^{(j)}, \frac{5|\tilde x_k^{(j)} - \tilde x_k^{(m_j)}|}{32}\right)\setminus B(\tilde x_k^{(j)}, 2\tilde r_k^{(j)}).
	\end{equation*}
	Inequality \eqref{eq:relative_radius_limsup} guarantees that each $\tilde A_k^{(j)}$ has nonempty interior whenever $k$ is sufficiently large. Moreover, setting $\tilde A_k = \bigcup_{j=0}^{N - 1}\tilde A_k^{(j)}$ for $k$ sufficiently large we have $\tilde E_k\subset  \tilde A_k \cup K$, where 
	\begin{equation*}
		\tilde E_k
		= B_{4M}\setminus \bigcup_{j = 0}^{N-1}B(\tilde x_k^{(j)}, 2\tilde r_k^{(j)}) 
	\end{equation*}
	is the image of $E_k$ under the rescaling $x = d_k y$. Using the change of variable $x = d_ky$ we have
	\begin{equation}
	\label{eq:purple}
	\begin{split}
		\int_{E_k}I_\mu[e^{\lambda u_k}\chi_\Omega]e^{\lambda u_k}
		& = \int_{\tilde E_k}I_\mu[e^{\lambda \tilde u_k}\chi_{\tilde \Omega_k}]e^{\lambda \tilde u_k}\\
		& \leq \int_{\tilde A_k}I_\mu[e^{\lambda \tilde u_k}\chi_{\tilde \Omega_k}]e^{\lambda \tilde u_k} + \int_{K}I_\mu[e^{\lambda \tilde u_k}\chi_{\tilde \Omega_k}]e^{\lambda \tilde u_k}\\
		& = \int_{\tilde A_k}I_\mu[e^{\lambda \tilde u_k}\chi_{\tilde \Omega_k}]e^{\lambda \tilde u_k} + \circ(1),
	\end{split}
	\end{equation}
	where the final equality holds by H\"older's inequality, the fact that $(I_\mu[e^{\lambda \tilde u_k}\chi_{\tilde \Omega_k}])_{k = 1}^\infty$ is bounded in $L^{4/\mu}(K)$, and the fact that $\tilde u_k\to-\infty$ uniformly on $K$. To estimate the remaining integral on the right-most side of \eqref{eq:purple}, first observe that if $k$ is sufficiently large then for every $i\neq j$ we have
		\begin{equation}
		\label{eq:B38_nonintersecting}
			B\left(\tilde x_k^{(i)}, \frac{3|\tilde x_k^{(i)} - \tilde x_k^{(m_i)}|}{8}\right)
			\cap B\left(\tilde x_k^{(j)}, \frac{3|\tilde x_k^{(j)} - \tilde x_k^{(m_j)}|}{8}\right)
			= \emptyset.
		\end{equation}
		Therefore, for any $y\in B\left(\tilde x_k^{(j)}, 5|\tilde x_k^{(j)} - \tilde x_k^{(m_j)}|/(16)\right)\setminus B(\tilde x_k^{(j)}, \tilde r_k^{(j)})$  
		assumption \eqref{eq:termination_condition_satisfied} gives
		\begin{equation*}
		\begin{split}
			\tilde u_k(y) & +2\log|y - \tilde x_k^{(j)}|\\
			& = \tilde u_k(y) +2\log\min_{0\leq i\leq N - 1}|y - \tilde x_k^{(i)}|\\
			& 
			\leq C. 
		\end{split}
		\end{equation*}
		Lemma \ref{lemma:spherical_harnack} guarantees the existence of a universal constant $\beta\in (0, 1)$ and a constant $\overline C>0$ such that 
		\begin{equation}
		\label{eq:jth_spherical_harnack}
			\sup_{\bdy B_r} \tilde u_k(\tilde x_k^{(j)} + \cdot)
			\leq \overline C + \beta\inf_{\bdy B_r}\tilde u_k(\tilde x_k^{(j)} + \cdot) - 2(1 - \beta)\log r
		\end{equation}
		for all $r\in [2\tilde r_k^{(j)}, 5|\tilde x_k^{(j)} - \tilde x_k^{(m_j)}|/(32)]$. Moreover, for any $C_1>1$ (e.g., $C_1 = 2$) Lemma \ref{lemma:rescaled_sup+inf_inequality} guarantees the existence of a $k$-independent constant $C_2$ such that 
		\begin{equation}
		\label{eq:jth_rescaled_supinf}
			\tilde u_k(\tilde x_k^{(j)}) + C_1\inf_{B_r}\tilde u_k(\tilde x_k^{(j)} + \cdot) + 2(1 + C_1)\log r \leq C_2
		\end{equation}
		whenever $0< r< 3|\tilde x^{(j)} - \tilde x^{(m_j)}|/8$. 
		Combining estimates \eqref{eq:jth_spherical_harnack} and \eqref{eq:jth_rescaled_supinf} and since $\tilde u_k$ is superharmonic we find that there is a $k$-independent constant $C$ for which 
		\begin{equation*}
			e^{\tilde u_k(y)}
			\leq C(\tilde \delta_k^{(j)})^{\frac{2\beta}{C_1}}|y - \tilde x_k^{(j)}|^{-2(1 + \frac \beta{C_1})}
		\end{equation*}
		for all $y\in \tilde A_k^{(j)}$ and therefore 
		\begin{equation*}
		\begin{split}
			\int_{\tilde A_k^{(j)}}e^{u_k}
			& \leq \left(\tilde \delta_k^{(j)}\right)^{\frac{2\beta}{C_1}}\int_{\bb R^2\setminus B(\tilde x_k^{(j)}, 2\tilde r_k^{(j)})}|y - \tilde x_k^{(j)}|^{-2(1 + \frac \beta{C_1})}\; \d y\\
			& \leq C\left(\frac{\tilde \delta_k^{(j)}}{\tilde r_k^{(j)}}\right)^{\frac{2\beta}{C_1}}\\
			& = \circ(1). 
		\end{split}
		\end{equation*}
		Finally, since $(\|I_\mu[e^{\lambda \tilde u_k}\chi_{\tilde \Omega_k}]\|_{L^{4/\mu}(\tilde \Omega_k)})_{k= 1}^\infty$ is bounded in $\bb R$ we obtain 
		\begin{equation*}
			\int_{\tilde A_k}I_\mu[e^{\lambda \tilde u_k}\chi_{\tilde \Omega_k}]e^{\lambda u_k}
			\leq \|I_\mu[e^{\lambda \tilde u_k}\chi_{\tilde \Omega_k}]\|_{L^{4/\mu}(\tilde \Omega_k)}\sum_{j = 0}^{N - 1}\left(\int_{\tilde A_k^{(j)}}e^{\tilde u_k}\right)^\lambda\\
			= \circ(1).
		\end{equation*}
		Bringing this estimate back to \eqref{eq:purple} establishes \eqref{eq:light_purple} thereby completing the verification of \eqref{eq:bigg_balls_limit}. The verification that 
		\begin{equation*}
			\lim_k\int_{B_{2Md_k}}V_kI_\mu[e^{\lambda u_k}\chi_{B_{2Md_k}}]e^{\lambda u_k}
			= \sum_{j = 0}^{N -1}\beta_j
		\end{equation*}
		follows by a similar argument so we omit the details. Equalities \eqref{eq:verify_for_N1_lemma} are established. 		
	\item \label{case:is_no_M} Assume there is no $M\geq 1$ for which \eqref{eq:comparable_centers} holds for all distinct indices $i, j\in\{0, \ldots, N  - 1\}$ and all $k$. There is a proper subset $J\subset \{0, \ldots,N - 1\}$ containing $\{0, 1\}$ and there is a constant $M\geq 1$  for which 
	\begin{equation*}
		|x_k^{(j)}|\leq M d_k
		\qquad \text{ for all }j\in J \text{ and all }k
	\end{equation*}
	and
	\begin{equation}
	\label{eq:non_comparable_centers}
		\lim_k\frac{|x_k^{(j)}|}{d_k} = +\infty
		\qquad \text{ for }j\in \{0, \ldots, N - 1\}\setminus J. 
	\end{equation}
	We assume without losing generality that $J = \{0, 1, \ldots, n - 1\}$ for some $n\in \{2, \ldots, N - 1\}$. We will apply the lemma with $N - n + 1$ (induction hypothesis) to the centers $x_k'^{(0)} = x_k^{(0)} = 0$, $\{x_k^{(j)}\}_{j = n}^{N -1}$, the radii $r_k'^{(0)} = 2Md_k$, $\{r_k^{(j)}\}_{j = n}^{N - 1}$ and the masses $\beta_0' = \sum_{j = 0}^{n - 1}\beta_j$, $\{\beta_j\}_{j = n}^{N - 1}$ to obtain 
	\begin{equation*}
		\int_{B_R}f_k\to \beta_0' + \sum_{j = n}^{N- 1}\beta_j
		= \sum_{j = 0}^{N - 1}\beta_j. 
	\end{equation*}
	In the remainder of the proof of \ref{case:is_no_M} we verify that the hypotheses of the lemma hold with $N - n + 1$ and with these centers, radii and masses. To verify that item \ref{item:rk_rate} holds it suffices to show that $\frac{r_k'^{(0)}}{\delta_k'^{(0)}}= \frac{2Md_k}{\delta_k^{(0)}}\to \infty$. The argument for doing so is as in \ref{case:A_exists}. To verify \eqref{eq:the_balls_disjoint}, it suffices to show that
	\begin{equation}
	\label{eq:induction_balls_disjoint}
		B(0, 2Md_k)\cap B(x_k^{(j)}, r_k^{(j)}) = \emptyset 
		\qquad \text{ for }j = n, \ldots, N - 1. 
	\end{equation}
	We do so by way of contradiction. If $j\in \{n, \ldots, N - 1\}$ is an index for which \eqref{eq:induction_balls_disjoint} fails then $r_k^{(j)} + 2Md_k \geq |x_k^{(j)}|$ and therefore \eqref{eq:non_comparable_centers} implies that $\frac{ r_k^{(j)}}{d_k}\to\infty$. Fixing any such $j$ we have 
	\begin{equation*}
		1\leq \frac{r_k^{(j)}}{|x_k^{(j)}|} + \frac{2Md_k}{|x_k^{(j)}|}
		= \frac{r_k^{(j)}}{|x_k^{(j)}|} + \circ(1)
		\leq 1 + \circ(1),  
	\end{equation*}
	where the final inequality holds since assumption \eqref{eq:the_balls_disjoint} applied with $i = 0$ guarantees that $r_k^{(j)} + r_k^{(0)}\leq |x_k^{(j)}|$. This shows that $B(0, Md_k)\subset B(x_k^{(j)}, 2r_k^{(j)})$ whenever $k$ is sufficiently large. 
	In particular, this containment combined with assumption \eqref{eq:the_balls_disjoint} (applied with $i = 0$) guarantees that $B(0, r_k^{(0)})\subset B(x_k^{(j)}, 2r_k^{(j)})\setminus B(x_k^{(j)}, r_k^{(j)})$. Therefore, the second equality in assumption \eqref{eq:self_energy_limits} applied with $j = 0$ and the first equality in \eqref{eq:self_energy_limits} applied with $j= j$ gives
	\begin{equation*}
	\begin{split}
		\beta_0 + \circ(1)
		&= \int_{B(0, r_k^{(0)})}V_kI_\mu[e^{\lambda u_k}\chi_{B(0, r_k^{(0)})}]e^{\lambda u_k}\\
		& \leq \int_{B(x_k^{(j)}, 2r_k^{(j)})\setminus B(x_k^{(j)}, r_k^{(j)})} V_k I_\mu[e^{\lambda u_k}\chi_{B(x_k^{(j)}, 2r_k^{(j)})\setminus B(x_k^{(j)}, r_k^{(j)})}]e^{\lambda u_k}\\
		& \leq \int_{B(x_k^{(j)}, 2r_k^{(j)})} V_k I_\mu[e^{\lambda u_k}\chi_{B(x_k^{(j)}, r_k^{(j)})}]e^{\lambda u_k} \\
		& \quad  - \int_{B(x_k^{(j)}, r_k^{(j)})} V_k I_\mu[e^{\lambda u_k}\chi_{B(x_k^{(j)}, 2r_k^{(j)})}]e^{\lambda u_k}\\
		& = \circ(1), 
	\end{split}
	\end{equation*}
	which contradicts the assumption $\beta_0>0$. This completes the verification of \eqref{eq:the_balls_disjoint}. To verify \eqref{eq:exists_vanishing_relative_radius} we only need to show that $\frac{2Md_k}{|x_k^{(j)}|}\to 0$ whenever $j\in \{n, \ldots, N - 1\}$. This follows from \eqref{eq:non_comparable_centers}. To verify \eqref{eq:termination_condition_satisfied} observe that if $x\in B_R\setminus B_{2Md_k}$ then for any $j\in \{0, \ldots, n - 1\}$ we have $|x - x_k^{(j)}| \geq
	|x|/2$. Therefore, 
	for any $x\in \overline B_R\setminus \left(B_{2Md_k}\cup(\bigcup_{j = n}^{N - 1}B(x_k^{(j)}, r_k^{(j)}))\right)$, using assumption \eqref{eq:termination_condition_satisfied} we have 
	\begin{equation*}
	\begin{split}
		u_k(x) & + 2\log\min\{|x|, \min_{n\leq j\leq N -1}|x - x_k^{(j)}|\}\\
		& \leq u_k(x) + 2\log\min_{0\leq j\leq N -1}|x - x_k^{(j)}| + 2\log 2\\
		& \leq C. 
	\end{split}
	\end{equation*}
	To verify \eqref{eq:self_energy_limits} we only need to verify that
	\begin{equation*}
	\begin{split}
		\lim_k\int_{B_{32Md_k}}&V_kI_\mu[e^{\lambda u_k}\chi_{B_{32Md_k}}]e^{\lambda u_k}\\
		& = \lim_k\int_{B_{2Md_k}}V_kI_\mu[e^{\lambda u_k}\chi_{B_{2Md_k}}]e^{\lambda u_k}\\
		&= \beta_0' = \sum_{j = 0}^{n - 1}\beta_j.		
	\end{split}
	\end{equation*}
	This is accomplished by the argument presented in \ref{case:A_exists}, so the details are omitted. To verify \eqref{eq:no_long_range} it suffices to show that
	\begin{equation}
	\label{eq:induction_no_long_range}
		\int_{B_{32Md_k}} I_\mu[e^{\lambda u_k}\chi_{\Omega\setminus B_{R/4}}]e^{\lambda u_k} = \circ(1). 
	\end{equation}
	This follows from an argument similar to the one carried out in \eqref{eq:far_away_pointwise}, \eqref{eq:far_away_integral} so we only give a brief description. First using the fact that $|x - z|\geq R/8$ whenever $x\in B_{32Md_k}$, $z\in \Omega\setminus B_{R/4}$ and $k$ is large we show that for any such $x$, there holds $I_\mu[e^{\lambda u_k}\chi_{\Omega\setminus B_{R/4}}](x)\leq C(\mu, c_0, R)$. Using this pointwise estimate and H\"older's inequality we find that 
	\begin{equation*}
		\int_{B_{32Md_k}} I_\mu[e^{\lambda u_k}\chi_{\Omega\setminus B_{R/4}}]e^{\lambda u_k}
		\leq Cd_k^{\mu/2}
		= \circ(1). 
	\end{equation*}
	To verify \eqref{eq:assumption_no_close_range} it suffices to show that
	\begin{equation}
	\label{eq:induction_no_close_range}
		\int_{B_{30Md_k}}I_\mu[e^{\lambda u_k}\chi_{B(x_k^{(j)}, 2r_k^{(j)})}]e^{\lambda u_k}
		= \circ(1)
	\end{equation}
	for $j = n, \ldots, N - 1$. 
	Fix any such $j$ and observe that for $k$ large, for $x\in B_{30Md_k} $ and for $z\in B(x_k^{(j)}, 2r_k^{(j)})$, \eqref{eq:relative_radius_limsup} and \eqref{eq:non_comparable_centers} guarantees that $|z - x|\geq 7|x_k^{(j)}|/9$. 
	Therefore, for $x\in B_{30Md_k}$ we have the pointwise estimate
	\begin{equation*}
	\begin{split}
		I_\mu[e^{\lambda u_k}\chi_{B(x_k^{(j)}, 2r_k^{(j)})}](x)
		& 
		\leq \|e^{u_k}\|_{L^1(\Omega)}^\lambda\left(\int_{\bb R^2\setminus B(x, 7|x_k^{(j)}|/9)}|x- z|^{-4}\; \d z\right)^{1 - \lambda}\\
		& \leq C|x_k^{(j)}|^{-\mu/2}. 
	\end{split}
	\end{equation*}
	Using this estimate together with H\"older's inequality and assumption \eqref{eq:non_comparable_centers} we obtain 
	\begin{equation*}
	\begin{split}
		\int_{B_{30Md_k}}& I_\mu[e^{\lambda u_k}\chi_{B(x_k^{(j)}, 2r_k^{(j)})}]e^{\lambda u_k}\\
		& \leq \|e^{u_k}\|_{L^1(\Omega)}^\lambda \|I_\mu[e^{\lambda u_k}\chi_{B(x_k^{(j)}, 2r_k^{(j)})}]\|_{L^{4/\mu}(B(30Md_k))}\\
		& 
		\leq C\left(\frac{d_k}{|x_k^{(j)}|}\right)^{\mu/2}\\
		& = \circ(1). 
	\end{split}
	\end{equation*}
	\end{enumerate}
\end{enumerate}
\end{proof}
With Lemmata \ref{lemma:N_bubble_selection} and \ref{lemma:neck_energy_vanish} in hand we now give the short proof of Proposition \ref{prop:energy_quantization}. 
\begin{proof}[Proof of Proposition \ref{prop:energy_quantization}]
Lemma \ref{lemma:blow_up_minimal_energy} guarantees that $V(0)>0$ so in view of assumption \eqref{eq:uniformly_down_away_from_origin}, after decreasing $R$ if necessary, we assume the existence of $a>0$ for which $a\leq V_k(x)$ for all $x\in \overline B_R$ and all $k\in \bb N$. Let $(\rho_k)_{k= 1}^\infty\subset (0, \infty)$ satisfy $\rho_k\to\infty$ and let $N\in \bb N$ and $\{(x_k^{(j)})_{k = 1}^\infty: j = 0, \ldots, N - 1\}$ be any positive integer and any collection of sequences in $B_R$ respectively whose existence is guaranteed by Lemma \ref{lemma:N_bubble_selection}. The assertion of Proposition \ref{prop:energy_quantization} follows by applying Lemma \ref{lemma:neck_energy_vanish} with $\delta_k^{(j)} = e^{-u_k(x_k^{(j)})/2}$, $r_k^{(j)} = \rho_k\delta_k^{(j)}$ and $\beta_j = 8\pi$ for all $j$. 
\end{proof}
%
\appendix\section{}
\label{s:appendix}
The following classification of solutions to problem \eqref{eq:local_entire} is established in  \cite{ChenLi1991}.
\begin{oldtheorem}
\label{theorem:CL_classification}
Every solution to problem \eqref{eq:local_entire} is of the form $u(x) = U_0(\delta(x -x_0)) + 2\log\delta$ for some $(x_0, \delta)\in \bb R^2\times (0, \infty)$, where 
\begin{equation}
\label{eq:CL_bubble}
	U_0(x) = \log\frac{8}{(1 + |x|^2)^2}. 
\end{equation}
In particular, every solution $u$ to problem \eqref{eq:local_entire} satisfies $\|e^u\|_{L^1(\bb R^2)} = 8\pi$. 
\end{oldtheorem}

The following classification of solutions to problem \eqref{eq:entire_PDE_and_L1} is established in \cite{Gluck2025classification}.  
\begin{oldtheorem}
\label{oldtheorem:classification}
Let $\mu\in (0, 2)$ and let $\lambda$ be as in \eqref{eq:lambda}. If $u\in L^1_{\loc}(\bb R^2)$ is a distributional solution to \eqref{eq:entire_PDE_and_L1} then there is $(x_0, \delta)\in \bb R^2\times (0, \infty)$ for which 
\begin{equation}
\label{eq:rescaled_bubble} 
	u(x) = U(\delta(x- x_0)) + 2\log \delta,
\end{equation} 
where 
\begin{equation}
\label{eq:unscaled_bubble}
	U(x)
	= -2\log(1 + |x|^2) + \frac 2{4 - \mu}\log\left(\frac{4(2 -\mu)}\pi\right)
	\qquad\text{ for }x\in \bb R^2. 
\end{equation}
In particular, for any such $u$ both of the following equalities hold:
\begin{equation}
\label{eq:asserted_energies}
	\int_{\bb R^2}e^u = (4(2 - \mu))^{\frac{2}{4 -\mu}}\pi^{\frac{2 - \mu}{4 - \mu}}
	\qquad \text{ and }\qquad
	\int_{\bb R^2}I_\mu[e^{\lambda u}]e^{\lambda u} = 8\pi. 
\end{equation}
\end{oldtheorem}
%
The proof of the following lemma can be found in Lemma 5.4 of \cite{Gluck2020classification}. 
\begin{lemma}
\label{lemma:bubble_potential}
The equality 
\begin{equation*}
	\log\left(\frac 2{1 + |x|^2}\right)
	= \frac{2}{|\bb S^n|}\int_{\bb R^n}\log\left(\frac{\sqrt2}{|x- y|}\right)\left(\frac 2{1 + |y|^2}\right)^n\; \d y
\end{equation*}
holds for all $x\in \bb R^n$. 
\end{lemma}

\begin{proof}[Proof of Lemma \ref{lemma:sup+inf_driving_estimate}]
It suffices to prove the lemma in the case $B_\rho(x_0) = B_1 \subset \omega$ as the general case can be recovered by considering $v(x) = u(x_0 +\rho x)$ for $x\in \tilde\omega := \frac{\omega - x_0}{\rho}\supset B_1$. Accordingly, let us assume $B_\rho(x_0) = B_1 \subset \omega$ and let $G(x,y)$ denote the Dirichlet Green's function for $-\lap$ on $B_1$. Green's representation formula gives
\begin{equation*}
	u(x)
	\geq\int_{B_1}G(x,y) f(y)\; \d y + \inf_{\bdy B_1}u
\end{equation*}
for all $x\in B_1$. 
From this inequality and the explicit expression 
\begin{equation*}
	G(x,y)
	= -\frac 1{2\pi}\left(\log|x - y| - \log\abs{|x|(y - \frac{x}{|x|^2})}\right), 
\end{equation*}
we obtain 
\begin{equation}
\label{eq:use_similar_triangles}
\begin{split}
	u(x) - \inf_\omega u
	& \geq u(x) - \inf_{B_1}u\\
	& \geq -\frac 1{2\pi}\int_{B_1}\log\frac{|x- y|}{\abs{|x|y - \frac x{|x|}}}f(y)\; \d y\\
	& = -\frac 1{2\pi}\int_{B_1}\log|y| f(y)\l \; \d y.  
\end{split}
\end{equation}
Integration by parts in \eqref{eq:use_similar_triangles} and using the equality 
\begin{equation*}
	\lim_{s\to 0}\log s\int_{B_s}f(y)\; \d y
	= \lim_{s\to 0}|B_s|\log s\fint_{B_s}f(y)\; \d y
	= 0 
\end{equation*}
gives
\begin{equation*}
\begin{split}
	u(x) - \inf_\omega u
	& \geq -\frac 1{2\pi}\int_0^1\log s\int_{\bdy B_s}f(y)\; \d S_y\; \d s\\
	& 
	= -\frac 1{2\pi}\left[\left.\log s\int_{B_s}f(y)\; \d y\right|_0^1 - \int_0^1\int_{B_s}f(y)\; \d y\; \frac{\d s}s\right]\\
	& = \frac 1{2\pi}\int_0^1\int_{B_s}f(y)\; \d y\; \frac{\d s}s, 
\end{split}
\end{equation*}
for all $x\in B_1$. Choosing $x = 0$ yields 
\begin{equation*}
	u(0) - \inf_\omega u
	\geq \frac 1{2\pi}\int_0^1\int_{B_s}f(x)\; \d x \; \frac{\d s}s. 
\end{equation*}
Upon rescaling and since $f\geq 0$, for any $r\in (0, \rho)$ we obtain 
\begin{equation*}
\begin{split}
	u(x_0) - \inf_\omega u
	& \geq \frac 1{2\pi}\int_0^\rho\int_{B_s(x_0)}f(x)\; \d x \; \frac{\d s}s\\
	& \geq \frac 1{2\pi}\int_r^\rho\int_{B_s(x_0)}f(x)\; \d x \; \frac{\d s}s\\
	& \geq \frac 1{2\pi}\int_{B_r(x_0)}f(x)\; \d x \int_r^\rho \; \frac{\d s}s\\
	& = \frac 1{2\pi}\int_{B_r(x_0)}f(x)\; \d x \; \log\frac\rho r. 
\end{split} 
\end{equation*}
\end{proof}
\begin{lemma}
\label{lemma:selection_process}
Let $(\tilde x, \rho)\in \bb R^n\times (0, \infty)$ and let  $\varphi\in C^0(\overline B_\rho(\tilde x))$ be a positive function. For any $a>0$ there exists $x\in B(\tilde x, \rho)$ and $r = r(x)>0$ for which both 
\begin{equation*}
	\varphi(x)\geq \left(\frac \rho{2r}\right)^a \varphi(\tilde x)
\end{equation*}
and 
\begin{equation*}
	\varphi(x)\geq \left(\frac 12\right)^a\max_{\overline B_r(x)}\varphi(y). 
\end{equation*}
\end{lemma}
\begin{proof}
Let $x\in B_\rho(\tilde x)$ be a maximizer of the function $\psi:\overline B_\rho(\tilde x)\to [0, \infty)$ defined by 
\begin{equation*}
	\psi(y)
	= (\rho - |y - \tilde x|)^a\varphi(y). 
\end{equation*}
and set $r = r(x) = (\rho - |x - \tilde x|)/2$. The inequality $\psi(x)\geq \psi(\tilde x)$ gives
\begin{equation*}
	\varphi(x)
	\geq \left(\frac\rho{\rho - |x - \tilde x|}\right)^a\varphi(\tilde x)
	= \left(\frac{\rho}{2r}\right)^a\varphi(\tilde x), 
\end{equation*}
which is the first of the asserted estimates. To show the second of the asserted inequalities, observe that for any $y\in \overline B_r(x)$ we have $|y - \tilde x|\geq |x - \tilde x| - r$ and therefore, 
\begin{equation*}
\begin{split}
	\varphi(x)
	& \geq \left(\frac{\rho - |y - \tilde x|}{\rho - |x - \tilde x|}\right)^a \varphi(y)\\
	& \geq \left(\frac r{\rho - |x - \tilde x|}\right)^a \varphi(y)\\
	& = 2^{-a}\varphi(y). 
\end{split}
\end{equation*}
The second of the asserted inequalities follows. 
\end{proof}
%
%
\bibliographystyle{alpha}

\end{document}